\documentclass{amsart}
\usepackage{amssymb}

\newtheorem{thm}{Theorem}[section]
\newtheorem{prop}[thm]{Proposition}
\newtheorem{cor}[thm]{Corollary}

\newtheorem{rema}[thm]{Remark}
\newtheorem{defn}[thm]{Definition}
\newtheorem{conj}[thm]{Conjecture}

\theoremstyle{remark}

\numberwithin{equation}{section}

\newcommand{\asqrt}[0]{a_{\scriptscriptstyle{\square}}}
\newcommand{\bsqrt}[0]{b_{\scriptscriptstyle{\square}}}
\newcommand{\nor}[0]{\Upsilon}
\newcommand{\sou}[0]{\Delta}

\newcommand{\Z}{\mathbb{Z}_+}

\newcommand{\dtheta}[0]{\frac{\partial}{\partial
\theta}}
\newcommand{\dz}[0]{\frac{\partial}{\partial z}}

\newcommand{\Lx}[0]{x^{j + 1}
\frac{\partial}{\partial x} + (\frac{j +
1}{2})\varphi x^j \frac{\partial}{\partial
\varphi}}
   
\newcommand{\Gx}[0]{x^j \left( \frac{\partial}{\partial \varphi} -
\varphi \frac{\partial}{\partial x}\right)}

\newcommand{\twoLow}[0]{2w \frac{\partial}{\partial w} + \rho
\frac{\partial}{\partial \rho}} 
\newcommand{\Lw}[0]{w^{j + 1} \frac{\partial}{\partial w} + (\frac{j +
1}{2})\rho w^j \frac{\partial}{\partial \rho}} 
\newcommand{\Gw}[0]{w^j \left( \frac{\partial}{\partial \rho} -
\rho \frac{\partial}{\partial w}\right)}

\begin{document}

\title[$N=1$ SG-VOSAs]{The notion of $N=1$ 
supergeometric vertex operator superalgebra and the isomorphism theorem}

\author{Katrina Barron}
\address{Department of Mathematics, 255 Hurley Hall, University of Notre
Dame, Notre Dame  IN 46556}
\email{kbarron@nd.edu}
\thanks{The author was supported in part by an AAUW
Educational Foundation American Dissertation Fellowship and an NSF
Mathematical Sciences  Postdoctoral Research Fellowship.}




\begin{abstract} We introduce the notion of {\it $N=1$ supergeometric 
vertex operator superalgebra} motivated by the geometry underlying
genus-zero,  two-dimensional, holomorphic $N=1$ superconformal field theory. 
We then show, assuming the convergence of certain projective factors, that
the category of such objects is isomorphic to the category of $N=1$ 
Neveu-Schwarz vertex operator superalgebras. 
\end{abstract}

\maketitle

\section{Introduction and Preliminaries} 

\subsection{Introduction}
 
In this paper, we introduce the notion of {\it $N=1$ supergeometric vertex
operator superalgebra} and prove that the category of such objects is 
isomorphic to the category of $N=1$ Neveu-Schwarz vertex operator 
superalgebras with (or without) odd formal variables \cite{B vosas}. 
The notion of $N=1$ supergeometric vertex operator superalgebra is based
on the geometric structure of the moduli space of $N=1$ superspheres with
tubes along with a sewing operation on this moduli space and an action
of the symmetric groups \cite{B memoirs}.  This moduli space is the
fundamental  geometric structure underlying genus-zero, two-dimensional,
holomorphic $N=1$ superconformal field theory. 

Conformal field theory (or more specifically, string theory) and
related theories such as superconformal field theories (cf. \cite{BPZ},
\cite{FS}, \cite{V}, and \cite{S}) are the most promising attempts at
developing a  physical theory that combines all fundamental interactions of
particles,  including gravity. The geometry of this theory extends the use
of Feynman diagrams, describing the interactions of point particles whose
propagation in time sweeps out a line in space-time, to one-dimensional
``particles'' (strings) whose propagation in time sweeps out a
two-dimensional surface.  For two-dimensional genus-zero holomorphic
conformal field theory, algebraically, these interactions can be described
by products of vertex operators or more precisely, by vertex operator
algebras (cf. \cite{FLM}).  

Independently {}from the physical theory of conformal fields, vertex 
operator algebras arose naturally {}from the study of representations of
infinite-dimensional Lie algebras and in the construction of
representations of the Monster finite  simple group \cite{FLM}.  However, 
it wasn't until 1990 that a rigorous  mathematical interpretation of the
geometry of particle interactions and the algebraic description of those
interactions via vertex operator algebras was realized.  In \cite{H thesis}
and \cite{H book}, motivated by the geometric notions arising in conformal
field theory, Huang gives a precise geometric interpretation of the notion
of vertex operator algebra thus giving a rigorous basis for the geometric
and algebraic ``sewing" together of particle interactions.

In \cite{Fd}, Friedan describes the extension of the physical model of
conformal field theory to that of superconformal field theory with one
fermionic component ($N = 1$), and the notion of a superstring whose
propagation in time sweeps out a supersurface is introduced.  Whereas
conformal field theory attempts to describe the interactions of
bosons, $N = 1$ superconformal field theory attempts to describe the
interactions of boson-fermion pairs and has been discovered to be 
essential for a full description of all particle interactions.  This
theory requires an operator $D$ such that $D^2 = \partial / \partial z$. 
Such an operator arises naturally in supergeometry.

Within the framework of supergeometry (cf. \cite{D}, \cite{R} and
\cite{CR}) and motivated by $N = 1$ superconformal field theory, in
this paper, we extend the work of Huang to an $N=1$ supergeometric 
interpretation of $N=1$ Neveu-Schwarz vertex operator superalgebras.  This
paper relies heavily on the results of \cite{B memoirs}, and the notion
and notation of $N=1$ Neveu-Schwarz vertex operator superalgebra ($N=1$
NS-VOSA) as presented in \cite{B vosas}.  In \cite{B memoirs}, we studied
the moduli space of $N=1$ superspheres with tubes and defined a sewing
operation on this moduli space.  This sewing operation defines a
partial-operad structure on the moduli space  (cf. \cite{M}, \cite{HL1},
\cite{HL2}), and in \cite{B memoirs}, we thoroughly analyze this structure
both algebraically and analytically.  The most important aspect of this
work is that given two canonical $N=1$ superspheres with tubes
representing equivalence classes in the moduli space, we can completely
determine the resulting sewn $N=1$ supersphere as being in the
equivalence class of a specific canonical $N=1$ supersphere.  This is done
by explicitly solving the sewing equation and thus constructing the
uniformizing function which takes the resulting sewn $N=1$ supersphere to
the canonical $N=1$ supersphere representative.

In this paper, assuming the convergence of a certain series, denoted
$\Gamma$, (see Conjecture \ref{Gamma converges}), we define an algebra
over this partial operad of the moduli space of $N=1$ superspheres with
tubes and the sewing operation and extended to include the projective
factors associated to $\Gamma$, 
and call such an algebra an {\it $N=1$ supergeometric vertex operator
superalgebra ($N=1$ SG-VOSA)}.  We then prove that the category of such
objects is isomorphic to the category of (algebraic) $N=1$ NS-VOSAs.   
The fact that $N=1$ SG-VOSAs are defined over the geometric structure  
of the moduli space with sewing, and that the sewing operation is 
completely soluble, provides much more information for the  
corresponding algebraic structure of an $N=1$ SG-VOSA, and thus by
isomorphism of categories, for the more familiar structure of an $N=1$
NS-VOSA as well.  This includes a full understanding of the effects  
of a change of coordinates for $N=1$ NS-VOSA as studied in \cite{B 
change}, as well as the development of alternate notions of 
superconformality and relations to $N=2$ superconformal field theory as
studied in \cite{B alternate}.  In addition, this fully realized genus-zero
case is the first step toward understanding higher genus $N=1$
superconformal field theory.

This paper is organized as follows.  In the remainder of this section,
Section 1, we give a brief summary of some basic definitions.  In Sections
2 and 3, we give an overview of the results of \cite{B memoirs}.  This
overview is only meant to recall these results for the reader who is
already familiar with  \cite{B memoirs}.  Certainly for a full
understanding of these results, it is essential that one refer to \cite{B
memoirs}.

In Section 4, we discuss $(1/2)\mathbb{Z}$-graded vector spaces with
finite-dimensional homogeneous subspaces and define the
$t^{1/2}$-contraction of $n$-ary functions {}from such a vector
space to its algebraic completion.  This contraction of functions can
be thought of as being analogous to the sewing operation on the moduli
space of superspheres with tubes in that both contraction and sewing
define the formal substitution maps in the definition of a partial
(or partial pseudo-) operad (cf.  \cite{HL1}, \cite{HL2}, \cite{H
book}).   In the second subsection of Section 4, we introduce the notion
of {\it $N=1$ supergeometric vertex operator superalgebra ($N=1$ SG-VOSA)}
assuming the convergence of the projective factors for the central
charge.  

In Section 5, we recall {}from \cite{B vosas} the notion of $N=1$ 
Neveu-Schwarz vertex operator superalgebra ($N=1$ NS-VOSA) and several
consequences of the definition including the properties of associativity
and supercommutativity.

In Section \ref{alg {}from geom}, we construct an $N=1$ NS-VOSA with odd
formal variables {}from an $N=1$ SG-VOSA, and in Section \ref{geom {}from
alg}, we construct an $N=1$ SG-VOSA {}from an $N=1$ NS-VOSA with odd
formal variables.   

In Section \ref{iso}, we prove that the category of $N=1$ SG-VOSAs is
isomorphic to the category of $N=1$ NS-VOSAs with odd formal variables.

The results of this paper were first announced in
\cite{B announce} as part of the authors Ph.D. thesis \cite{B thesis}.   

\subsection{Grassmann algebras and the $N=1$ Neveu-Schwarz algebra}
\label{Grassmann}

Let $\mathbb{Z}_2$ denote the integers modulo two.  For a 
$\mathbb{Z}_2$-graded vector space $V = V^0 \oplus V^1$, define the 
{\it sign function} $\eta$ on the homogeneous subspaces of $V$ by 
$\eta(v) = i$ for $v \in V^i$, $i = 0,1$.  If $\eta(v) = 0$, we say 
that $v$ is {\it even}, and if $\eta(v) = 1$, we say that $v$ is 
{\it odd}. 

A {\it superalgebra} is an (associative) algebra $A$ (with identity $1
\in A$), such that: (i) $A$ is a $\mathbb{Z}_2$-graded algebra; (ii)
$ab = (-1)^{\eta(a) \eta(b)} ba$ for $a,b$ homogeneous in $A$.

The exterior algebra over a vector space $U$, denoted $\bigwedge(U)$, 
has the structure of a superalgebra.  Fix $U_L$ to be an 
$L$-dimensional vector space over $\mathbb{C}$ for $L \in \mathbb{N}$ 
with fixed basis $\{ \zeta_1, ..., \zeta_L \}$ such that $U_L \subset
U_{L+1}$.  We denote $\bigwedge(U_L)$ by $\bigwedge_L$ and call this the
{\it Grassmann algebra on $L$ generators}.  Note that $\bigwedge_L 
\subset \bigwedge_{L+1}$, and taking the direct limit as $L \rightarrow
\infty$, we have the {\it infinite Grassmann algebra} denoted by
$\bigwedge_\infty$.  We use the notation $\bigwedge_*$ to denote a
Grassmann algebra, finite or infinite.    

The $\mathbb{Z}_2$-grading of $\bigwedge_*$ is given explicitly by
\begin{eqnarray*}
\mbox{$\bigwedge_*^0$} \! &=& \! \Bigl\{a \in \mbox{$\bigwedge_*$} \; 
\big\vert \; a = \sum_{(i) \in I_*} a_{(i)}\zeta_{i_{1}}\zeta_{i_{2}} 
\cdots \zeta_{i_{2n}}, \; a_{(i)} \in \mathbb{C}, \; n \in \mathbb{N}
\Bigr\}\\ 
\mbox{$\bigwedge_*^1$} \! &=& \! \Bigl\{a \in \mbox{$\bigwedge_*$} \; 
\big\vert \; a = \sum_{(j) \in J_*} a_{(j)}\zeta_{j_{1}}\zeta_{j_{2}} 
\cdots \zeta_{j_{2n + 1}}, \; a_{(j)} \in \mathbb{C}, \; n \in \mathbb{N}
\Bigr\},
\end{eqnarray*}
where
\begin{eqnarray*}
I_* \! \!  &=&  \! \! \bigl\{ (i) = (i_1, i_2, \ldots, i_{2n}) \; | \; i_1 <
i_2 < \cdots < i_{2n}, \; i_l \in \{1, 2, ..., *\}, \; n \in
\mathbb{N} \bigr\}, \\ 
J_*  \! \! &=&  \! \! \bigl\{(j) = (j_1, j_2, \ldots, j_{2n + 1}) \; | \;
j_1 < j_2 < \cdots < j_{2n + 1}, \; j_l \in \{1, 2, ..., *\}, \;  n \in
\mathbb{N} \bigr\},
\end{eqnarray*}
with $\{1,2,...,*\}$ denoting $\{1,2,...,L\}$ if $\bigwedge_*$ is the
finite-dimensional Grassmann algebra $\bigwedge_L$ and denoting the 
positive integers $\mathbb{Z}_+$ if $\bigwedge_* = \bigwedge_\infty$. We
can also decompose $\bigwedge_*$ into {\it body}, $(\bigwedge_*)_B = 
\{ a_{(\emptyset)} \in \mathbb{C} \}$,  and {\it soul} 
\[(\mbox{$\bigwedge_*$})_S \; = \; \Bigl\{a \in \mbox{$\bigwedge_*$} \; \big\vert \; 
a = \! \! \! \sum_{ \begin{scriptsize} \begin{array}{c}
(k) \in I_* \cup J_*\\
k \neq (\emptyset)
\end{array} \end{scriptsize}} \! \! \!
a_{(k)} \zeta_{k_1} \zeta_{k_2} \cdots \zeta_{k_n}, \; a_{(k)} \in \mathbb{C} 
\Bigr\}\] 
subspaces such that $\bigwedge_* = (\bigwedge_*)_B \oplus 
(\bigwedge_*)_S$.  For $a \in \bigwedge_*$, we write $a = a_B + 
a_S$ for its body and soul decomposition.

A $\mathbb{Z}_2$-graded vector space $\mathfrak{g}$ is said to be a
{\it Lie superalgebra} if it has a bilinear operation $[\cdot,\cdot]$
such that for $u,v$ homogeneous in $\mathfrak{g}$: (i) $[u,v] \in
\mathfrak{g}^{(\eta(u) + \eta(v)) \mathrm{mod} \; 2}$; (ii)
skew-symmetry holds $[u,v] = -(-1)^{\eta(u)\eta(v)}[v,u]$; (iii) the
Jacobi identity holds $(-1)^{\eta(u)\eta(w)}[[u,v],w] +
(-1)^{\eta(v)\eta(u)}[[v,w],u] +  (-1)^{\eta(w)\eta(v)}[[w,u],v] = 0$.

Given a Lie superalgebra $\mathfrak{g}$ and a superalgebra 
$A$, the space $(A^0 \otimes \mathfrak{g}^0) \oplus (A^1 \otimes 
\mathfrak{g}^1)$ is a Lie algebra with bracket given by $[au , bv] =
(-1)^{\eta(b) \eta(u)} ab [u,v]$ for $a,b \in A$ and $u,v \in 
\mathfrak{g}$ homogeneous.  We call $(A^0 \otimes \mathfrak{g}^0) \oplus
(A^1 \otimes \mathfrak{g}^1)$ the {\it $A$-envelope of $\mathfrak{g}$}. 
Similarly, given two superalgebras $A$ and $\hat{A}$, we can form the
{\it $A$-envelope of $\hat{A}$} given by $(A^0 \otimes \hat{A}^0) \oplus
(A^1 \otimes \hat{A}^1)$ which is naturally an algebra.

For any $\mathbb{Z}_2$-graded associative algebra $A$ and for $u,v \in A$ 
of homogeneous sign, we can define $[u,v] = u v - (-1)^{\eta(u)\eta(v)} v 
u$, making $A$ into a Lie superalgebra.  The algebra of endomorphisms of 
$A$, denoted $\mbox{End} \; A$, has a natural $\mathbb{Z}_2$-grading 
induced {}from that of $A$, and defining $[X,Y] = X Y - (-1)^{\eta(X)
\eta(Y)} Y X$ for $X,Y$ homogeneous in $\mbox{End} \; A$, this gives 
$\mbox{End} \; A$ a Lie superalgebra structure.  An element $D \in 
(\mbox{End} \; A)^i$, for $i \in \mathbb{Z}_2$, is called a {\it 
superderivation of sign $i$} (denoted $\eta(D) = i$) if $D$ satisfies the 
super-Leibniz rule $D(uv) = (Du)v + (-1)^{\eta(D) \eta(u)} uDv$
for $u,v \in A$ homogeneous.

The Virasoro algebra is the Lie algebra with basis consisting of the
central element $d$ and $L_n$, for $n \in \mathbb{Z}$, and commutation
relations  
\begin{equation}\label{Virasoro relation}
[L_m ,L_n] = (m - n)L_{m + n} + \frac{1}{12} (m^3 - m) \delta_{m + n 
, 0} \; d ,
\end{equation}
for $m, n \in \mathbb{Z}$.  The $N=1$ Neveu-Schwarz Lie superalgebra,
denoted $\mathfrak{ns}$, is a super-extension of the Virasoro algebra by
the odd elements $G_{n + 1/2}$, for $n \in \mathbb{Z}$, with 
supercommutation relations 
\begin{eqnarray}
\left[ G_{m + \frac{1}{2}},L_n \right] &=& \Bigl(m - \frac{n - 1}{2} \Bigr) G_{m
+ n + \frac{1}{2}}  \label{Neveu-Schwarz relation1} \\  
\left[ G_{m + \frac{1}{2}} , G_{n - \frac{1}{2}} \right] &=& 2L_{m +
n} + \frac{1}{3} (m^2 + m) \delta_{m + n , 0} \; d \label{Neveu-Schwarz relation2}
\end{eqnarray}
in addition to (\ref{Virasoro relation}).

\subsection{Superanalytic and superconformal superfunctions}

Let $U$ be a subset of $\bigwedge_*$, and write $U = U^0 \oplus U^1$
for the decomposition of $U$ into even and odd subspaces.  Let $z$ 
be an even variable in $U^0$ and $\theta$ an odd variable in $U^1$.  We
call $H: U \longrightarrow \bigwedge_*$, mapping $(z,\theta) \mapsto
H(z,\theta)$, a {\it $\bigwedge_*$-superfunction on $U$ in (1,1)-variables}.

Let $z_B$ be a complex variable and $h(z_B)$ a complex
analytic function in some open set $U_B \subset \mathbb{C}$.  For $z$ a
variable in $\bigwedge_*^0$, we define $h(z)$ to be the Taylor expansion
about the body of $z = z_B + z_S$. Then $h(z)$ is well defined (i.e.,
convergent) in the open neighborhood $\{z = z_B + z_S \in \bigwedge_*^0 \; |
\; z_B \in U_B \} = U_B \times(\bigwedge_*^0)_S \subseteq
\bigwedge_*^0$.   Since $h(z)$ is algebraic in each $z_{(i)}$, for $(i)
\in I_*$, it follows that $h(z)$ is complex analytic in each of the
complex variables $z_{(i)}$.

For $n \in \mathbb{N}$, we introduce the notation $\bigwedge_{*>n}$ to
denote a finite Grassmann algebra $\bigwedge_L$ with $L > n$ or an 
infinite Grassmann algebra.  We will use the corresponding index 
notations for the corresponding indexing sets $I_{*>n}$ and $J_{*>n}$.
A {\it superanalytic $\bigwedge_{*>0}$-superfunction in
$(1,1)$-variables} $H$ is a $\bigwedge_{*>0}$-superfunction in
$(1,1)$-variables of the form
\begin{eqnarray*}
H(z, \theta) \! \! &=& \! \! (f(z) + \theta \xi(z), \psi(z) + \theta
g(z)) \\ 
&=& \! \! \!  \Biggl(\sum_{(i) \in I_{ * - 1}} \!
\! \! f_{(i)} (z) \zeta_{i_1} \zeta_{i_2} \cdots \zeta_{i_n} + \;
\theta \! \! \! \sum_{(j) \in J_{ * - 1}} \! \! \!
\xi_{(j)} (z) \zeta_{j_1} \zeta_{j_2} \cdots \zeta_{j_{2n + 1}},
\Biggr. \\ 
& & \qquad \quad \Biggl. \sum_{(j) \in J_{ * - 1}} \!  \! \!
\psi_{(j)} (z) \zeta_{j_1} \zeta_{j_2} \cdots \zeta_{j_{2n + 1}} + \;
\theta \! \! \! \sum_{(i) \in I_{ * - 1}} \! \! \!
g_{(i)} (z) \zeta_{i_1} \zeta_{i_2} \cdots \zeta_{i_n} \Biggr)
\end{eqnarray*}
where $f_{(i)}(z_B)$, $g_{(i)}(z_B)$, $\xi_{(j)}(z_B)$, and
$\psi_{(j)}(z_B)$ are all complex analytic in some non-empty open
subset $U_B \subseteq \mathbb{C}$.  If each $f_{(i)}(z_B)$,
$g_{(i)}(z_B)$, $\xi_{(j)}(z_B)$,  and $\psi_{(j)}(z_B)$ is complex
analytic in $U_B \subseteq \mathbb{C}$, then $H(z,\theta)$ is well
defined (i.e., convergent) for $\{(z,\theta) \in \bigwedge_{*>0} \; | \;
z_B \in U_B \} = U_B \times (\bigwedge_{*>0})_S$.  Consider the topology
on $\bigwedge_{*}$ given by the product of the usual topology on 
$(\bigwedge_{*})_B = \mathbb{C}$ and the trivial topology on 
$(\bigwedge_{*})_S$.  This topology on $\bigwedge_{*}$ is called 
the {\it DeWitt topology}.  The natural domain of any superanalytic
$\bigwedge_{*>0}$-superfunction is an open set in the DeWitt topology 
on $\bigwedge_{*>0}$.

Since $1/a = \sum_{n \in \mathbb{N}} (-1)^n a_S^n/a_B^{n + 1}$ is well
defined if and only if $a_B \neq 0$, the set of invertible elements in
$\bigwedge_*$, denoted $\bigwedge_*^\times$, is given by
$\bigwedge_*^\times = \{a \in \mbox{$\bigwedge_*$} \; | \; a_B \neq 0
\}$. Note that $\bigwedge_*^\times \subset \bigwedge_*$ is open in the
DeWitt topology.

We define the (left) partial derivatives $\partial/\partial z$ and
$\partial/\partial \theta$ acting on superfunctions which are
superanalytic in some DeWitt open neighborhood $U$ of $(z,\theta) \in
\bigwedge_{*>0}$ by
\begin{eqnarray*}
\Delta z \left(\dz H(z,\theta) \right) + O((\Delta z)^2) &=& H(z +
\Delta z, \theta) - H(z,\theta) \\
\Delta \theta \left(\dtheta H(z,\theta) \right) &=& H(z, \theta +
\Delta \theta) - H(z,\theta) 
\end{eqnarray*}
for all $\Delta z \in \bigwedge_{*>0}^0$ and $\Delta \theta \in
\bigwedge_{*>0}^1$ such that $z + \Delta z \in U^0 = U_B \times
(\bigwedge_{*>0}^0)_S$ and $\theta + \Delta \theta \in U^1 =
\bigwedge_{*>0}^1$.  Note that $\partial/\partial z$ and $\partial/
\partial \theta$ are endomorphisms of the  superalgebra of superanalytic
$(1,1)$-superfunctions, and in fact, are  even and odd superderivations,
respectively.

If $h(z_B)$ is complex analytic in an open neighborhood of the complex
plane, then $h(z_B)$ has a Laurent series expansion in $z_B$, given by 
$h(z_B) = \sum_{l \in \mathbb{Z}} c_l z_B^l$, for $c_l \in \mathbb{C}$, 
and we have $h(z) = \sum_{n \in \mathbb{N}} (z_S^n/n!) h^{(n)}(z_B) =
\sum_{l \in \mathbb{Z}} c_l (z_B + z_S)^l = \sum_{l \in \mathbb{Z}} c_l z^l$
where $(z_B + z_S)^l$, for $l \in \mathbb{Z}$, is always understood to 
mean expansion in positive powers of the second variable, in this
case $z_S$.  Thus if $H$ is a $\bigwedge_{*>0}$-superfunction in 
$(1,1)$-variables which is superanalytic in a (DeWitt) open 
neighborhood, $H$ can be expanded as
\begin{equation}\label{Laurent series}
H(z, \theta) = \biggl(\sum_{l \in \mathbb{Z}} a_l z^l + \theta \sum_{l
\in \mathbb{Z}} n_l z^l, \sum_{l \in \mathbb{Z}} m_l z^l + \theta
\sum_{l \in \mathbb{Z}} b_l z^l \biggr)
\end{equation}   
for $a_l, b_l \in \bigwedge_{ * - 1}^0$ and $m_l, n_l \in
\bigwedge_{ * - 1}^1$.  

Define $D$ to be the odd superderivation $D = \partial / \partial \theta +
\theta \partial / \partial z$ acting on superanalytic
$\bigwedge_{*>0}$-superfunctions in $(1,1)$-variables. Then $D^2 = \partial /
\partial z$, and if $H(z,\theta) = (\tilde{z},\tilde{\theta})$ is
superanalytic in some  DeWitt open subset, then $D$ transforms under
$H(z,\theta)$ by $D = (D\tilde{\theta})\tilde{D} + (D\tilde{z} -
\tilde{\theta} D \tilde{\theta})\tilde{D}^2$.  We define a {\it
superconformal $(1,1)$-superfunction} on a DeWitt open  subset $U$ of
$\bigwedge_{*>0}$ to be a superanalytic superfunction $H$ under which $D$
transforms homogeneously of degree one.  Thus a superanalytic function
$H(z,\theta) =  (\tilde{z}, \tilde{\theta})$ is superconformal if and only
if, in  addition to being superanalytic, $H$ satisfies
$D\tilde{z} - \tilde{\theta} D\tilde{\theta} = 0$,
for $D \tilde{\theta}$ not identically zero, thus transforming $D$ by 
$D = (D\tilde{\theta})\tilde{D}$.  

\subsection{Super-Riemann surfaces}

A {\em DeWitt $(1,1)$-dimensional manifold over 
$\bigwedge_*$} is a topological space $X$ with a countable basis which 
is locally homeomorphic to an open subset of $\bigwedge_*$ in the  
DeWitt topology.  A {\em DeWitt $(1,1)$-chart on $X$ over $\bigwedge_*$} 
is a pair $(U, \Omega)$ such that $U$ is an open subset of $X$ and $\Omega$ 
is a homeomorphism of $U$ onto an open subset of $\bigwedge_*^0$ in the 
DeWitt topology.  A {\em superanalytic atlas of DeWitt $(1,1)$-charts on $X$ 
over $\bigwedge_{* > 0}$} is a family of charts 
$\{(U_{\alpha}, \Omega_{\alpha})\}_{\alpha \in A}$ satisfying: 
(i) each $U_{\alpha}$ is open in $X$, and $\bigcup_{\alpha \in A} U_{\alpha}
= X$; (ii) each $\Omega_{\alpha}$ is a homeomorphism {}from $U_{\alpha}$ to a
(DeWitt) open set in $\bigwedge_{* >0}$
such that $\Omega_{\alpha} \circ \Omega_{\beta}^{-1}: \Omega_{\beta}(U_\alpha
\cap U_\beta) \longrightarrow \Omega_{\alpha}(U_\alpha \cap U_\beta)$
is superanalytic for all non-empty $U_{\alpha} \cap U_{\beta}$.  Such an
atlas is called {\em maximal} if, given any chart $(U, 
\Omega)$ such that $\Omega \circ \Omega_{\beta}^{-1} : \Omega_{\beta} (U \cap
U_\beta) \longrightarrow \Omega (U \cap U_\beta)$ 
is a superanalytic homeomorphism for all $\beta$, then $(U, \Omega)
\in \{(U_{\alpha}, \Omega_{\alpha})\}_{\alpha \in A}$.  A {\em DeWitt
$(1,1)$-supermanifold over $\bigwedge_{* >0}$} is a DeWitt
$(1,1)$-dimensional topological space $M$ together with a maximal
superanalytic atlas of DeWitt $(1,1)$-charts over
$\bigwedge_{* > 0}$.  

Given a DeWitt $(1,1)$-supermanifold $M$ over $\bigwedge_{* > 0}$,
define an equivalence relation $\sim$ on M by letting $p \sim q$
if and only if there exists $\alpha \in A$ such that $p,q \in U_\alpha$
and $\pi_B (\Omega_\alpha (p)) = \pi_B (\Omega_\alpha (q))$
where $\pi_B$ is the canonical projection of $\bigwedge_{*>0}$ onto its body.
Let $p_B$ denote the equivalence class of $p$ under this 
equivalence relation.  Define the {\it body} $M_B$ of $M$ to be the 
one-dimensional complex manifold with analytic structure given by the 
coordinate charts $\{((U_\alpha)_B, (\Omega_\alpha)_B) \}_{\alpha \in A}$ 
where $(U_\alpha)_B = \{ p_B \; | \; p \in U_\alpha \}$, and
$(\Omega_\alpha)_B : (U_\alpha)_B \longrightarrow \mathbb{C}$ is given
by $(\Omega_\alpha)_B (p_B) = \pi_B \circ \Omega_\alpha (p)$.

For any DeWitt $(1,1)$-supermanifold $M$, its body $M_{B}$ is a Riemann
surface. A {\em super-Riemann surface over $\bigwedge_{*>0}$} is a
DeWitt $(1,1)$-supermanifold over $\bigwedge_{*>0}$ with coordinate
atlas $\{(U_{\alpha}, \Omega_{\alpha})\}_{\alpha \in A}$ such that the
coordinate transition functions $\Omega_{\alpha} \circ
\Omega_{\beta}^{-1}$ in addition to being superanalytic are also
superconformal for all non-empty $U_{\alpha} \cap U_{\beta}$. Since 
the condition that the coordinate transition functions be superconformal 
instead of merely superanalytic is such a strong condition (unlike in 
the nonsuper case), we again stress the distinction between a 
supermanifold which has {\it superanalytic} transition functions versus 
a super-Riemann surface which has {\it superconformal} transition 
functions.  

\section{The moduli space of $N=1$ superspheres with tubes and the sewing
operation}

In \cite{B memoirs}, we define the moduli space $SK(n)$ of $N=1$
superspheres with $n$ incoming tubes and one outgoing tube, define a sewing
operation on this space and give a detailed algebraic and analytic study of
this space with this sewing operation.  In this section, we recall some
results of this work that we will need in order to define the notion of
$N=1$ supergeometric vertex operator superalgebra and to prove the
Isomorphism Theorem.

\subsection{The moduli space of $N=1$ superspheres with tubes}

By {\em $N=1$ supersphere} we will mean a (superconformal) 
super-Riemann surface over $\bigwedge_{*>0}$ such that its body is a 
genus-zero one-dimensional connected compact complex manifold.  {}From
now on we will refer to such an object as a {\it supersphere}.

A {\em supersphere with $1+n$ tubes} for $n \in \mathbb{N}$, is a
supersphere $S$ with 1 negatively oriented point $p_0$ and $n$
positively oriented points $p_1,...,p_n$ (we call them {\em punctures}) on
$S$ which all have distinct bodies and with local superconformal 
coordinates $(U_0,\Omega_0),...,(U_n, \Omega_n)$ vanishing at the  punctures
$p_0,...,p_n$, respectively.  We denote this structure by
$(S;p_0,...,p_n;(U_0, \Omega_0),...,(U_n, \Omega_n))$ .  We will always
order the punctures so that the negatively oriented puncture is $p_0$. 

Let $S_1 = (S_1;p_0,...,p_m;(U_0, \Omega_0),...,(U_m, \Omega_m))$
be a supersphere with $1+m$ tubes, for $m \in \mathbb{N}$, and let
$S_2 = (S_2;q_0,...,q_n;(V_0, \Xi_0),...,(V_n, \Xi_n))$ be a supersphere
with $1+n$ tubes, for $n \in \mathbb{N}$.  A map $F : S_1 \rightarrow
S_2$ will be said to be superconformal if $\Xi_\beta \circ F \circ
\Omega_\alpha^{-1}$ is superconformal for all  charts $(U_\alpha,
\Omega_\alpha)$ of $S_1$, for all charts $(V_\beta, \Xi_\beta)$ of $S_2$,
and for all $(z,\theta) \in \Omega_\alpha (U_\alpha)$ such that $F \circ
\Omega_\alpha^{-1}  (z,\theta) \in V_\beta$.  If $m=n$ and there is a
superconformal  isomorphism $F : S_1 \rightarrow S_2$ such that for each
$j = 0,...,  n$, we have that $F(p_j) = q_j$ and $ \Omega_j |_{W_j} = 
\Xi_j \circ F |_{W_j}$, for $W_j$ some DeWitt neighborhood of $p_j$, then
we say that these two  superspheres with $1 + n$ tubes are {\it
superconformally equivalent} and $F$ is a {\it superconformal
equivalence} {}from $S_1$ to $S_2$. Thus the superconformal equivalence
class of a supersphere with tubes depends only on the supersphere, the
punctures, and the germs of the local coordinate maps vanishing at the
punctures.

The collection of all superconformal equivalence classes of
superspheres over $\bigwedge_{*>0}$ with $1+n$ tubes, for $n \in 
\mathbb{N}$, is called the {\it moduli space of superspheres over 
$\bigwedge_{*>0}$ with $1+n$ tubes}.  The collection of all 
superconformal equivalence classes of superspheres over 
$\bigwedge_{*>0}$ with tubes is called the {\it moduli space of 
superspheres over $\bigwedge_{*>0}$ with tubes}.

Let $S\hat{\mathbb{C}}$ be the supersphere with superconformal 
structure given by the covering of local coordinate neighborhoods 
$\{ U_\sou, U_\nor \}$ and the local coordinate maps
$\sou : U_\sou \longrightarrow \bigwedge_{*>0}$ and 
$\nor: U_\nor  \longrightarrow  \bigwedge_{*>0}$, 
which are homeomorphisms of $U_\sou$ onto $\bigwedge_{*>0}$ and
$U_\nor$ onto $\bigwedge_{*>0}$, respectively, such that 
$\sou \circ \nor^{-1} : \bigwedge_{*>0}^\times \longrightarrow
\bigwedge_{*>0}^\times$ is given by  $(z,\theta)  \mapsto 
(1/z, i \theta/z) = \sou \circ \nor^{-1}(z,\theta) = I(z,\theta)$.
Thus the body of $S\hat{\mathbb{C}}$ is the Riemann sphere, 
$(S\hat{\mathbb{C}})_B = \hat{\mathbb{C}} = \mathbb{C} \cup \{\infty\}$,  
with coordinates $z_B$ near 0 and $1/z_B$ near $\infty$.   We will
call $S\hat{\mathbb{C}}$ the {\em super-Riemann sphere} and
will refer to $\nor^{-1} (0)$ as {\em the point at $(\infty, 
0)$} or just {\em the point at infinity} and to 
$\sou^{-1}(0)$ as {\em the point at $(0, 0)$} or just {\em the point at 
zero}.  

Let $T: S\hat{\mathbb{C}} \longrightarrow S\hat{\mathbb{C}}$ be a
superconformal automorphism.  Then defining $T_\sou = \sou \circ T \circ
\sou^{-1}$, we have that $T$ must be a superprojective
transformation satisfying 
\[ T_\sou : \mbox{$\bigwedge_{*>0}$} \smallsetminus \bigl( \{(- d/c)_B\}
\times (\mbox{$\bigwedge_{*>0}$})_S \bigr) \longrightarrow  \mbox{$\bigwedge_{*>0}$} 
\smallsetminus \bigl(\{(a/c)_B \} \times (\mbox{$\bigwedge_{*>0}$})_S \bigr) 
\hspace{.5in}\]  
\[\hspace{1.4in} (z, \theta) \mapsto \left( \frac{az + b}{cz + d} +
\theta \frac{\gamma z + \delta}{(cz + d)^2} , \frac{\gamma z + \delta}{cz
+ d} + \theta \frac{1 + \frac{1}{2} \delta \gamma }{cz + d} \right) , \]
for $a,b,c,d \in \bigwedge_{ * - 1}^0$, $\gamma, \delta \in \bigwedge_{ *
- 1}^1$ and $ad - bc = 1$, and $T$ is determined uniquely by $T_\sou$.

In \cite{CR}, Crane and Rabin prove the uniformization theorem for
super-Riemann surfaces. For super-Riemann surfaces with genus-zero compact 
bodies, the theorem states that any such super-Riemann surface is
superconformally equivalent to the super-Riemann sphere
$S\hat{\mathbb{C}}$.  Thus as is proved in \cite{B memoirs}, any
supersphere with $1+n$ tubes is superconformally equivalent to
\[(S\hat{\mathbb{C}}; \infty, (z_1,\theta_1),...,(z_{n-1}, \theta_{n-1}),
0; H_0, H_1,...,H_n),\]
where $(z_0, \theta_0) = \infty = \nor^{-1}(0)$ is the coordinate of the
outgoing puncture, $(z_n, \theta_n) = 0 = \sou^{-1}(0)$ is the coordinate
of the last incoming puncture, $(z_i,\theta_i) \in
(\bigwedge_{*>0})^\times$ with $(z_i)_B \neq (z_j)_B$ for $i\neq j$, and
$H_i$ is a coordinate {}from an open set around the $i$-th puncture to
$\bigwedge_{*>0}$ which maps the puncture to zero.

We now temporarily switch to a formal algebraic setting in order to
better understand superconformal local coordinates.  Let $R$ be a
superalgebra, let $x$ be an even formal variable, and let 
$\varphi$ be an odd formal variable.  By this we mean that $x$ commutes 
with all formal variables and all superalgebra elements and $\varphi$ 
anticommutes with all odd formal variables and all odd superalgebra 
elements but commutes with even elements.  Note that $R((x))[\varphi]$ and
$R((x^{-1}))[\varphi]$ are superalgebras.  For $j \in \mathbb{Z}$, consider 
the even superderivations 
\begin{equation}\label{L notation}
L_j(x,\varphi) = - \biggl( \Lx \biggr)
\end{equation}
and the odd superderivations
\begin{equation}\label{G notation}
G_{j -\frac{1}{2}} (x,\varphi) = - \Gx
\end{equation}
in $\mbox{Der} (R[[x,x^{-1}]]))[\varphi])$.  Then $L_j(x,\varphi)$ and 
$G_{j-1/2}(x,\varphi)$ give a representation of the $N=1$ Neveu-Schwarz
Lie superalgebra, $\mathfrak{ns}$, with central charge zero.

For any formal power series $H \in R[[x,x^{-1}]][\varphi]$, we say that 
$H(x,\varphi) = (\tilde{x},\tilde{\varphi})$ is {\it superconformal} if
$D\tilde{x} = \tilde{\varphi}D\tilde{\varphi}$ for
$D = \partial/\partial \varphi + \varphi \partial/\partial x$. 

Let $(R^0)^\times$ denote the invertible elements in $R^0$.
Let $(R^0)^\infty$ be the set of all sequences $A = \{A_j\}_{j \in \Z}$ of
even  elements in $R$, let $(R^1)^\infty$ be the set of all sequences 
$M = \{M_{j - 1/2}\}_{j \in \Z}$ of odd elements in $R$, and let $R^\infty =
(R^0)^\infty \oplus (R^1)^\infty$.  In \cite{B thesis} and \cite{B memoirs},
we show that there is a bijection between the set of formal superconformal
power series $H \in R[[x,x^{-1}]][\varphi]$ which vanish at zero, i.e.,
formal superconformal power series of the form 
\begin{multline}\label{formal Laurent series}
H(x, \varphi) = \biggl(\asqrt^2 \Bigl(x + \sum_{j \in \mathbb{Z}_+} a_j
x^{j+1} + \varphi \sum_{j \in \mathbb{Z}_+} n_j x^j \Bigr), \\ 
\asqrt \Bigl( \sum_{j \in \mathbb{Z}_+} m_{j -\frac{1}{2}} x^j + \varphi
(1+\sum_{j \in \mathbb{Z}_+} b_j x^j ) \Bigr)\biggr)
\end{multline}   
with $\asqrt \in (R^0)^\times$ and $(a,m) = \{(a_j,m_{j-\frac{1}{2}})\}_{j
\in \mathbb{Z}_+}$, and satisfying $D\tilde{x} = \tilde{\varphi} D
\tilde{\varphi}$, and formal power series of the form
\begin{equation}
\exp\Biggl( \! - \! \sum_{j \in \Z} \Bigl( A_j L_j(x,\varphi) + M_{j -
\frac{1}{2}} G_{j -\frac{1}{2}} (x,\varphi) \Bigr) \! \Biggr) \cdot
(\asqrt^2 x, \asqrt \varphi) . 
\end{equation}
This defines a bijection 
\begin{eqnarray}
E : R^\infty &\longrightarrow& R^\infty\\
(A,M) = \{(A_j, M_{j-\frac{1}{2}})\}_{j \in \mathbb{Z}_+} &\mapsto& (a,m)
= \{(a_j,m_{j-\frac{1}{2}})\}_{j \in \mathbb{Z}_+}.
\nonumber
\end{eqnarray}
As in \cite{B memoirs}, this allows us to define a map $\tilde{E}$ {}from
$R^\infty$ to  the set of all formal superconformal power series in
$xR[[x]][\varphi]$ with leading even coefficient of $\varphi$ equal to one,
by defining
\begin{eqnarray}
\varphi \tilde{E}^0(A,M)(x,\varphi) &=& \varphi \biggl( x + \sum_{j \in
\Z} E_j(A,M) x^{j + 1} \biggr), \label{Etilde1}\\
\varphi \tilde{E}^1(A,M)(x,\varphi) &=& \varphi \sum_{j \in \Z} E_{j -
\frac{1}{2}}(A,M) x^j , \label{Etilde2}
\end{eqnarray}
and letting $\tilde{E}(A,M) (x, \varphi)$ be the unique formal
superconformal power series with even coefficient of $\varphi$ equal to one 
such that the even and odd components of $\tilde{E}$ satisfy (\ref{Etilde1})
and (\ref{Etilde2}), respectively.  For $\asqrt \in (R^0)^\times$, we define
a map $\hat{E}$ {}from $(R^0)^\times \times R^\infty$ to the set of all
formal superconformal power series in $xR[[x]][\varphi]$ with invertible
leading even coefficient of $\varphi$, by defining 
\begin{eqnarray}
\hat{E}(\asqrt,A,M)(x,\varphi) &=&  (\hat{E}^0(\asqrt, A, M)
(x,\varphi), \hat{E}^1(\asqrt, A, M) (x,\varphi)) \\ \nonumber
&=& (\asqrt^2 \tilde{E}^0(A,M)(x,\varphi), \asqrt
\tilde{E}^1(A,M)(x,\varphi).  
\end{eqnarray}
The following proposition is proved in \cite{B memoirs}.

\begin{prop}\label{above}(\cite{B memoirs})
The map $\hat{E}$ {}from $(R^0)^{\times} \times R^{\infty}$ to the
set of all formal superconformal power series $H(x, \varphi) \in
xR[[x]][\varphi]$ of the form
\begin{equation}\label{Ehat} 
\varphi H(x,\varphi) = \varphi \Biggl(\asqrt^2 \Bigl( x + \sum_{j \in \Z}
a_j x^{j + 1} \Bigr), \asqrt \sum_{j \in \Z} m_{j - \frac{1}{2}} x^j
\Biggr)  
\end{equation} 
and with even coefficient of $\varphi$ equal to $\asqrt$ 
for $(\asqrt,a,m) \in (R^0)^{\times} \times R^\infty$, is a
bijection.

The map $\tilde{E}$ {}from $R^{\infty}$ to the set of
formal superconformal power series of the form (\ref{Ehat}) with
$\asqrt = 1$ and even coefficient of $\varphi$ equal to 1 is also a 
bijection. 

In particular, we have inverses $\tilde{E}^{-1}$ and $\hat{E}^{-1}$.
\end{prop}

Let
\begin{multline*}
\mathcal{H}  =  \bigl\{ (A,M) \in \mbox{$\bigwedge_\infty^\infty$} 
\; | \; \tilde{E}(A,M)(z,\theta) \; \mbox{is an absolutely convergent power }
\bigr. \\
\bigl. \mbox{series in some neighborhood of } (z,\theta) = 0 \bigr\},
\end{multline*} 
and for $n \in \Z$, let
\[ SM^{n - 1} = \bigl\{ \bigl((z_1, \theta_1),...,(z_{n-1}, \theta_{n-1})\bigr) \; | \;
(z_i, \theta_i) \in \mbox{$\bigwedge_\infty^\times$}, \; (z_i)_B \neq (z_j)_B , \;
\mbox{for} \; i \neq j \bigr\} . \] 
Note that for $n=1$, the set $SM^0$ has exactly one element.

For convenience, define the superconformal shift
\begin{eqnarray*}
s_{(z,\theta)} : \mbox{$\bigwedge_\infty$} &\rightarrow&
\mbox{$\bigwedge_\infty$} \\
(w,\rho) &\mapsto&  (w - z -\rho\theta,\rho - \theta).
\end{eqnarray*} 
We also define the linear operator $\asqrt^{-2L_0(x,\varphi)} =
\asqrt^{2x\frac{\partial}{\partial x} + \varphi \frac{\partial}{\partial
\varphi}}$ on $R[[x,x^{-1}]][\varphi]$, by 
\begin{equation}
\asqrt^{-2L_0(x,\varphi)} \cdot c \varphi^k x^n = c \asqrt^{2n + k}
\varphi^k x^n ,
\end{equation}
for $c \in R$, $k = 0,1$, and $n \in \mathbb{Z}$, and by extending
linearly.

For any supersphere with $1 + n$ tubes over $\bigwedge_\infty$, we can write
the power series  expansion of the local coordinate at the
$i$-th puncture
$(z_i,\theta_i)$,  for $i = 1,...,n$, as 
\begin{eqnarray*}
H_i (w, \rho) &=& \exp
\Biggl(\! - \! \sum_{j \in \Z} \Bigl( A^{(i)}_j L_j(x,\varphi)
+  M^{(i)}_{j - \frac{1}{2}} G_{j -\frac{1}{2}} \Bigr) \! \Biggr) \cdot\\
& & \hspace{1.4in} \left. \cdot (\asqrt^{(i)})^{-2L_0(x,\varphi)} \cdot  (x, 
\varphi) \right|_{(x,\varphi) = s_{(z_i,\theta_i)}(w,\rho) }\\
&=& \hat{E}(\asqrt^{(i)}, A^{(i)}, M^{(i)}) \circ s_{(z_i,
\theta_i)}(w,\rho) 
\end{eqnarray*}
for $(\asqrt^{(i)}, A^{(i)}, M^{(i)}) \in (\bigwedge_\infty^0)^\times 
\times \mathcal{H}$, and we can write the power series expansion of the
local coordinate at $\infty$ as  
\begin{eqnarray*}
H_0 (w, \rho) &=& \exp \Biggl(\sum_{j \in \Z} \Bigl( A^{(0)}_j
L_{-j}(w,\rho) +  M^{(0)}_{j - \frac{1}{2}} G_{-j + \frac{1}{2}}(w,\rho)
\Bigr) \! \Biggr) \cdot \Bigl(\frac{1}{w}, \frac{i \rho}{w}\Bigr)  \\ 
&=& \tilde{E}(A^{(0)}, -iM^{(0)}) \Bigl(\frac{1}{w}, \frac{i \rho}{w}\Bigr)
\end{eqnarray*}
for $(A^{(0)},-iM^{(0)}) \in \mathcal{H}$.  But $(A^{(0)},-iM^{(0)}) \in 
\mathcal{H}$ if and only if $(A^{(0)},M^{(0)}) \in \mathcal{H}$. {}From \cite{B memoirs}, we have the following proposition.

\begin{prop}\label{moduli2}(\cite{B memoirs})
The moduli space of superspheres over $\bigwedge_\infty$ with $1 + n$
tubes, for $n \in \Z$,  can be identified with the set   
\begin{equation}
SK(n) = SM^{n-1} \times \mathcal{H} \times \bigl((\mbox{$\bigwedge_\infty^0$})^\times 
\times \mathcal{H}\bigr)^n .
\end{equation}
The moduli space of superspheres with one tube can be identified 
with the set
\[SK(0) = \bigl\{(A,M) \in \mathcal{H} \; | \; (A_1, M_{\frac{1}{2}}) =
(0,0) \bigr\} .\] 
\end{prop}

Let $SK = \bigcup_{n \in \mathbb{N}} SK(n)$ denote the moduli space of
superspheres over $\bigwedge_\infty$ with tubes.  The actual elements of
$SK$ give the data for a {\it canonical supersphere} representative of a
given equivalence class of superspheres with tubes modulo superconformal
equivalence.  Any element of $SK(n)$, for 
$n \in \Z$, can be written as 
\[\bigl((z_1, \theta_1),...,(z_{n-1}, \theta_{n-1}); (A^{(0)}, M^{(0)}),
(\asqrt^{(1)}, A^{(1)}, M^{(1)}), ...,(\asqrt^{(n)}, A^{(n)}, M^{(n)} 
)\bigr)  \]
where $(z_1, \theta_1),...,(z_{n-1}, \theta_{n-1}) \in SM^{n-1}$,
$(A^{(0)}, M^{(0)}) \in \mathcal{H}$, and 
\[(\asqrt^{(1)}, A^{(1)}, M^{(1)}),...,(\asqrt^{(n)}, A^{(n)}, 
M^{(n)}) \in (\mbox{$\bigwedge_\infty^0$})^\times 
\times \mathcal{H} . \]
Thus for an element $Q \in SK$, we can think of $Q$ as consisting of the
above data, or as being a canonical supersphere with tubes corresponding
to that data. 

Note that above we are working over $\bigwedge_\infty$.  But we can
always restrict variables and coefficients to be in a finite dimensional
Grassmann subalgebra of $\bigwedge_\infty$ if we so wish as long as care
is taken to see that if we employ any partial derivatives with respect to
odd variables that these partials are well defined. 

\subsection{A left action of the symmetric group $S_n$ on the moduli space
$SK(n)$}\label{symmetric group}

Let $S_n$ be the group of  permutations on $n$ letters, for $n \in
\Z$.   In \cite{B memoirs}, we defined an action of $S_n$ on
$SK(n)$.  We would like to take this opportunity to point out that
although this action was stated to be a right action, it is in fact a
left action.  We define this left action by first defining an action of
$S_{n - 1}$ on $SK(n)$ given by permuting the ordering of the first $n -
1$ positively oriented  punctures and their local coordinates and then
extending this to an action of $S_n$ on $SK(n)$, by defining the action
of the transposition $(n-1 \; n)$ on $SK(n)$.

To define the action of $S_{n-1}$ on $SK(n)$, we note that an ordering 
or labeling of the first $n-1$ punctures of a supersphere with $n$-tubes with
the first $n-1$ punctures at the points $p_1, p_2,...,p_{n-1}$ is a
bijection
\begin{eqnarray*}
l: \{1,2,...,n-1\} \longrightarrow \{p_1,p_2,...,p_{n-1}\}. 
\end{eqnarray*}
Then a left action on any ordering $l$ is given by $\sigma \cdot l (i) = l
(\sigma^{-1}(i))$, for $\sigma \in S_{n-1}$ and $i = 1,..., n-1$.  Thus,
for $\sigma \in S_{n-1}$, and $Q \in SK(n)$ given by
\[Q = \bigl((z_1, \theta_1),...,(z_{n-1}, \theta_{n-1}); (A^{(0)},M^{(0)}), 
(\asqrt^{(1)}, A^{(1)}, M^{(1)}),..., (\asqrt^{(n)}, A^{(n)}, M^{(n)})\bigr), \]
$\sigma$ acting on $Q$, denoted $\sigma \cdot Q$, is given by
\begin{multline*}
\sigma \cdot Q = \bigl((z_{\sigma^{-1}(1)}, \theta_{\sigma^{-1}(1)}), ...,
(z_{\sigma^{-1}(n-1)}, \theta_{\sigma^{-1}(n-1)}); (A^{(0)}, M^{(0)}),\\
(\asqrt^{(\sigma^{-1}(1))}, A^{(\sigma^{-1}(1))}, M^{(\sigma^{-1}(1))}),
...,(\asqrt^{(\sigma^{-1}(n-1))}, A^{(\sigma^{-1}(n-1))},
M^{(\sigma^{-1}(n-1))})\\ 
(\asqrt^{(n)}, A^{(n)}, M^{(n)})\bigr) .  
\end{multline*}

To extend this to a left action of $S_n$ on $SK(n)$, we first note that
$S_n$ is generated by the symmetric group on the first $n-1$ letters
$S_{n - 1}$ and the transposition $(n-1 \; n)$.  We can let $(n-1 \;
n)$ act on $SK(n)$ by permuting the $(n- 1)$-st and $n$-th punctures
and their local coordinates for a canonical supersphere with $1 +
n$ tubes but the resulting supersphere with $1+n$ tubes is not
canonical.  To obtain the superconformally equivalent canonical
supersphere, we have to translate the new $n$-th puncture to 0.  This
translation will not change the local coordinates at positively
oriented punctures but will change the local coordinates at the
negatively oriented puncture $\infty$.  In \cite{B memoirs} we show that
the resulting canonical supersphere with tubes is given by  
\begin{eqnarray*}
& & \hspace{-.4in} (n-1 \; n) \cdot Q \\
&=& \! \! \Bigl( \infty,(z_1, \theta_1),...,(z_{n-2}, \theta_{n-2}), 0, 
(z_{n-1},\theta_{n-1}); (A^{(0)},M^{(0)}), (\asqrt^{(1)}, A^{(1)},
M^{(1)}), \Bigr.\\
& & \hspace{.15in} \Bigl. ..., (\asqrt^{(n-2)}, A^{(n-2)}, M^{(n-2)}), 
(\asqrt^{(n)}, A^{(n)}, M^{(n)}), (\asqrt^{(n-1)}, A^{(n-1)}, M^{(n-1)})
\Bigr) \\  
&=& \! \! \Bigl(\infty,
s_{(z_{n-1},\theta_{n-1})}(z_1,\theta_1),s_{(z_{n-1},\theta_{n-1})}
(z_2,\theta_2),..., s_{(z_{n-1},\theta_{n-1})}(z_{n-2},
\theta_{n-2}),\Bigr.\\
& & s_{(z_{n-1},\theta_{n-1})}(0),0 ; (\tilde{A}^{(0)},
\tilde{M}^{(0)}), (\asqrt^{(1)}, A^{(1)}, M^{(1)}), ..., 
(\asqrt^{(n-2)}, A^{(n-2)}, M^{(n-2)}),\\ 
& & \hspace{1.3in} (\asqrt^{(n)}, A^{(n)}, M^{(n)}), (\asqrt^{(n-1)},
A^{(n-1)}, M^{(n-1)}) \Bigr) \in SK(n)  
\end{eqnarray*}  
where $(\tilde{A}^{(0)},\tilde{M}^{(0)}) \in \bigwedge_\infty^\infty$ is 
given by
\begin{eqnarray}
& & \hspace{-.2in} \tilde{E}(\tilde{A}^{(0)}, -i \tilde{M}^{(0)}) 
\Bigl(\frac{1}{w},\frac{i\rho}{w}\Bigr) \label{Sn action}\\
\hspace{.3in} &=& \exp \Biggl(\sum_{j \in \Z} \left(\tilde{A}^{(0)}_j
L_{-j}(w,\rho) + \tilde{M}^{(0)}_{j - \frac{1}{2}} G_{-j + \frac{1}{2}}
(w,\rho) \right) \!\Biggr) \! \cdot \! \Bigl(\frac{1}{w},\frac{i\rho}{w}
\Bigr)  \nonumber\\
&=&  \exp \Biggl(\sum_{j \in \Z} \left(A^{(0)}_j L_{-j}(x,\varphi) + 
M^{(0)}_{j - \frac{1}{2}} G_{-j + \frac{1}{2}}(x,\varphi)
\right) \! \Biggr) \! \cdot \nonumber \\
& & \biggl. \hspace{2.2in} \cdot \Bigl(\frac{1}{x},\frac{i\varphi}{x}\Bigr) 
\biggr|_{(x,\varphi) = s_{(-z_{n-1}, -\theta_{n-1})}(w,\rho)} \nonumber \\
&=& \exp \left(-z_{n-1} L_{-1}(w,\rho) -
\theta_{n-1} G_{-\frac{1}{2}}(w,\rho) \right) \cdot \nonumber\\
& & \hspace{.6in} \cdot \exp \Biggl(\sum_{j \in \Z} \left(A^{(0)}_j L_{-j}(w,\rho) + 
M^{(0)}_{j - \frac{1}{2}} G_{-j + \frac{1}{2}}(w,\rho)
\right)\Biggr) \! \cdot \! \Bigl(\frac{1}{w},\frac{i\rho}{w}\Bigr) . \nonumber
\end{eqnarray}

\subsection{Supermeromorphic superfunctions on $SK$ and
supermeromorphic tangent spaces of $SK$}

A {\it supermeromorphic superfunction on $SK(n)$}, for $n \in \Z$, is a
superfunction $F : SK(n) \rightarrow \bigwedge_\infty$ of the form
\begin{eqnarray}
\hspace{.4in} F(Q) \! \! &=& \! \! F\bigl((z_1, \theta_1),...,(z_{n-1}, 
\theta_{n-1}); (A^{(0)},M^{(0)}), (\asqrt^{(1)}, A^{(1)}, M^{(1)}), ..., 
\label{meromorphic} \bigr.\\
& & \bigl. \hspace{2.6in}(\asqrt^{(n)}, A^{(n)}, M^{(n)})\bigr) \nonumber \\
&=& \! \! F_0\bigl( (z_1,\theta_1),...,(z_{n-1}, \theta_{n-1}); (A^{(0)},M^{(0)}), 
(\asqrt^{(1)}, A^{(1)}, M^{(1)}),..., \bigr.\nonumber \\
& & \bigl. (\asqrt^{(n)},A^{(n)}, M^{(n)})\bigr) \times \biggl(\prod^{n-1}_{i
= 1} z_i^{-s_i} \prod_{1 \leq i<j \leq n-1} (z_i - z_j - \theta_i \theta_j)^{-s_{ij}}\biggr)
\nonumber 
\end{eqnarray} 
where $s_i$ and $s_{ij}$ are nonnegative integers and 
\[F_0\bigl((z_1,\theta_1),...,(z_{n-1}, \theta_{n-1}); (A^{(0)},M^{(0)}),
(\asqrt^{(1)}, A^{(1)}, M^{(1)}),...,(\asqrt^{(n)},A^{(n)}, M^{(n)})\bigr)\]  
is a polynomial in the $z_i$'s, $\theta_i$'s, $\asqrt^{(i)}$'s,
$(\asqrt^{(i)})^{-1}$'s, $A^{(i)}_j$'s, and $M^{(i)}_{j - 1/2}$'s.
For $n=0$ a {\it supermeromorphic superfunction on $SK(0)$} is a
polynomial in the components of elements of $SK(0)$, i.e., a polynomial
in the $A^{(0)}_j$'s, and $M^{(0)}_{j - 1/2}$'s.  For $F$ of the form
(\ref{meromorphic}), we say that $F$ has a pole of order $s_{ij}$ at 
$(z_i,\theta_i) = (z_j,\theta_j)$.

The set of all canonical supermeromorphic superfunctions on $SK(n)$ is
denoted by $SD(n)$.  Note that $SD(n)$ has a natural $\mathbb{Z}_2$-grading
since any supermeromorphic superfunction $F$ can be decomposed into an
even component $F^0$ and an odd component $F^1$ where $F^0(Q) \in
\bigwedge_\infty^0$, for all $Q \in SK$ and $F^1(Q) \in
\bigwedge_\infty^1$ for all $Q \in SK$. 

An {\it even supermeromorphic tangent vector} of $SK$ at 
a point $Q \in SK$ is a linear map $\mathcal{L}_Q : SD(n) \longrightarrow
\bigwedge_\infty$ such that for $F_1, F_2 \in SD(n)$, the map 
$\mathcal{L}_Q$ satisfies $\mathcal{L}_Q (F_1 F_2) = \mathcal{L}_Q F_1
\cdot F_2(Q) +  F_1(Q) \cdot \mathcal{L}_Q F_2$, and an {\it odd
supermeromorphic tangent vector} of $SK$ at a point $Q \in SK$ is a linear
map $\mathcal{L}_Q : SD(n) \longrightarrow \bigwedge_\infty$ such that
for $F_1$ of homogeneous sign in $SD(n)$, we have $\mathcal{L}_Q (F_1
F_2) = \mathcal{L}_Q F_1 \cdot F_2(Q) + (-1)^{ \eta(F_1)} F_1(Q) \cdot
\mathcal{L}_Q F_2$.   The set of all supermeromorphic tangent vectors at 
$Q$ is the {\it supermeromorphic tangent space of $SK$ at $Q$} and is
denoted $T_QSK$ (or $T_QSK(n)$ when $Q \in SK(n)$).  Since $T_QSK$ is
naturally $\mathbb{Z}_2$-graded, the sign function $\eta$ is well defined
on $T_QSK$. In \cite{B thesis} and \cite{B memoirs}, we prove that any
supermeromorphic tangent vector $\mathcal{L}_Q$ can be expressed as 
\begin{multline}\label{tangent}
\mathcal{L}_Q = \sum_{i = 1}^{n-1} \biggl( \Bigl. \alpha_i
\frac{\partial}{\partial z_i} \Bigr|_Q + \beta_i \Bigl.
\frac{\partial}{\partial \theta_i} \Bigr|_Q \biggr) + \sum_{i =
1}^n \delta_i \biggl. \frac{\partial}{\partial \asqrt^{(i)}}
\biggr|_Q \\
+ \; \sum_{i = 0}^n \sum_{j \in \Z}
\Biggl(\gamma_i^j \biggl. \frac{\partial}{\partial A_j^{(i)}} \biggr|_Q
+ \nu_i^j \biggl. \frac{\partial}{\partial M_{j - \frac{1}{2}}^{(i)}}
\biggr|_Q \Biggr) 
\end{multline}  
where
\[\alpha_i = \mathcal{L}_Q z_i, \quad \beta_i = \mathcal{L}_Q \theta_i,
\quad \delta_i = \mathcal{L}_Q \asqrt^{(i)}, \quad \gamma_i^j = \mathcal{L}_Q
A_j^{(i)}, \quad \nu_i^j = \mathcal{L}_Q M_{j - \frac{1}{2}}^{(i)} . \]
Of course if $\eta(\mathcal{L}_Q) = 0$, then $\alpha_i, \delta_i, \gamma_i^j
\in \bigwedge_\infty^{\eta(\mathcal{L}_Q)}$ and $\beta_i, \nu_i^j \in
\bigwedge_\infty^{(\eta(\mathcal{L}_Q) + 1) \; \mathrm{mod} \; 2}$.

\subsection{The sewing operation}

In \cite{B memoirs} we introduce the sewing operation for superspheres with
tubes.  This gives rise to the sewing equation, normalization conditions and
boundary conditions derived {}from taking two canonical superspheres with
tubes, sewing them together and determining conditions on the uniformizing
function which maps the resulting supersphere to a super-Riemann sphere. 
This resulting super-Riemann sphere is then superconformal to a canonical 
supersphere, i.e., an element in $SK$, via a superprojective  transformation.
Thus for any $m \in \Z$, $n \in \mathbb{N}$, and any  positive integer $i
\leq m$, the sewing operation for superspheres  with tubes defined in
\cite{B memoirs} induces an operation 
\begin{eqnarray*}
 _i\infty_0 : SK(m) \times SK(n) &\rightarrow& SK(m + n -  1)\\
(Q_1, Q_2) &\mapsto& Q_1 \; _i\infty_0 \; Q_2 .  
\end{eqnarray*}

For $m \in \Z$, let $Q_1 \in SK(m)$ be given by
\[\bigl((z_1, \theta_1),...,(z_{m-1}, \theta_{m-1});
(A^{(0)},M^{(0)}), (\asqrt^{(1)}, A^{(1)}, M^{(1)}),...,(\asqrt^{(m)},
A^{(m)}, M^{(m)})\bigr) \]
and for $n \in \mathbb{N}$, let $Q_2 \in SK(n)$ be given by  
\[\bigl((z_1', \theta_1'),...,(z_{n-1}', \theta_{n-1}');
(B^{(0)},N^{(0)}), (\bsqrt^{(1)}, B^{(1)}, N^{(1)}),...,(\bsqrt^{(n)},
B^{(n)}, N^{(n)})\bigr) .\]
For convenience, we will sometimes denote the puncture at $0$ of $Q_1$
by $(z_m, \theta_m)$.  Let
\begin{eqnarray*}
H^{(1)} (w, \rho)\!  &=& \! \exp \Biggl(\!  - \! \sum_{j \in \Z} \biggl( A^{(i)}_j
L_j(w,\rho)  +  M^{(i)}_{j - \frac{1}{2}} G_{j - \frac{1}{2}}(w,\rho)
\biggr) \! \Biggr) \cdot \\
& & \hspace{2.4in} \cdot (\asqrt^{(i)})^{-2L_0(w,\rho)}\cdot (w, \rho ) \\
&=& \hat{E}(\asqrt^{(i)}, A^{(i)}, M^{(i)}) (w,\rho), \\
H^{(2)} (w, \rho) \! &=& \! \exp \Biggl(\sum_{j \in \Z} \biggl( B^{(0)}_j
L_{-j}(w,\rho) +  N^{(0)}_{j - \frac{1}{2}} G_{-j + \frac{1}{2}}(w,\rho)
\biggr) \! \Biggr) \cdot \Bigl(\frac{1}{w}, \frac{i \rho}{w} \Bigr) \\
&=& \tilde{E}(B^{(0)},-i N^{(0)}) \Bigl(\frac{1}{w}, \frac{i \rho}{w} \Bigr). 
\end{eqnarray*}
Then the local coordinate vanishing at the $i$-th puncture of the canonical
supersphere represented by $Q_1$ is $H^{(1)} \circ s_{(z_i,\theta_i)}(w,\rho)$,
and the local coordinate vanishing at the puncture at $\infty$ of the 
canonical supersphere represented by $Q_1$ is $H^{(2)}(w,\rho)$. 

Recall the standard coordinate atlas for the super-Riemann sphere, 
$S\hat{\mathbb{C}}$, given by $\{(U_\sou,\sou),(U_\nor,\nor)\}$ with coordinate
transition given by $I = \sou \circ \nor^{-1}$.  For a superconformal
equivalence $F : S\hat{\mathbb{C}} \rightarrow S\hat{\mathbb{C}}$ define
$F_\sou = \sou \circ F \circ \sou^{-1}$ and $F_\nor = \nor \circ F \circ
\nor^{-1}$.   Also define the DeWitt open (respectively, closed) discs
about zero and of radius $r >0$,  by
\[ \mathcal{B}_0^r = \{(w,\rho) \in \mbox{$\bigwedge_\infty$} \; | \;
|w_B| < r \} \quad (\mbox{resp.,} \; \;  \bar{\mathcal{B}}_0^r =
\{(w,\rho) \in \mbox{$\bigwedge_\infty$} \; | \; |w_B| \leq r \}) .\]

The following theorem summarizes several results {}from \cite{B thesis}
and \cite{B memoirs}, and provides exact formulas for calculating the
resulting canonical supersphere {}from the sewing together of two
canonical superspheres.

\begin{thm}\label{actual sewing}(\cite{B memoirs})
The i-th tube of $Q_1$ can be sewn with the 0-th tube of $Q_2$ if and
only if there exist $r_1, r_2 \in \mathbb{R}_+$, with $r_1 > r_2$ such 
that the series $(H^{(1)})^{-1}(w,\rho)$ and $(H^{(2)})^{-1}(w,\rho)$
are absolutely convergent in $\mathcal{B}_0^{r_1}$ and
$\mathcal{B}_0^{1/r_2} \smallsetminus (\{0 \} \times (\bigwedge_\infty)_S)$,
respectively, $(-z_i, -\theta_i), (z_k - z_i - \theta_k \theta_i, 
\theta_k - \theta_i) \notin (H^{(1)})^{-1} (\mathcal{B}_0^{r_1})$,
for $k = 1,..., m - 1$, $k \neq i$, and $0, (z_l', \theta_l') \notin
(H^{(2)})^{-1}(\mathcal{B}_0^{1/r_2} \smallsetminus (\{0 \} \times
(\mbox{$\bigwedge_\infty$})_S))$, for $l = 1,...,n-1$.  In
this case, there exist unique  bijective superconformal functions
$F^{(1)} (w,\rho)$ and $F^{(2)} (w,\rho)$ defined on 
\[S\hat{\mathbb{C}} \smallsetminus \sou^{-1} \circ s_{(z_i,\theta_i)}^{-1} 
\circ (H^{(1)})^{-1} (\bar{\mathcal{B}}_0^{r_2})\] 
and 
\[U_\sou \smallsetminus \sou^{-1} \circ (H^{(2)})^{-1}
(\bar{\mathcal{B}}_0^{1/r_1} \smallsetminus (\{0\} \times 
\mbox{$(\bigwedge_\infty)_S$})) \subset S\hat{\mathbb{C}},\]  
respectively, satisfying the normalization conditions:
\begin{eqnarray}
F^{(1)}_\nor  (0) &=& 0, \label{normalize 1} \\
\lim_{w \rightarrow \infty} \frac{\partial}{\partial \rho} (F^{(1)}_\sou)^1
(w,\rho) &=& 1 \label{normalize 2} \\
F^{(2)}_\sou (0) &=& 0 \label{normalize 3} 
\end{eqnarray}
and such that in $(H^{(1)} \circ s_{(z_i, \theta_i)})^{-1} (\mathcal{B}_0^{r_1} \smallsetminus
\bar{\mathcal{B}}_0^{r_2})$, we have 
\begin{equation}\label{F relation}
F^{(1)}_\sou  (w, \rho) = F^{(2)}_\sou  \circ (H^{(2)})^{-1} \circ I \circ H^{(1)}
\circ s_{(z_i,\theta_i)} (w, \rho) . 
\end{equation}
Furthermore, if we let 
\[H_0^{(1)}(w,\rho) = \tilde{E}(A^{(0)}, -iM^{(0)}) \Bigl(\frac{1}{w},
\frac{i\rho}{w}\Bigr) ,\]
denote the local coordinate of the puncture vanishing at $\infty$ of the
canonical supersphere represented by $Q_1$; let
\[H_k^{(1)}(w,\rho) = \hat{E}(\asqrt^{(k)}, A^{(k)}, M^{(k)}) (w, \rho), \]
so that the local coordinate vanishing at the $k$-th puncture of the 
canonical supersphere represented by $Q_1$ is given by $H_k^{(1)} \circ 
s_{(z_k,\theta_k)}$, for $k = 1,...,m$, $k \neq i$; and let 
\[H_l^{(2)}(w,\rho) = \hat{E}(\bsqrt^{(l)}, B^{(l)}, N^{(l)}) (w, \rho), \] 
so that the local coordinate vanishing at the $l$-th puncture of the 
canonical supersphere represented by $Q_2$ is given by $H_l^{(2)} \circ 
s_{(z_l',\theta_l')}$, for $l = 1,...,n$; then we have the following:

(1) When $i = m$, and $n > 0$, the punctures of the canonical supersphere
represented by $Q_1 \; _m\infty_0 \; Q_2$ are 
\[\infty, \; F^{(1)}_\sou  (z_1, \theta_1),..., F^{(1)}_\sou (z_{m-1},
\theta_{m-1}), \;  F^{(2)}_\sou (z_1', \theta_1'),...,
 F^{(2)}_\sou (z_{n-1}', \theta_{n-1}'), \; 0;\]
and the local coordinates vanishing at these punctures are
\[H_0^{(1)} \circ (F^{(1)}_\sou)^{-1} (w, \rho),\]    
\[H_1^{(1)} \circ s_{(z_1,\theta_1)} \circ (F^{(1)}_\sou)^{-1} (w,\rho), 
..., H_{m-1}^{(1)} \circ s_{(z_{m-1},\theta_{m-1})} \circ (F^{(1)}_\sou)^{-1}
(w,\rho), \] 
\[H_1^{(2)} \circ s_{(z_1',\theta_1')} \circ (F^{(2)}_\sou)^{-1}
(w,\rho), ..., H_{n-1}^{(2)} \circ s_{(z_{n-1}',\theta_{n-1}')} \circ 
(F^{(2)}_\sou)^{-1} (w,\rho) \]
\[H_{n}^{(2)} \circ  (F^{(2)}_\sou)^{-1}(w,\rho) , \]
respectively.

(2)  When $i = m = 1$ and $n = 0$, the canonical supersphere represented
by $Q_1 \; _1\infty_0 \; Q_2$ has only the one puncture at $\infty$ with 
local coordinate given by
\begin{equation}\label{coordinate with only one}
H_0^{(1)} \circ (F^{(1)}_\sou)^{-1} \circ s_{(a',m')}^{-1} (w,\rho),
\end{equation}
where $(a',m') \in \bigwedge_\infty$ is the unique element such that
$Q_1 \; _1\infty_0 \; Q_2$ represents a canonical supersphere in $SK(0)$,
i.e., such that the expansion of (\ref{coordinate with only one}) has
even coefficient of $w^{-2}$ and odd coefficient of $w^{-1}$ equal to zero.
Thus
\begin{multline*}
Q_1 \; _1\infty_0 \; Q_2 = \bigl( (\tilde{E}^{-1})^0 (H_0^{(1)} \circ 
(F^{(1)}_\sou)^{-1} \circ s_{(a',m')}^{-1} \circ I^{-1} (w,\rho)),  \\
i(\tilde{E}^{-1})^1 (H_0^{(1)} \circ
(F^{(1)}_\sou)^{-1} \circ s_{(a',m')}^{-1} \circ I^{-1} (w,\rho)) \bigr). 
\end{multline*}

(3)  When $i = m > 1$, and $n = 0$, and writing $F^{(1)}_\sou (z_{m - 1},
\theta_{m-1}) = p \in \bigwedge_\infty$, the punctures of the canonical 
supersphere represented by $Q_1 \; _i\infty_0 \; Q_2$ are
\[\infty, \; s_p \circ F^{(1)}_\sou (z_1, \theta_1),..., s_p \circ F^{(1)}_\sou 
(z_{m-2}, \theta_{m-2}), \; 0; \]
and the local coordinates vanishing at these punctures are
\[H_0^{(1)} \circ (F^{(1)}_\sou)^{-1} \circ  s_p^{-1} (w, \rho), \]
\[H_1^{(1)} \circ s_{(z_1,\theta_1)} \circ (F^{(1)}_\sou)^{-1} \circ 
s_p^{-1} (w,\rho), ..., H_{m-2}^{(1)} \circ s_{(z_{m-2},\theta_{m-2})} 
\circ (F^{(1)}_\sou)^{-1}\circ s_p^{-1} (w,\rho)\] 
\[H_{m-1}^{(1)} \circ s_{(z_{m-1},\theta_{m-1})} \circ (F^{(1)}_\sou)^{-1}
\circ s_p^{-1} (w,\rho),\] 
respectively. 

(4) When $i < m$ and $n \neq 0$, writing $F^{(1)}_\sou (0) = p \in 
\bigwedge_\infty$, the punctures of the canonical supersphere represented 
by $Q_1 \; _i\infty_0 \; Q_2$ are
\begin{multline*}
\infty, \; s_p \circ F^{(1)}_\sou (z_1,\theta_1),..., s_p \circ 
F^{(1)}_\sou (z_{i-1}, \theta_{i-1}), \; s_p \circ F^{(2)}_\sou (z_1',\theta_1'), 
...,\\
s_p \circ F^{(2)}_\sou (z_{n-1}', \theta_{n-1}'), 
- F^{(1)}_\sou (0), \; s_p \circ F^{(1)}_\sou (z_{i+1}, \theta_{i+1}), ...,\\
s_p \circ F^{(1)}_\sou (z_{m-1}, \theta_{m-1}), \; 0 ; 
\end{multline*}
and the local coordinates vanishing at these punctures are
\[H_0^{(1)} \circ (F^{(1)}_\sou)^{-1} \circ s_p^{-1} (w, \rho), \]
\[H_1^{(1)} \circ s_{(z_1,\theta_1)} \circ (F^{(1)}_\sou)^{-1} \circ 
s_p^{-1} (w,\rho),..., H_{i-1}^{(1)} \circ s_{(z_{i-1},\theta_{i-1})} 
\circ (F^{(1)}_\sou)^{-1} \circ s_p^{-1} (w,\rho), \]    
\[H_1^{(2)} \circ s_{(z_1',\theta_1')} \circ (F^{(2)}_\sou)^{-1} \circ
s_p^{-1} (w,\rho), ...,  H_{n-1}^{(2)} \circ s_{(z_{n-1}',\theta_{n-1}')} 
\circ  (F^{(2)}_\sou)^{-1} \circ s_p^{-1} (w,\rho),\]  
\[H_n^{(2)} \circ (F^{(2)}_\sou)^{-1} \circ s_p^{-1} (w,\rho), \] 
\begin{multline*}
H_{i+1}^{(1)} \circ s_{(z_{i+1},\theta_{i+1})} \circ
(F^{(1)}_\sou)^{-1}  \circ s_p^{-1} (w,\rho),...,  \\
 H_{m-1}^{(1)} \circ s_{(z_{m-1},\theta_{m-1})} 
\circ (F^{(1)}_\sou)^{-1} \circ s_p^{-1} (w,\rho), 
\end{multline*} 
\[H_m^{(1)} \circ (F^{(1)}_\sou)^{-1}\circ s_p^{-1} (w,\rho),\]
respectively. 

(5) When $i < m$ and $n = 0$, writing $F^{(1)}_\sou (0) = p \in 
\bigwedge_\infty$, the punctures of the canonical supersphere represented 
by $Q_1 \; _i\infty_0 \; Q_2$ are
\begin{multline*}
\infty, \; s_p \circ F^{(1)}_\sou (z_1,\theta_1),..., s_p \circ 
F^{(1)}_\sou (z_{i-1}, \theta_{i-1}), \; s_p \circ F^{(1)}_\sou (z_{i+1}, 
\theta_{i+1}), ..., \\
s_p \circ F^{(1)}_\sou (z_{m-1}, \theta_{m-1}), \; 0 ; 
\end{multline*}
and the local coordinates vanishing at these punctures are
\[H_0^{(1)} \circ (F^{(1)}_\sou)^{-1} \circ s_p^{-1} (w, \rho), \]
\[H_1^{(1)} \circ s_{(z_1,\theta_1)} \circ (F^{(1)}_\sou)^{-1} \circ 
s_p^{-1} (w,\rho),..., H_{i-1}^{(1)} \circ s_{(z_{i-1},\theta_{i-1})} 
\circ (F^{(1)}_\sou)^{-1} \circ s_p^{-1} (w,\rho), \]    
\begin{multline*}
H_{i+1}^{(1)} \circ s_{(z_{i+1},\theta_{i+1})} \circ (F^{(1)}_\sou)^{-1} 
\circ s_p^{-1} (w,\rho),...,  \\
H_{m-1}^{(1)} \circ s_{(z_{m-1},\theta_{m-1})} 
\circ (F^{(1)}_\sou)^{-1} \circ s_p^{-1} (w,\rho), 
\end{multline*}
\[H_m^{(1)} \circ (F^{(1)}_\sou)^{-1}\circ s_p^{-1} (w,\rho), \]
respectively. 

Moreover, the components $F_\sou^{(1)}$ and $F_\sou^{(2)}$ of the
uniformizing function  can be expanded in formal power series such that 
\begin{eqnarray} \label{Fs1}
\hspace{.4in} F_\sou^{(1)} (x,\varphi) &=& \left. \bar{F}^{(1)}
\circ s_{(z_i,\theta_i)} (x,\varphi)
\right|_{(\alpha_0^{1/2}, \mathcal{A}, \mathcal{M}, \mathcal{B}, \mathcal{N}) = 
(\asqrt^{(i)},A^{(i)},M^{(i)},B^{(0)},N^{(0)})}\\
F_\sou^{(2)} (x,\varphi) &=& \left. \bar{F}^{(2)} (x,\varphi)
\right|_{(\alpha_0^{1/2}, \mathcal{A}, \mathcal{M}, \mathcal{B}, \mathcal{N}) = 
(\asqrt^{(i)},A^{(i)},M^{(i)},B^{(0)},N^{(0)})} , \label{Fs2}
\end{eqnarray}
where
\begin{eqnarray}
\ \ \ \ \ \ \ \ \bar{F}^{(1)} (x,\varphi) &= & \exp \Biggl(\sum_{j \in
\Z} \biggl( \Psi_{-j} L_{-j}(x,\varphi)  + \Psi_{-j+ \frac{1}{2}}
G_{-j + \frac{1}{2}}(x,\varphi)
\biggr) \! \Biggr) \cdot (x,\varphi) \label{F1}\\
\bar{F}^{(2)} (x,\varphi) &=& \exp \left( - \Psi_0 2L_0(x,\varphi) \right) 
\cdot (\asqrt^{(i)})^{ 2L_0 (x,\varphi)} \cdot \nonumber\\
& & \qquad \exp \Biggl(\! - \! 
\sum_{j \in \Z} \biggl( \Psi_j L_j(x,\varphi) + \Psi_{j-\frac{1}{2}}
G_{j - \frac{1}{2}}(x,\varphi) \biggr) \! \Biggr) \cdot (x,\varphi) \label{F2} 
\end{eqnarray} 
are formal power series in 
\[x \mathbb{C} [\alpha_0^{\frac{1}{2}}, \alpha_0^{-\frac{1}{2}}]
[[\mathcal{A}, \mathcal{B}]] [\mathcal{M}, \mathcal{N}] [[x^{-1}]] [\varphi] \quad
\mbox{and} \quad x\mathbb{C}[\alpha_0^{\frac{1}{2}},
\alpha_0^{-\frac{1}{2}}] [[\mathcal{A}, \mathcal{B}]] [\mathcal{M}, \mathcal{N}]
[[x]] [\varphi], \]
respectively, and where 
\begin{equation}\label{where psi lives}
(\Psi_j, \Psi_{j - \frac{1}{2}}) = (\Psi_j, \Psi_{j - \frac{1}{2}})
(\alpha_0^\frac{1}{2}, \mathcal{A}, \mathcal{M}, \mathcal{B}, \mathcal{N}) 
\end{equation} 
are a unique pair of sequences in $\mathbb{C} [\alpha_0^{1/2}, \alpha_0^{-1/2}][[\mathcal{A},
\mathcal{B}]] [\mathcal{M}, \mathcal{N}]$, for $j \in \mathbb{Z}$, such that
\begin{eqnarray}
(\Psi_j, \Psi_{j - \frac{1}{2}}) \! \! &=& \! \! (- \mathcal{A}_j , - \mathcal{M}_{j -
\frac{1}{2}}) + \mathcal{P}_j (\alpha_0^\frac{1}{2}, \mathcal{A}, \mathcal{M},
\mathcal{B}, \mathcal{N}) , \label{psi condition 1} \\  
\hspace{.2in} (\Psi_{-j}, \Psi_{- j + \frac{1}{2}}) \! \! &=& \! \! (- \alpha_0^{-j} 
\mathcal{B}_j, - \alpha_0^{-j + \frac{1}{2}}  \mathcal{N}_{j - \frac{1}{2}}) + 
\mathcal{P}_{-j} (\alpha_0^\frac{1}{2}, \mathcal{A}, \mathcal{M}, \mathcal{B}, 
\mathcal{N}),  \label{psi condition 2} \\  
\Psi_0 \! \! &=& \! \! 0 + \mathcal{P}_0 (\alpha_0^\frac{1}{2}, \mathcal{A}, \mathcal{M},
\mathcal{B}, \mathcal{N}) , \label{psi condition 3}  
\end{eqnarray}
each $\mathcal{P}_j (\alpha_0^{1/2}, \mathcal{A}, \mathcal{M},
\mathcal{B}, \mathcal{N})$, for $j \in \mathbb{Z}$, contains only terms with 
total degree at least one in the $\mathcal{A}_k$'s and $\mathcal{M}_{k - 1/2}$'s, 
for $k \in \Z$, and with total degree at least one in the $\mathcal{B}_k$'s and 
$\mathcal{N}_{k - 1/2}$'s, for $k \in \Z$.

In particular, letting $t^{1/2}$ be a complex variable, and letting
\[\Psi_j(t^\frac{1}{2}) = \Psi_j(t^{-\frac{1}{2}}\asqrt^{(i)}, A^{(i)},
M^{(i)}, B^{(0)}, N^{(0)}),\]  
for $j \in \frac{1}{2} \mathbb{Z}$.  Then $\Psi_j (t^{1/2}) \in
\bigwedge_\infty [[t^{-1/2}, t^{1/2}]]$ and if $Q_1 \; _i\infty_0 \; Q_2$
exists, then the series $\Psi_j(t^{1/2})$, for $j \in \frac{1}{2} 
\mathbb{Z}$, are convergent when $|t^{1/2}| \leq 1$, and the values  of
these convergent series are equal to the coefficients $\Psi_j$ of the 
power series (\ref{Fs1}) and (\ref{Fs2}).

Finally, using this solution of the uniformizing function $F$ given by
$F^{(1)}$ and $F^{(2)}$, we obtain the canonical element of $SK(m + n -1)$
representing the sewn spheres by taking $\hat{E}^{-1}$ of each of the
local coordinates given in (1) - (5).
\end{thm}

\begin{rema}\label{operad remark} 
{\em The moduli space of $N=1$ superspheres with tubes, $SK$, along 
with the sewing operation and the action of the symmetric group defined 
in Section \ref{symmetric group}, is an example of a partial operad (cf.
\cite{M}, \cite{HL1}, \cite{HL2}, \cite{H book}).}
\end{rema}

The element of $\mathcal{H}$ with all components equal to $0$ will be 
denoted by $(\mathbf{0},\mathbf{0})$ or just $\mathbf{0}$.  Note that 
the local coordinate chart corresponding to $(1,\mathbf{0},\mathbf{0}) 
\in (\bigwedge_\infty^0)^\times \times \mathcal{H}$ is the identity 
map on $\bigwedge_\infty$ if the puncture is at $0$; the shift 
$s_{(z_i,\theta_i)}(w,\rho) =  (w - z_i -\rho\theta_i,\rho - \theta_i)$ 
if the puncture is at $(z_i,\theta_i)$; and $I(w,\rho) = (1/w,
i\rho/w)$ if the puncture is at $\infty$.  We will sometimes refer to 
such coordinates as {\it standard local coordinates}. 
{}From Theorem \ref{actual sewing}, we see that $SK$ has a 
unit under the sewing operation given by $e = (\mathbf{0}, (1,
\mathbf{0})) \in SK(1)$, i.e., for $Q \in SK(n)$ and $0<i\leq n$, the
$i$-th puncture of $Q$ can always be sewn with the $0$-th puncture of
$e$, the first puncture of $e$ can always be sewn with the $0$-th
puncture of $Q$, and we have $Q \; _i\infty_0 \; e = Q = e \;
_1\infty_0 \; Q$.   {}From the geometry of sewing defined in \cite{B
memoirs}, the following  associativity of the sewing operation is obvious.

\begin{prop}\label{sewing associativity}(\cite{B memoirs})
Let $l \in \Z$ and $m, n \in \mathbb{N}$ such that $l + m -1 \in \Z$, and
let $Q_1 \in SK(l)$, $Q_2 \in SK(m)$, $Q_3 \in SK(n)$, and $i,j \in \Z$
such that $1 \leq i \leq l$, and $1 \leq j \leq l + m -1$.  The iterated
sewings $(Q_1 \; _i\infty_0 \; Q_2) \; _j\infty_0 \; Q_3$ exist if and
only if one of the following holds:

(i) $j<i$ and the sewings $(Q_1 \; _j\infty_0 \; Q_3) \; _{i + n -
1}\infty_0 \; Q_2$ exist, in which case
\[(Q_1 \; _i\infty_0 \; Q_2) \; _j\infty_0 \; Q_3 = (Q_1 \; _j\infty_0
\; Q_3) \; _{i + n - 1}\infty_0 \; Q_2 ;\] 

(ii) $j \geq i + m$ and the sewings $(Q_1 \; _{j - m + 1}\infty_0
\; Q_3) \; _i\infty_0 \; Q_2$ exist, in which case
\[(Q_1 \; _i\infty_0 \; Q_2) \; _j\infty_0 \; Q_3 = (Q_1 \; _{j - m +
1}\infty_0 \; Q_3) \; _i\infty_0 \; Q_2 ;\] 

or

(iii) $i \leq j < i + m$ and the sewings $Q_1 \; _i\infty_0 \;
(Q_2 \; _{j - i + 1}\infty_0 \; Q_3)$ exist, in which case 
\[(Q_1 \; _i\infty_0 \; Q_2) \; _j\infty_0 \; Q_3 = Q_1 \; _i\infty_0
\; (Q_2 \; _{j - i + 1}\infty_0 \; Q_3) .\]
\end{prop}

\subsection{Superspheres with one, two, and three tubes}

\begin{prop}\label{how to get any supersphere}(\cite{B memoirs})
Any element $Q \in SK$ can be obtained by sewing the following types of
elements of $SK$ \\
(i) $(\mathbf{0}) \in SK(0)$, \\
(ii) $ \bigl((A^{(0)},M^{(0)}), (\asqrt^{(1)}, A^{(1)},M^{(1)})\bigr) \in SK(1)$, \\
(iii) $\bigl((z, \theta);\mathbf{0}, (1, \mathbf{0}), (1, \mathbf{0})\bigr) \in
SK(2)$. \\   
\end{prop}

It is clear {}from the definition of sewing that $SK(1)$ is a partial monoid.
In the following proposition, we give some subgroups of $SK(1)$.

\begin{prop}\label{t and s composition}(\cite{B memoirs})
Let $s,t \in \bigwedge_\infty^0$, and $(A,M) \in \bigwedge_\infty^\infty$.
Assume 
\[s(A,M), t(A,M) \in \mathcal{H}.\]  
Then $(s + t)(A,M) \in \mathcal{H}$, both
\[\bigl(\mathbf{0}, (1,t(A,M)\bigr) \; _1\infty_0 \; \bigl(\mathbf{0}, (1,s(A,M)\bigr)\]
and
\[ \bigl(s(A,M),(1, \mbox{\bf 0})\bigr) \; _1\infty_0 \; \bigl(t(A,M),(1, \mathbf{0})\bigr) \]
exist, and we have
\begin{eqnarray}
\bigl(\mathbf{0}, (1,(s + t)(A,M))\bigr) \! &=& \! \bigl(\mathbf{0}, (1,t(A,M)\bigr) \;
_1\infty_0 \; \bigl(\mathbf{0}, (1,s(A,M)\bigr), \label{sequencesew1} \\
\bigl((s + t)(A,M), (1, \mathbf{0})\bigr) \! &=& \! \bigl(s(A,M),(1, \mbox{\bf
0})\bigr) \; _1\infty_0 \; \bigl(t(A,M),(1, \mathbf{0})\bigr) . \label{sequencesew2} 
\end{eqnarray}
In particular, for $(A,M) \in \mathcal{H}$, the sets 
\begin{equation}\label{first subgroup}
\bigl\{\bigl(\mathbf{0},(1,t(A,M)\bigr) \; | \; t \in \mbox{$\bigwedge_\infty^0$},
\;  t(A,M) \in \mathcal{H} \bigr\}
\end{equation} 
and
\begin{equation}\label{second subgroup}
\bigl\{\bigl(t(A,M), (1,\mathbf{0})\bigr) \; | \; t \in \mbox{$\bigwedge_\infty^0$}, \; 
t(A,M) \in \mathcal{H} \bigr\}
\end{equation} 
are subgroups of $SK(1)$.  In addition we have the subgroups given by
taking $A = \mathbf{0}$ or $M = \mathbf{0}$ in (\ref{first subgroup})
and (\ref{second subgroup}). 
\end{prop}

\begin{prop}\label{for comm and assoc}(\cite{B memoirs})
Let $z_1,z_2 \in \bigwedge_\infty^0$ such that $|(z_1)_B| > |(z_2)_B| >
|(z_1)_B - (z_2)_B| > 0$.  Then
\begin{eqnarray} 
& & \hspace{-.4in} \bigl((z_1,\theta_1), (z_2,\theta_2); 
\mathbf{0}, (1,\mathbf{0}), (1,\mathbf{0}), (1,\mathbf{0})\bigr)  = \nonumber \\
\hspace{.3in} &=& \! \! \bigl((z_2,\theta_2); \mathbf{0}, (1,\mathbf{0}), 
(1,\mathbf{0})\bigr) \; _1\infty_0 \; \bigl((z_1 - z_2 -\theta_1 \theta_2 , \theta_1 - 
\theta_2);  \bigr.\label{associativity1}\\  
& & \bigl.\hspace{3in} \mathbf{0}, (1,\mathbf{0}), (1,\mathbf{0})\bigr) \nonumber \\
&=& \! \! \bigl((z_1,\theta_1); \mathbf{0}, (1,\mathbf{0}), (1,\mathbf{0})\bigr) \;
_2\infty_0 \; \bigl((z_2, \theta_2); \mathbf{0}, (1,\mathbf{0}),
(1,\mathbf{0})\bigr) .\label{associativity2}
\end{eqnarray}
\end{prop}

\begin{rema} {\em Proposition \ref{for comm and assoc} is used in 
Section \ref{alg {}from geom} to prove the associativity property for an
$N=1$ Neveu-Schwarz vertex operator superalgebra with odd formal 
variables obtained {}from an $N=1$ supergeometric vertex operator 
superalgebra.  Thus Proposition \ref{for comm and assoc} can be thought of 
as a geometric version of this algebraic relation.}  
\end{rema}

Let $(a,m) \in \bigwedge_\infty^\infty$, and $k,l \in \Z$, and define
\begin{multline}\label{define tau series}
\bigl(A(a,k), M(m,l - 1/2)\bigr) =  \Bigl( \bigl\{A_j \; | \; 
A_k = a, \; A_j = 0, \; \mbox{for} \; j \neq k \bigr\}_{j \in \Z} , \Bigr. \\
\Bigl. \bigl\{M_{j - \frac{1}{2}} \; | \; M_{l - \frac{1}{2}} = m, \; 
M_{j - \frac{1}{2}} = 0, \; \mbox{for} \; j \neq l \bigr\}_{j \in \Z} \Bigr) 
\end{multline} 
which is an element of $\bigwedge_\infty^\infty$.

\begin{prop}\label{for derivatives}(\cite{B memoirs})
For $(z,\theta), (z_0, \theta_0) \in \bigwedge_\infty$ such that
$0 < |(z_0)_B| < |z_B|$, we have 
\begin{multline*}
\bigl((z_0 + z +  \theta_0 \theta, \theta + \theta_0); \mathbf{0},
(1,\mathbf{0}), (1,\mathbf{0})\bigr) = \\
\bigl((z, \theta); \mathbf{0}, (1,\mathbf{0}), (1,\mathbf{0}) \bigr) \; 
_1\infty_0 \; \bigl((A(-z_0,1), M(-\theta_0,1/2)), (1,\mathbf{0}) 
\bigr) .   
\end{multline*}
\end{prop}

Let $\epsilon \in \bigwedge_\infty^1$. {}From the definition of $SK(0)$,
we see that $(\mathbf{0}, M(\epsilon,3/2)) \in SK(0)$.  Let
$F$ be any element of $SD(1)$, and $(z, \theta) \in \bigwedge_\infty$.
We define a linear functional on $SD(1)$ by 
\begin{equation}\label{define G}
\mathcal{G}_e (z, \theta) F = \Bigl.  \frac{d}{d \epsilon} F
\Bigl( \bigl((z,\theta) ; \mathbf{0}, (1,\mathbf{0}), (1,\mathbf{0})\bigr)
\; _1\infty_0 \; \bigl(\mathbf{0}, M(\epsilon, 3/2) \bigr) \Bigr) 
\Bigr|_{\epsilon = 0} . 
\end{equation}   

\begin{prop}\label{the linear functional G}(\cite{B memoirs})
The linear functional $\mathcal{G}_e (z, \theta)$ is in $T_e SK(1)$, and 
\begin{multline}\label{tau} 
\mathcal{G}_e (z, \theta) = \sum_{k = 0}^1 \sum_{j \in \Z}
z^{-(2k - 1)j - 2 + k} \Biggl. \frac{\partial}{\partial M_{j -
\frac{1}{2}}^{(k)}} \Biggr|_e  \\
 + \; 2 \theta \Biggl(z^{-2} \Biggl. \frac{1}{2}
\frac{\partial}{\partial \asqrt^{(1)}} \Biggr|_e + \sum_{k = 0}^1 
\sum_{j \in \Z} z^{-(2k - 1)j - 2} \Biggl. \frac{\partial}{\partial
A_j^{(k)}} \Biggr|_e \Biggr) .  
\end{multline}
\end{prop}

\subsection{An $N=1$ Neveu-Schwarz algebra structure of central charge 
zero on the supermeromorphic tangent  space of $SK(1)$ at its
identity}\label{tangent section}

In \cite{B memoirs} we define a bracket operation on a subspace $\hat{T}_e
SK(1)$ of the supermeromorphic tangent space $T_e SK(1)$ of $SK(1)$ at the
identity $e$.  This subspace $\hat{T}_e SK(1)$ is the subspace of $T_e
SK(1)$ consisting of all finite linear combinations of $- \left.
\frac{\partial}{\partial a}
\right|_e$, $- \Bigl. \frac{\partial}{\partial A_j^{(0)}}\Bigr|_e$,
$- \Bigl. \frac{\partial}{\partial M_{j - 1/2}^{(0)}}
\Bigr|_e$, $- \Bigl. \frac{\partial}{\partial A_j^{(1)}}\Bigr|_e$,
and $\Bigl. - \frac{\partial}{\partial M_{j - 1/2}^{(1)}}\Bigr|_e$, 
for $j \in \Z$.  In \cite{B memoirs}, we prove the following proposition.

\begin{prop}\label{NS bracket}(\cite{B memoirs})
The vector space $\hat{T}_e SK(1)$ with the bracket operation defined
in \cite{B memoirs} is the $N=1$ Neveu-Schwarz algebra with central
charge zero.  The basis is given by 
\begin{eqnarray}
\mathcal{L} (j) &=& \biggl. - \frac{\partial}{\partial A_{-j}^{(0)}}
\biggr|_e, \quad \mbox{for} \; \; - j \in \Z ,\\
\mathcal{L} (j) &=& \biggl. - \frac{\partial}{\partial A_j^{(1)}}
\biggr|_e, \quad \mbox{for} \; \; j \in \Z ,\\
\mathcal{L} (0) &=& \Bigl. - \frac{1}{2} \frac{\partial}{\partial a}
\Bigr|_e, \\ 
\mathcal{G} (j + \frac{1}{2}) &=& \Biggl. - \frac{\partial}{\partial M_{-j
- \frac{1}{2}}^{(0)}} \Biggr|_e, \quad \mbox{for} \; \; - j \in \Z ,\\
\mathcal{G} (j - \frac{1}{2}) &=& \Biggl. - \frac{\partial}{\partial M_{j
- \frac{1}{2}}^{(1)}} \Biggr|_e, \quad \mbox{for} \; \; j \in \Z .
\end{eqnarray}
\end{prop}

\section{The sewing identities on modules for the $N=1$ Neveu-Schwarz 
algebra}

In this section we recall {}from \cite{B memoirs} the sewing identities
on modules for the $N=1$ Neveu-Schwarz algebra, $\mathfrak{ns}$, arising 
{}from the sewing operation on the moduli space of superspheres with
tubes. 

\subsection{Modules for the $\mathfrak{ns}$}\label{modules}

For any representation of $\mathfrak{ns}$, we shall use $L(j)$, $G(j -
1/2)$ and $c \in \mathbb{C}$ to denote the representation images
of $L_j$, $G_{j - 1/2}$ and $d$, respectively.  

Let  $V = \coprod_{k \in (1/2) \mathbb{Z}} V_{(k)}$
be a module for $\mathfrak{ns}$ of central charge $c \in \mathbb{C}$ 
(i.e., $dv = cv$ for $v \in V$) such that for $v \in V_{(k)}$, we have 
$L(0)v = kv$.  Let $P(k)$ be the projection {}from $V$ to $V_{(k)}$.  For
any formal  variable $t^{1/2}$ and $l \in 2\mathbb{Z} \smallsetminus
\{0\}$, we  define $(t^{1/2})^{lL(0)} \in (\mathrm{End} \; V) [[t^{1/2},
t^{-1/2}]]$  by $(t^{1/2})^{lL(0)}v = (t^{1/2})^{kl}v$ for $v \in
V_{(k)}$, or equivalently $(t^{1/2})^{lL(0)} v= \sum_{k \in
(1/2)\mathbb{Z}} P(k)  (t^{1/2})^{kl}v$ for any $v \in V$.

Let $V_P$ be a vector space over $\mathbb{C}$ with basis $\{P_k \; | \; 
k \in (1/2) \mathbb{Z}\}$.  Let $T(\mathfrak{ns} \oplus V_P)$ be the 
tensor algebra generated by the direct sum of $\mathfrak{ns}$ and $V_P$, 
and let $\mathcal{I}$ be the ideal of $T(\mathfrak{ns} \oplus V_P)$ 
generated by
\begin{multline*}
\Bigl\{L_m \otimes L_n - L_n \otimes L_m - [L_m, L_n], \;  
L_m \otimes G_{n - \frac{1}{2}} - G_{n - \frac{1}{2}} \otimes L_m -
\bigl[L_m , G_{n - \frac{1}{2}}\bigr], \\
G_{m + \frac{1}{2}} \otimes G_{n - \frac{1}{2}} + G_{n - \frac{1}{2}} 
\otimes G_{m + \frac{1}{2}}  - \bigl[G_{m + \frac{1}{2}}, G_{n - 
\frac{1}{2}}\bigr], \; L_n \otimes d - d \otimes L_n,  \\
\Bigl. G_{n - \frac{1}{2}} \otimes d - d \otimes G_{n - \frac{1}{2}}, 
\; P_j \otimes P_k - \delta_{j,k} P_j, \; P_j \otimes L_n - L_n \otimes 
P_{j + n}, \Bigr. \\
\bigl. P_j \otimes G_{n - \frac{1}{2}} - G_{n - \frac{1}{2}}
\otimes P_{j + n -\frac{1}{2}}, \; P_j \otimes d - d \otimes P_j \; \bigr|
\; m,n 
\in \mathbb{Z}, \; j,k \in \mbox{$\frac{1}{2}\mathbb{Z}$} \Bigr\} .
\end{multline*}

Then $U_P(\mathfrak{ns}) = T(\mathfrak{ns} \oplus V_P) / \mathcal{I}$ is an
associative superalgebra and the universal enveloping algebra
$U(\mathfrak{ns})$ of $\mathfrak{ns}$ is a subalgebra.  Linearly
$U_P(\mathfrak{ns}) = U(\mathfrak{ns}) \otimes V_P$.

For any even formal variable $t^{1/2}$ and $l \in 2 \mathbb{Z} 
\smallsetminus \{0\}$, we define
\[(t^\frac{1}{2})^{lL_0} = \sum_{k \in \frac{1}{2}\mathbb{Z}} P_k 
(t^\frac{1}{2})^{kl} \; \in \; U_P(\mathfrak{ns}) [[t^\frac{1}{2},
t^{-\frac{1}{2}}]] .\] Then for $l,n \in 2\mathbb{Z} \smallsetminus
\{0\}$, and $m \in
\mathbb{Z}$, 
\begin{eqnarray*}
(t^\frac{1}{2})^{lL_0} (t^\frac{1}{2})^{nL_0} &=& (t^\frac{1}{2})^{(l + n)L_0},  \\
(t^\frac{1}{2})^{lL_0} L_m &=& L_m (t^\frac{1}{2})^{-lm} (t^\frac{1}{2})^{lL_0}, \\
(t^\frac{1}{2})^{lL_0} G_{m - \frac{1}{2}} &=& G_{m - \frac{1}{2}} 
(t^\frac{1}{2})^{-l(m - \frac{1}{2})} t^{lL_0},
\end{eqnarray*} 
and $(t^{1/2})^{lL_0}$ commutes with $d$. {}From \cite{B memoirs}, we have the following proposition.   

\begin{prop}\label{universal enveloping}
Let $V$ be a module for $\mathfrak{ns}$ as above.  Then there is a unique
algebra homomorphism {}from $U_P (\mathfrak{ns})$ to $\mbox{End} \; V$
such that  $L_j$, $G_{j - 1/2}$, $d$ and $P_k$ are mapped to $L(j)$, $G(j
- 1/2)$, $c$ and $P(k)$, respectively, for $j \in \mathbb{Z}$ and
$k \in \frac{1}{2} \mathbb{Z}$.
\end{prop}

If $V_{(k)} = 0$ for $n$ sufficiently small, then we say that $V$ is a 
{\it positive energy module} for $\mathfrak{ns}$ and that the
corresponding representation is a {\it positive energy representation} of
$\mathfrak{ns}$. 

\subsection{The series $\Gamma$ and $\{\Psi_j,\Psi_{j-1/2}\}_{j \in
\mathbb{Z}}$}

Let $\mathcal{A} = \{\mathcal{A}_j \}_{j \in \mathbb{Z}_+}$ and
$\mathcal{B} = \{\mathcal{B}_j \}_{j \in \mathbb{Z}_+}$ be two sequences
of even formal variables, let $\mathcal{M} = \{\mathcal{M}_{j-1/2} \}_{j
\in \mathbb{Z}_+}$ and $\mathcal{N} = \{\mathcal{N}_{j-1/2} \}_{j \in
\mathbb{Z}_+}$ be two sequences of odd formal variables, and let
$\alpha_0^{1/2}$ be another even formal variable.

\begin{prop}(\cite{B memoirs}, \cite{BHL})\label{central charge sewing}
Let $\Psi_j$, for $j \in \frac{1}{2}\mathbb{Z}$, be the canonical
sequence of formal series given by Theorem \ref{actual sewing}.  There
exists a unique canonical formal series $\Gamma =
\Gamma(\alpha_0^{1/2}, \mathcal{A}, \mathcal{M}, \mathcal{B},
\mathcal{N}) \in (\mathbb{C}[\alpha_0^{1/2}, \alpha_0^{-1/2}]
[[\mathcal{A}, \mathcal{B}]] [\mathcal{M}, \mathcal{N}])^0$ 
such that
\begin{equation}\label{first terms of Gamma}
\Gamma =  \sum_{j \in \Z} \! \left( \! \biggl( \frac{j^3 - j}{12} \biggr)
\alpha_0^{-j} \mathcal{A}_j \mathcal{B}_j + \biggl( \frac{j^2 - j}{3} \biggr)
\alpha_0^{-j + \frac{1}{2}} \mathcal{N}_{j - \frac{1}{2}} \mathcal{M}_{j -
\frac{1}{2}} \right) + \Gamma_0
\end{equation}
where $\Gamma_0$ contains only terms with total degree at least three
in the $\mathcal{A}_j$'s, $\mathcal{M}_{j - 1/2}$'s, $\mathcal{B}_j$'s,
and $\mathcal{N}_{j - 1/2}$'s for $j \in \Z$ with each term
containing at least one of the $\mathcal{A}_j$'s or $\mathcal{M}_{j -
1/2}$'s and at least one of the $\mathcal{B}_j$'s or $\mathcal{N}_{j
- 1/2}$'s, and such that in the $\mathbb{C}[\alpha_0^{1/2}, 
\alpha_0^{-1/2}] [[\mathcal{A}, \mathcal{B}]][\mathcal{M},
\mathcal{N}]$-envelope of $U_P(\mathfrak{ns})$, we have
\begin{multline}\label{switching identity in NS}
e^{- \sum_{j \in \Z} (\mathcal{A}_j L_j + \mathcal{M}_{j - \frac{1}{2}}
G_{j - \frac{1}{2}})} \; (\alpha_0^\frac{1}{2})^{-2L_0} \; e^{-
\sum_{j \in \Z} (\mathcal{B}_j L_{-j} + \mathcal{N}_{j - \frac{1}{2}} G_{- j
+ \frac{1}{2}})} \\
= \; e^{\sum_{j \in \Z} (\Psi_{-j} L_{-j} + \Psi_{- j +
\frac{1}{2}} G_{- j + \frac{1}{2}})} \; e^{\sum_{j \in \Z} (\Psi_j L_j
+ \Psi_{j - \frac{1}{2}} G_{j - \frac{1}{2}})} \; e^{2\Psi_0 L_0} \;
(\alpha_0^\frac{1}{2})^{-2L_0} \; e^{\Gamma c} . 
\end{multline}
\end{prop}

\begin{rema}  {\em In \cite{B memoirs} we use the sewing identities on
the moduli space on $N=1$ superspheres with tubes to prove
Proposition \ref{central charge sewing} above.  In \cite{BHL}, we give 
a more straightforward and Lie-theoretic proof of Proposition 
\ref{central charge sewing} by proving a certain bijectivity property 
for the Campbell-Baker-Hausdorff formula in the theory of Lie algebras.  
This allows us to prove Proposition \ref{central charge sewing} 
directly for $U_P(\mathfrak{ns})$ rather than go through the
representation in terms of superderivations and then lift to
$U_P(\mathfrak{ns})$ as we did in \cite{B thesis}, \cite{B memoirs}.}  
\end{rema}

Let $V = \coprod_{k \in (1/2) \mathbb{Z}} V_{(k)}$ be a module for
$\mathfrak{ns}$, and let $L(j), G(j - 1/2) \in \mbox{End} \; V$ and $c
\in \mathbb{C}$ be the representation images of $L_j$, $G_{j - 1/2}$ and
$d$, respectively, such that $L(0)v = kv$ for $v \in V_{(k)}$.  

\begin{cor}(\cite{B memoirs})\label{normal order in End}
In the $\mathbb{C}[\alpha_0^{1/2}, \alpha_0^{-1/2}] [[\mathcal{A}, 
\mathcal{B}]] [\mathcal{M}, \mathcal{N}]$-envelope of $\mbox{\em End} \;
V$,  we have 
\begin{multline}\label{psi gamma corollary}
e^{- \sum_{j \in \Z} (\mathcal{A}_j L(j) + \mathcal{M}_{j - \frac{1}{2}}
G(j - \frac{1}{2}))} \cdot (\alpha_0^\frac{1}{2})^{-2L(0)} \cdot e^{-
\sum_{j \in \Z} (\mathcal{B}_j L(-j) + \mathcal{N}_{j - \frac{1}{2}} G(- j +
\frac{1}{2}))} 
\end{multline}
\begin{multline*}
= \; e^{\sum_{j \in \Z} (\Psi_{-j} L(-j) + \Psi_{- j + \frac{1}{2}} G(- j
+ \frac{1}{2}))} \cdot e^{\sum_{j \in \Z} (\Psi_j L(j) + \Psi_{j -
\frac{1}{2}} G(j - \frac{1}{2}))} \cdot \\
\cdot e^{2\Psi_0 L(0)} \cdot (\alpha_0^\frac{1}{2})^{-2L(0)} \cdot e^{\Gamma c} .
\end{multline*}
Furthermore, the formal series $\Gamma$ and $\Psi_j$, for $j \in
\frac{1}{2}\mathbb{Z}$, are actually in 
\[\mathbb{C}[\mathcal{A}][\mathcal{M}][\mathcal{B}][\mathcal{N}]
[[\alpha_0^{-\frac{1}{2}} ]] . \]  
Thus letting $(A,M)$ and $(B,N)$ be two sequences in $\bigwedge_*^\infty$,
letting $\asqrt \in (\bigwedge_*^0)^\times$, and letting $t^{1/2}$ be an even 
formal variable,  we have that $\Psi_j(t^{-1/2}\asqrt, A, M, B, N)$, for $j \in
\frac{1}{2}\mathbb{Z}$, and $\Gamma(t^{-1/2} \asqrt, A,M,B,N)$ are well
defined and belong to $\bigwedge_*[[t^{1/2}]]$.  If in addition, 
$V$ is a positive-energy module for $\mathfrak{ns}$,  then the following
identity holds in
$(\mbox{\em End} \; (\bigwedge_* \otimes_\mathbb{C} V))^0 [[t^{1/2}]]$.
\begin{multline}\label{positive energy gamma equation}
e^{- \sum_{j \in \Z} (A_j L(j) + M_{j - \frac{1}{2}}
G(j - \frac{1}{2}))} \cdot (t^{-\frac{1}{2}}\asqrt)^{-2L(0)} \cdot\\
\cdot e^{- \sum_{j \in \Z} (\mathcal{B}_j L(-j) + \mathcal{N}_{j - \frac{1}{2}} 
G(- j + \frac{1}{2}))} \cdot 
\end{multline}
\begin{multline*}
= \; e^{\sum_{j \in \Z} (\Psi_{-j}(t^{-\frac{1}{2}} \asqrt, A, M, B, N) L(-j) +
\Psi_{- j + \frac{1}{2}}(t^{-\frac{1}{2}} \asqrt, A, M, B, N) G(- j +
\frac{1}{2}))} \cdot \\
\cdot e^{\sum_{j \in \Z} (\Psi_j(t^{-\frac{1}{2}}\asqrt, A, M, B, N) L(j) +
\Psi_{j - \frac{1}{2}}(t^{-\frac{1}{2}} \asqrt, A, M, B, N) G(j -
\frac{1}{2}))}  \cdot\\
\cdot e^{2\Psi_0(t^{-\frac{1}{2}} \asqrt, A, M, B, N) L(0)} \cdot (t^{-\frac{1}{2}}
\asqrt)^{-2L(0)} \cdot e^{\Gamma(t^{-\frac{1}{2}} \asqrt, A,M,B,N)c} . 
\end{multline*}  
\end{cor}

Before introducing the definition of $N=1$ supergeometric vertex 
operator algebra, we present the following conjecture.  The Isomorphism
Theorem, Theorem \ref{supergeometric and superalgebraic}, only holds
when this conjecture is true, or in the trivial case in which the
representation of the Neveu-Schwarz algebra being considered has
zero central charge.

\begin{conj}\label{Gamma converges}
Let $Q_1 \in SK(m)$ and $Q_2 \in SK(n)$ be given by 
\begin{multline*}
Q_1 = ((z_1, \theta_1),...,(z_{m-1}, \theta_{m-1}); (A^{(0)},
M^{(0)}), (a^{(1)}, A^{(1)}, M^{(1)}),\\
...,(a^{(m)}, A^{(m)}, M^{(m)})), 
\end{multline*} 
and
\[Q_2 = ((z_1', \theta_1'),...,(z_{n-1}', \theta_{n-1}'); (B^{(0)},
N^{(0)}), (b^{(1)}, B^{(1)}, N^{(1)}),...,(b^{(n)}, B^{(n)}, N^{(n)}))
.\] 
If the $i$-th tube of $Q_1$ can be sewn with the 0-th tube of $Q_2$, then the
$t^{1/2}$-series
\[e^{\Gamma (t^{-\frac{1}{2}}a^{(i)}, A^{(i)}, M^{(i)}, B^{(0)},
N^{(0)})c}\] 
is absolutely convergent at $t^{1/2} = 1$.
\end{conj}

We have already noted the convergence property of the $\Psi_j$'s in
Theorem \ref{actual sewing}.

\subsection{The series $\{\Theta^{(1)}_j,\Theta^{(1)}_{j- 1/2}\}_{j \in
\mathbb{Z}}$ and $\{\Theta^{(2)}_j,\Theta^{(2)}_{j-1/2}\}_{j \in
\mathbb{Z}}$}

Let 
\begin{eqnarray*}
H^{(1)}_{\alpha_0^{1/2},\mathcal{A}, \mathcal{M}} (x, \varphi)  &=& 
\exp \Biggl(\!  - \! \sum_{j \in \Z}
\biggl( \mathcal{A}_j L_j(x,\varphi)  +  \mathcal{M}_{j - \frac{1}{2}}
G_{j - \frac{1}{2}}(x,\varphi) \biggr) \! \Biggr) \cdot \\
& & \hspace{2in} \cdot (\alpha_0^{1/2})^{-2L_0(x,\varphi)}\cdot (x,
\varphi ) \\ 
&=& \hat{E}(\alpha_0^{1/2},\mathcal{A}, \mathcal{M} )
(x,\varphi).
\end{eqnarray*}
Let 
\[(\tilde{x}, \tilde{\varphi}) = (H^{(1)}_{\alpha_0^{1/2},\mathcal{A}, 
\mathcal{M}})^{-1} (x,\varphi) \in (\alpha_0^{-1} x,\alpha_0^{-\frac{1}{2}}
\varphi) + x\mathbb{C}[x,\varphi][\alpha_0^{-\frac{1}{2}}][[\mathcal{A}]]
[\mathcal{M}] . \]
Let $w$ be another even formal variable and $\rho$ another odd formal
variable.  Then $s_{(x,\varphi)} \circ  H^{(1)}_{\alpha_0^{1/2},\mathcal{A},
\mathcal{M}}
\circ s_{(\tilde{x}, \tilde{\varphi})}^{-1} (\alpha_0^{-1} w, 
\alpha_0^{-1/2}\rho)$ is in 
\[w\mathbb{C} [x, \varphi] [\alpha_0^{\frac{1}{2}}, \alpha_0^{-\frac{1}{2}}] 
[[\mathcal{A}]] [\mathcal{M}][[w]] \oplus \rho \mathbb{C}[x, \varphi] 
[\alpha_0^{\frac{1}{2}}, \alpha_0^{-\frac{1}{2}}][[\mathcal{A}]][\mathcal{M}][[w]] 
,\] 
is superconformal in $(w,\rho)$, and the even coefficient of the monomial
$\rho$ is an element in $(1 + x\mathbb{C}[x][\alpha_0^{1/2},
\alpha_0^{-1/2}] [[\mathcal{A}]][\mathcal{M}] \oplus \varphi \mathbb{C}
[x][\alpha_0^{1/2}, \alpha_0^{-1/2}] [[\mathcal{A}]][\mathcal{M}])$.

Let $\Theta^{(1)}_j = \Theta^{(1)}_j(\alpha_0^{1/2},\mathcal{A}, 
\mathcal{M}, (x, \varphi)) \in \mathbb{C}[x, \varphi] [\alpha_0^{1/2},
\alpha_0^{-1/2}][[\mathcal{A}]][\mathcal{M}]$, for $j \in \frac{1}{2}
\mathbb{N}$, be defined by
\begin{multline}\label{define first Theta}
\Bigl(\exp(\Theta^{(1)}_0(\alpha_0^{\frac{1}{2}}, \mathcal{A}, \mathcal{M}, 
(x, \varphi)), \Bigl\{\Theta^{(1)}_j(\alpha_0^{\frac{1}{2}}, \mathcal{A}, 
\mathcal{M}, (x, \varphi)), \Bigr. \Bigr. \\
\Bigl. \Bigl. \Theta^{(1)}_{j - \frac{1}{2}} (\alpha_0^{\frac{1}{2}}, 
\mathcal{A}, \mathcal{M}, (x, \varphi)) \Bigr\}_{j \in \Z} \Bigr)
\end{multline}  
\[ = \;\hat{E}^{-1}(s_{(x,\varphi)} \circ H^{(1)}_{\alpha_0^{1/2},\mathcal{A}, 
\mathcal{M}} \circ s_{(\tilde{x},\tilde{\varphi})}^{-1}
(\alpha_0^{-1}w,\alpha_0^{-\frac{1}{2}}\rho)) .  
\hspace{1.5in} \]
In other words, the $\Theta_j^{(1)}$'s are determined uniquely by
\begin{multline*}
s_{(x,\varphi)} \circ H^{(1)}_{\alpha_0^{1/2},\mathcal{A}, \mathcal{M}}
\circ s_{(\tilde{x}, \tilde{\varphi})}^{-1} (\alpha_0^{-1} w, 
\alpha_0^{-\frac{1}{2}}\rho) \\  
= \exp \Biggl( \sum_{j \in  \Z} \biggl( \Theta^{(1)}_j \! \left( \Lw \right)
+ \Theta^{(1)}_{j - \frac{1}{2}} \Gw \biggr) \! \Biggr) \! \cdot \\
\exp \left(\! \Theta^{(1)}_0 \left( \twoLow \right) \! \right) \! \cdot
(w,\rho) . 
\end{multline*}  

This formal power series in $(w,\rho)$ gives the formal local 
superconformal coordinate at a puncture of the canonical supersphere 
obtained {}from the sewing together of two particular canonical 
superspheres with punctures.  Specifically, given a canonical supersphere
$Q_1$ representing an element in $SK(m)$
\begin{multline*}
Q_1 = ((z_1, \theta_1), \dots, (z_{m-1}, \theta_{m-1});
(A^{(0)}, M^{(0)}),  (\asqrt^{(1)}, A^{(1)}, M^{(1)}), \dots, \\
(\asqrt^{(m)}, A^{(m)}, M^{(m)}))
\end{multline*}
and a canonical supersphere $Q_2 \in SK(2)$ with standard local coordinates
given by
\[Q_2 = ((z,\theta); I(w,\rho), s_{(z,\theta)}(w,\rho), (w,\rho)) =
((z,\theta) ; \mathbf{0}, (1, \mathbf{0}), (1, \mathbf{0}))\] 
and if $|z_B| < r$ where $r$ is the radius of convergence of
$\hat{E}(\asqrt^{(m)}, A^{(m)}, M^{(m)}) (w,\rho) =
H^{(1)}_{(\asqrt^{(m)}, A^{(m)}, M^{(m)})}(w,\rho)$, then we can sew the
puncture at infinity of $Q_2$ to the $m$-th puncture of
$Q_1$.  In this case, the uniformizing function is given by $F_\sou^{(1)}
(w,\rho) = (w, \rho)$ and $F_\sou^{(2)} (w,\rho) = (\hat{E}(\asqrt^{(m)},
A^{(m)}, M^{(m)}))^{-1} (w, \rho) = (H^{(1)}_{(\asqrt^{(m)}, A^{(m)},
M^{(m)})})^{-1} (w,\rho)$.  Thus the resulting canonical supersphere is
given by case (1) of Theorem \ref{actual sewing}, and we have
\begin{eqnarray*}
Q_1 \; _m\infty_0 \; Q_2 &=& \bigl((z_1, \theta_1),...,(z_{m-1},
\theta_{m-1}), (H^{(1)}_m)^{-1}(z, \theta); (A^{(0)}, M^{(0)}), \\
& & \qquad (\asqrt^{(1)}, A^{(1)},M^{(1)}), \dots, (\asqrt^{(m-1)},
A^{(m-1)}, M^{(m-1)}), \\
& & \quad \hat{E}^{-1} \circ s_{(z,\theta)}
\circ H^{(1)}_m \circ s_{(H^{(1)}_m)^{-1}(z, \theta)}^{-1}
((\asqrt^{(m)})^{-2} w, (\asqrt^{(m)})^{-1}\rho) \bigr)\\
&=& \bigl( (z_1, \theta_1),...,(z_{m-1}, \theta_{m-1}),
(H^{(1)}_m)^{-1}(z, \theta); (A^{(0)}, M^{(0)}), \\
& & \qquad (\asqrt^{(1)}, A^{(1)},M^{(1)}), \dots, (\asqrt^{(m-1)},
A^{(m-1)}, M^{(m-1)}), \\
& & \hspace{1.8in} (\exp \Theta^{(1)}_0, \{ \Theta^{(1)}_j ,
\Theta^{(1)}_{j - \frac{1}{2}} \}_{j \in \mathbb{Z}_+}) \bigr)
\end{eqnarray*}
where $(\exp \Theta^{(1)}_0, \{ \Theta^{(1)}_j ,
\Theta^{(1)}_{j - \frac{1}{2}} \}_{j \in \mathbb{Z}_+})$ is given
formally by (\ref{define first Theta}), with
$(\alpha_0^{1/2},\mathcal{A},\mathcal{M}) = (\asqrt^{(m)},
A^{(m)},M^{(m)})$ and $(x,\varphi) = (z,\theta)$.

In Section 7, in order to prove that given an $N=1$ NS-VOSA we can
construct an $N=1$ SG-VOSA, we will need the following proposition
concerning the series $\Theta^{(1)}_j$, for $j \in \frac{1}{2} \mathbb{N}$.

\begin{prop}(\cite{B memoirs}) \label{first Theta prop in
universal enveloping algebra} In $(U(\mathfrak{ns}) [x,
\varphi][\alpha_0^{1/2}, 
\alpha_0^{-1/2}] [[\mathcal{A}]][\mathcal{M}])^0$, we have   
\begin{multline*}
\exp \Biggl(\! - \! \! \sum_{m = -1}^{\infty} \sum_{j \in \Z}
\binom{j+1}{m+1} \alpha_0^{-j} x^{j - m} \biggl(\!  \Bigl( \mathcal{A}_j
+ 2\left(\frac{j-m}{j+1} \right) \alpha_0^{\frac{1}{2}} x^{-1} \varphi
\mathcal{M}_{j - \frac{1}{2}} \Bigr)  L_m \\
+ \; x^{-1} \Bigl( \left(\frac{j-m}{j+1} \right) \alpha_0^{\frac{1}{2}}
\mathcal{M}_{j - \frac{1}{2}} + \varphi \frac{(j-m)}{2} \mathcal{A}_j
\Bigr)  G_{m + \frac{1}{2}} \biggl) \Biggr) 
\end{multline*}
\[= \; e^{\left(\!(\tilde{x} - \alpha_0^{-1}x) L_{-1}  + (\tilde{\varphi} -
\alpha_0^{-1/2}\varphi) G_{-\frac{1}{2}} \right)} \cdot e^{\left( - \! \sum_{j \in 
\Z} \left( \Theta^{(1)}_j L_j + \Theta^{(1)}_{j - \frac{1}{2}} G_{j - 
\frac{1}{2}} \right) \! \right)} \cdot e^{\left( - 2 \Theta^{(1)}_0 L_0 \right)},  \]
for $(\tilde{x}, \tilde{\varphi}) = (H_{\alpha_0^{1/2},\mathcal{A}, 
\mathcal{M}}^{(1)})^{-1} (x,\varphi)$.  
\end{prop}

Let $V = \coprod_{k \in (1/2) \mathbb{Z}} V_{(k)}$ be a module for
$\mathfrak{ns}$, and let $L(j), G(j - 1/2) \in \mbox{End} \; V$ and $c \in
\mathbb{C}$ be the representation images of $L_j$, $G_{j - 1/2}$ and $d$,
respectively, such that $L(0)v = kv$ for $v \in V_{(k)}$.  

\begin{cor}(\cite{B memoirs}) \label{first Theta identity in End}
In $((\mbox{\em End} \; V) [x, \varphi][\alpha_0^{1/2}, 
\alpha_0^{-1/2}] [[\mathcal{A}]][\mathcal{M}])^0$, we have 
\begin{multline}\label{first Theta corollary}
\exp \Biggl(\! - \! \! \sum_{m = -1}^{\infty} \sum_{j \in \Z}
\binom{j+1}{m+1} \alpha_0^{-j} x^{j - m}  \\
\Biggl. \biggl(\!  \Bigl( \mathcal{A}_j + 2\left(\frac{j-m}{j+1} \right)
\alpha_0^{\frac{1}{2}} x^{-1} \varphi \mathcal{M}_{j - \frac{1}{2}}
\Bigr)  L(m) \\
+ \; x^{-1} \Bigl( \left(\frac{j-m}{j+1} \right) \alpha_0^{\frac{1}{2}}
\mathcal{M}_{j - \frac{1}{2}} + \varphi \frac{(j-m)}{2} \mathcal{A}_j
\Bigr)  G(m + \frac{1}{2}) \biggl) \Biggr) 
\end{multline} 
\begin{multline*} 
= \;  e^{\left((\tilde{x} - \alpha_0^{-1}x) L(-1)  + (\tilde{\varphi} - 
\alpha_0^{-1/2}\varphi)
G(-\frac{1}{2}) \right)} \cdot e^{\left( - \sum_{j \in \Z} \left( \Theta^{(1)}_j
L(j) + \Theta^{(1)}_{j - \frac{1}{2}} G(j - \frac{1}{2}) \right) \right)} 
 \cdot \\
\cdot e^{\left( - 2\Theta^{(1)}_0 L(0) \right)} 
\end{multline*} 
for $(\tilde{x}, \tilde{\varphi}) = (H_{\alpha_0^{1/2},\mathcal{A}, 
\mathcal{M}}^{(1)})^{-1} (x,\varphi)$.

Furthermore, the formal series $\Theta^{(1)}_j (t^{ - 1/2} \alpha_0^{1/2}, \mathcal{A}, 
\mathcal{M}, (x,\varphi))$, for $j \in \frac{1}{2} \mathbb{N}$, are in 
$\mathbb{C} [x, \varphi][\mathcal{A}][\mathcal{M}][\alpha_0^{-1/2}]
[[t^{1/2}]]$.  Thus for $\asqrt \in (\bigwedge_*^0)^\times$, and $(A,M) \in
\bigwedge_*^\infty$, the series $\Theta^{(1)}_j (t^{ - 1/2} \asqrt,
A,M, (x,\varphi))$,  for $j \in \frac{1}{2}\mathbb{N}$, are well defined
and belong to $\bigwedge_*  [x, \varphi][[t^{1/2}]]$. 

Let $(\tilde{x}(t^{1/2}),\tilde{\varphi} (t^{1/2})) = 
H^{(1)}_{ t^{ - 1/2}\asqrt,A,M}(x,\varphi)$, and let $V$ be a
positive-energy module for $\mathfrak{ns}$.  Then in
$((\mbox{\em End} \; (\bigwedge_* \otimes_\mathbb{C} V))
[[t^{\frac{1}{2}}]] [[x,x^{-1}]][\varphi])^0$, the following identity
holds for
$\Theta^{(1)}_j =
\Theta^{(1)}_j (t^{-1/2}\asqrt, A,M,(x,\varphi))$. 
\begin{multline}\label{last first Theta equation}
\exp \Biggl(\! - \! \! \sum_{m = -1}^{\infty} \sum_{j \in \Z}
\binom{j+1}{m+1} t^j \asqrt^{-2j} x^{j - m}  \\
\biggl(\!  \Bigl( A_j + 2\left(\frac{j-m}{j+1} \right)
t^{-\frac{1}{2}} \asqrt x^{-1} \varphi M_{j - \frac{1}{2}}
\Bigr)  L(m) \\
+ \; x^{-1} \Bigl( \left(\frac{j-m}{j+1} \right) t^{-\frac{1}{2}}
\asqrt M_{j - \frac{1}{2}} + \varphi \frac{(j-m)}{2} A_j
\Bigr)  G(m + \frac{1}{2}) \biggl) \Biggr) 
\end{multline}
\begin{multline*}
= \; \exp \left( (\tilde{x}(t^\frac{1}{2}) - t \asqrt^{-2}x) L(-1) +
(\tilde{\varphi}(t^\frac{1}{2}) - t^{\frac{1}{2}}\asqrt^{-1} \varphi) 
G(-\frac{1}{2}) \right) \cdot \\
\exp \Biggl( \! - \! \sum_{j \in \Z} \Bigl( \Theta^{(1)}_j L(j) +
\Theta^{(1)}_{j - \frac{1}{2}} G(j - \frac{1}{2}) \Bigr) \! \Biggr) \! \cdot \exp
\left( - 2\Theta^{(1)}_0 L(0) \right) . 
\end{multline*}
\end{cor}

Let
\begin{eqnarray*}
H^{(2)} (x, \varphi) \! &=& \! \exp \Biggl(\sum_{j \in \Z} \biggl(
\mathcal{B}_j L_{-j}(x,\varphi) +  \mathcal{N}_{j - \frac{1}{2}} G_{-j +
\frac{1}{2}}(x,\varphi) \biggr) \! \Biggr) \cdot \Bigl(\frac{1}{x}, \frac{i
\varphi}{x} \Bigr) \\ 
&=& \tilde{E}(\mathcal{B},-i \mathcal{N}) \Bigl(\frac{1}{x},
\frac{i \varphi}{x} \Bigr). 
\end{eqnarray*}
Let
\[(\tilde{x}, \tilde{\varphi}) = (H^{(2)}_{\mathcal{B}, \mathcal{N}})^{-1}
\circ I (x,\varphi) \in (x,\varphi) + \mathbb{C}[x^{-1}, \varphi]
[[\mathcal{B}]][\mathcal{N}] . \]
Let $w$ be another even formal variable and $\rho$ another odd formal
variable.  Now write $s_{(x,\varphi)}(w,\rho) = (-x + w - \rho 
\varphi, \rho - \varphi)$.  We will use the convention that we should expand
$(-x + w -\rho \varphi)^j = (-x + w)^j - j\rho \varphi (-x + w)^{j-1}$ 
in positive powers of the second even variable $w$, for $j \in
\mathbb{Z}$. Then 
\[s_{(x,\varphi)} \circ I^{-1} \circ H^{(2)}_{\mathcal{B},\mathcal{N}}
\circ s_{(\tilde{x}, \tilde{\varphi})}^{-1} (w,\rho)  \in  
w\mathbb{C} [x^{-1} \!, \varphi][[\mathcal{B}]][\mathcal{N}][[w]] \oplus 
\rho \mathbb{C} [x^{-1} \!, \varphi][[\mathcal{B}]][\mathcal{N}][[w]] ,\]
is superconformal in $(w,\rho)$, and the even coefficient of the monomial
$\rho$ is an element in $(1 + x^{-1} \mathbb{C} [x^{-1},
\varphi][[\mathcal{B}]][\mathcal{N}])$.  

Let $\Theta^{(2)}_j = \Theta^{(2)}_j(\mathcal{B}, \mathcal{N}, (x, \varphi))
\in \mathbb{C}[x^{-1}, \varphi][[\mathcal{B}]][\mathcal{N}]$, 
for $j \in \frac{1}{2} \mathbb{N}$, be defined by
\begin{equation}\label{define second Theta}
\Bigl(\exp(\Theta^{(2)}_0(\mathcal{B}, \mathcal{N}, (x, \varphi)), \Bigl\{
\Theta^{(2)}_j(\mathcal{B},\mathcal{N}, (x, \varphi)),
\Theta^{(2)}_{j - \frac{1}{2}} (\mathcal{B},
\mathcal{N}, (x, \varphi)) \Bigr\}_{j \in \Z} \Bigr) 
\end{equation}
\[ = \; \hat{E}^{-1}(s_{(x,\varphi)} \circ I^{-1} \circ
H^{(2)}_{\mathcal{B}, \mathcal{N}} \circ s_{(\tilde{x},
\tilde{\varphi})}^{-1} (w,\rho)) . \hspace{1.6in}  \] 
In other words, the $\Theta_j^{(2)}$'s are determined uniquely by
\begin{multline*}
s_{(x,\varphi)} \circ I^{-1} \circ H^{(2)}_{\mathcal{B}, 
\mathcal{N}} \circ s_{(\tilde{x}, \tilde{\varphi})}^{-1} (w,\rho) \\ 
= \exp \Biggl( \sum_{j \in  \Z} \biggl( \Theta^{(2)}_j \! \left( \Lw \right)
+ \Theta^{(2)}_{j - \frac{1}{2}} \Gw \biggr) \! \Biggr) \! \cdot \\
\exp \left(\! \Theta^{(2)}_0 \left( \twoLow \right) \! \right) \! 
\cdot (w,\rho) . 
\end{multline*}

This formal power series in $(w,\rho)$ gives the formal local 
superconformal coordinate at a puncture of the canonical supersphere 
obtained {}from the sewing together of two particular canonical 
superspheres with punctures.  Specifically, given a canonical supersphere
$Q_1$ representing an element in $SK(2)$ with standard local coordinates
given by 
\[Q_1 = ((z,\theta); I(w,\rho), s_{(z,\theta)}(w,\rho), (w,\rho)) =
((z,\theta) ; \mathbf{0}, (1, \mathbf{0}), (1, \mathbf{0})) \in SK(2)\]
and a canonical supersphere $Q_2 \in SK(1)$ given by
\[Q_1 = ((B^{(0)}, N^{(0)}),  (\bsqrt^{(1)}, B^{(1)}, N^{(1)})) \in
SK(1) ,\]
and if $|z_B| > r$ where $r$ is the radius of convergence of
$\tilde{E} (B^{(0)}, -iN^{(0)}) (1/w,i\rho/w) = H_{(B^{(0)},
N^{(0)})}^{(2)}(w,\rho)$, then we can sew the puncture at infinity of
$Q_2$ to the second puncture of $Q_1$.  In this case, the uniformizing
function is given by $F_\sou^{(1)} (w,\rho) =  (H_{(B^{(0)},
N^{(0)})}^{(2)})^{-1} \circ I (w,\rho) $ and
$F_\sou^{(2)} (w,\rho) = (w,\rho)$.  Thus the resulting canonical
supersphere is given by case (1) of Theorem \ref{actual sewing}, and we
have
\begin{eqnarray*}
Q_1 \; _2\infty_0 \; Q_2 &=& \bigl((H_{(B^{(0)},N^{(0)})}^{(2)})^{-1}
\circ I  (z, \theta); (B^{(0)}, N^{(0)}), \\
& & \quad \hat{E}^{-1} \circ s_{(z,\theta)}
\circ I^{-1} \circ H^{(2)}_{(B^{(0)}, N^{(0)})} \circ
s_{(H^{(2)}_{(B^{(0)}, N^{(0)})})^{-1}(z,\theta)}^{-1} (w, \rho),\\
& & \hspace{2.6in} (\bsqrt^{(1)}, B^{(1)}, N^{(1)}) \bigr)\\ 
&=& \bigl((H_{(B^{(0)},N^{(0)})}^{(2)})^{-1}
\circ I  (z, \theta); (B^{(0)}, N^{(0)}), \\
& & \hspace{.8in} (\exp \Theta^{(2)}_0, \{\Theta^{(2)}_j , \Theta^{(2)}_{j
-
\frac{1}{2}} \}_{j \in \mathbb{Z}_+}), (\bsqrt^{(1)},
B^{(1)}, N^{(1)}) \bigr)
\end{eqnarray*}
where $(\exp \Theta^{(2)}_0, \{ \Theta^{(2)}_j ,
\Theta^{(2)}_{j - \frac{1}{2}} \}_{j \in \mathbb{Z}_+})$ is given
formally by (\ref{define second Theta}), with
$(\mathcal{B},\mathcal{N}) = (B^{(0)},N^{(0)})$ and $(x,\varphi) =
(z,\theta)$.

In Section 7, in order to prove that given an $N=1$ NS-VOSA we can
construct an $N=1$ SG-VOSA, we will need the following proposition
concerning the series $\Theta^{(2)}_j$, for $j \in \frac{1}{2} \mathbb{N}$.

\begin{prop} (\cite{B memoirs})\label{second Theta prop in universal
enveloping algebra} In $(U(\mathfrak{ns}) [x^{-1},
\varphi][[\mathcal{B}]][\mathcal{N}])^0$, we have   
\begin{multline*}
\exp \! \biggl( \sum_{m = -1}^{\infty} \sum_{j \in \Z}
\binom{-j + 1}{m + 1} x^{-j - m} \biggl( \Bigl( \mathcal{B}_j + 2\varphi 
\mathcal{N}_{j - \frac{1}{2}} \Bigr)  L_m \\   
+ \; \Bigl( \mathcal{N}_{j - \frac{1}{2}} + \varphi x^{-1}
\frac{(-j-m)}{2}  \mathcal{B}_j \Bigr) G_{m + \frac{1}{2}} \biggr)
\biggr) 
\end{multline*}   
\[= \; e^{\left((\tilde{x} - x) L_{-1}  + (\tilde{\varphi} - \varphi)
G_{-\frac{1}{2}} \right)} \cdot e^{\left( - \sum_{j \in \Z} \left( \Theta^{(2)}_j
L_j + \Theta^{(2)}_{j - \frac{1}{2}} G_{j - \frac{1}{2}} \right) \right)}
\cdot e^{\left( - 2 \Theta^{(2)}_0 L_0 \right)}, \hspace{.4in} \]
where $(\tilde{x}, \tilde{\varphi}) = (H_{\mathcal{B}, 
\mathcal{N}}^{(2)})^{-1} \circ I (x,\varphi)$.  
\end{prop}

Let $V = \coprod_{k \in (1/2) \mathbb{Z}} V_{(k)}$ be a module for
$\mathfrak{ns}$, and let $L(j), G(j - 1/2) \in \mbox{End} \; V$ and $c \in
\mathbb{C}$ be the representation images of $L_j$, $G_{j - 1/2}$ and $d$,
respectively, such that $L(0)v = kv$ for $v \in V_{(k)}$.

\begin{cor} (\cite{B memoirs})\label{second Theta identity in End}
In $((\mbox{\em End} \; V) [x^{-1}, \varphi] [[\mathcal{B}]]
[\mathcal{N}])^0$, we have 
\begin{multline}\label{second Theta corollary}
\exp \! \biggl(\sum_{m = -1}^{\infty} \sum_{j \in \Z}
\binom{-j + 1}{m + 1} x^{-j - m} \biggl( \Bigl( \mathcal{B}_j + 2\varphi 
\mathcal{N}_{j - \frac{1}{2}} \Bigr) L(m) \\
 + \; \Bigl( \mathcal{N}_{j - \frac{1}{2}} + \varphi x^{-1}
\frac{(-j-m)}{2} \mathcal{B}_j \Bigr)  G(m + \frac{1}{2}) \biggr) \biggr) 
\end{multline} 
\[= \; e^{\left((\tilde{x} - x) L(-1)  + (\tilde{\varphi} - \varphi)
G(-\frac{1}{2}) \right)} \cdot e^{\left( - \sum_{j \in \Z} \left( \Theta^{(2)}_j
L(j) + \Theta^{(2)}_{j - \frac{1}{2}} G(j - \frac{1}{2}) \right) \right)}
\cdot e^{\left( -2\Theta^{(2)}_0 L(0) \right)} ,\]
for $(\tilde{x}, \tilde{\varphi}) = (H_{\mathcal{B}, \mathcal{N}}^{(2)})^{-1} 
\circ I (x,\varphi)$.

Furthermore, the formal series $\Theta^{(2)}_j (\{t^k \mathcal{B}_k, t^{k - 1/2}
\mathcal{N}_{k - 1/2} \}_{k \in \Z}, (x,\varphi))$, for $j \in
\frac{1}{2} \mathbb{N}$, are in $\mathbb{C} [x^{-1}, \varphi]
[\mathcal{B}][\mathcal{N}][[t^{1/2}]]$.  Thus for $(B,N) \in \bigwedge_*^\infty$,
the series 
\[\Theta^{(2)}_j (\{t^k B_k, t^{k - \frac{1}{2}} N_{k - \frac{1}{2}} \}_{k
\in \Z}, (x,\varphi))\] 
are well defined and belong to $\bigwedge_* [x^{-1}, \varphi]
[[t^{1/2}]]$. 

Let $(\tilde{x}(t^{1/2}),\tilde{\varphi} (t^{1/2})) = (H_{ \{t^j B_j, 
t^{j - 1/2} N_{j - 1/2} \}_{j \in \Z}}^{(2)})^{-1} \circ I (x,\varphi)$, 
and let $V$ be a positive-energy module for $\mathfrak{ns}$.  Then in 
$((\mbox{\em End} \; (\bigwedge_* \otimes_\mathbb{C} V))
[[t^{1/2}]][[x,x^{-1}]][\varphi])^0$, the  following identity holds for
$\Theta^{(2)}_j = \Theta^{(2)}_j (\{t^k  B_k, t^{k - 1/2} N_{k - 1/2}
\}_{k \in \Z}, (x,\varphi))$. 
\begin{multline}\label{last second Theta equation}
\exp  \biggl( \sum_{m = -1}^{\infty} \sum_{j \in \Z} \binom{-j + 1}{m +
1} x^{-j - m} \biggl( \Bigl( t^j B_j +  2\varphi t^{j - \frac{1}{2}} N_{j
- \frac{1}{2}} \Bigr) L(m) \\
+ \; \Bigl( t^{j - \frac{1}{2}}  N_{j - \frac{1}{2}} + \varphi x^{-1}
\frac{(-j-m)}{2} t^j B_j \Bigr)  G(m + \frac{1}{2}) \biggr) 
\biggr)
\end{multline} 
\begin{multline*}
= \; \exp \left( (\tilde{x}(t^\frac{1}{2}) - x) L(-1) +
(\tilde{\varphi}(t^\frac{1}{2}) - \varphi) G(-\frac{1}{2}) \right)
\cdot \\  
\exp \Biggl( \! - \! \sum_{j \in \Z} \Bigl( \Theta^{(2)}_j L(j) +
\Theta^{(2)}_{j - \frac{1}{2}} G(j - \frac{1}{2}) \Bigr) \! \Biggr) \! \cdot \exp
\left( - 2\Theta^{(2)}_0 L(0) \right) . 
\end{multline*} 
\end{cor}

\section{$N=1$ supergeometric vertex operator superalgebras}

In this section we introduce the notion of $N=1$ supergeometric vertex
operator superalgebra.

\subsection{Linear algebra of $\frac{1}{2} \mathbb{Z}$-graded
$\bigwedge_*$-modules with finite-dimensional homogeneous weight
subspaces}\label{Linear algebra} 

Let
\[V = \coprod_{k \in \frac{1}{2} \mathbb{Z}} V_{(k)} = \coprod_{k \in
\frac{1}{2} \mathbb{Z}} V_{(k)}^0 \oplus \coprod_{k \in \frac{1}{2}
\mathbb{Z}} V_{(k)}^1 = V^0 \oplus V^1 \]  
with
\[\dim V_{(k)} < \infty \quad \mbox{for} \quad k \in \frac{1}{2}
\mathbb{Z} ,\] 
be a $\frac{1}{2} \mathbb{Z}$-graded (by weight) $\bigwedge_*$-module with
finite-dimensional homogeneous weight spaces $V_{(k)}$ which is also
$\mathbb{Z}_2$-graded.  Let 
\[V' = \coprod_{k \in \frac{1}{2} \mathbb{Z}} V_{(k)}^* \]
be the graded dual space of $V$, 
\[\bar{V} = \prod _{k \in \frac{1}{2} \mathbb{Z}} V_{(k)} = (V')^* \] 
the algebraic completion of $V$, and $\langle \cdot , \cdot \rangle$
the natural pairing between $V'$ and $\bar{V}$.  For any $n \in
\mathbb{N}$, let 
\[\mathcal{SF}_V(n) = \mbox{Hom}_{\bigwedge_*}(V^{\otimes n}, \bar{V})
.\] 
For any $m \in \Z$, $n \in \mathbb{N}$, and any positive integer $i \leq
m$, we define the {\it $t^{1/2}$-contraction} 
\begin{eqnarray*} 
 _i*_0 : \mathcal{SF}_V(m) \times \mathcal{SF}_V(n) & \rightarrow &
\mbox{Hom} (V^{\otimes (m + n - 1)}, V[[t^{\frac{1}{2}}, t^{-
\frac{1}{2}}]])  \\ 
(f,g) &\mapsto& (f \; _i*_0 \; g)_{t^{1/2}} , 
\end{eqnarray*}
by
\begin{multline}\label{t-contraction}
(f \; _i*_0 \; g)_{t^{1/2}} (v_1 \otimes \cdots \otimes v_{m + n -
1}) \\
= \sum_{k \in \frac{1}{2} \mathbb{Z}}f(v_1 \otimes \cdots \otimes
v_{i - 1} \otimes P(k)(g(v_i \otimes \cdots \otimes v_{i + n - 1}))
\otimes v_{i + n} \otimes \cdots \otimes v_{m + n - 1}) t^k 
\end{multline}
for all $v_1,...,v_{m + n - 1} \in V$, where for any $k \in
\frac{1}{2} \mathbb{Z}$, the map $P(k) : V \rightarrow V_{(k)}$ is the
canonical projection map.  This is an obvious extension of the notion of
$t$-contraction given in \cite{H book}.  

Let
\begin{equation}\label{V basis}
\{e^{(k)}_{l^{(k)}} \in V_{(k)} \; | \; k \in \frac{1}{2} \mathbb{Z}, 
\; l^{(k)} = 1,..., \dim V_{(k)} \} 
\end{equation} 
be a homogeneous basis of $V$, and
\begin{equation}\label{dual basis}
\{(e^{(k)}_{l^{(k)}})^* \in V_{(k)}^* \; | \; k \in \frac{1}{2} 
\mathbb{Z}, \; l^{(k)} = 1,..., \dim V_{(k)} \} 
\end{equation} 
the corresponding dual basis of $V'$.  In terms of these bases, the
definition of $t^{1/2}$-contraction can be written as
\begin{multline*}
(f \;  _i*_0 \; g)_{t^{1/2}} (v_1 \otimes \cdots \otimes v_{m + n -
1}) \\
= \sum_{k \in \frac{1}{2} \mathbb{Z}} \sum_{i^{(k)} = 1}^{\dim V_{(k)}}
f(v_1 \otimes \cdots \otimes v_{i - 1} \otimes e^{(k)}_{i^{(k)}}
\otimes v_{i + n} \otimes \cdots \otimes v_{m + n - 1}) \cdot \\
\cdot \langle (e^{(k)}_{i^{(k)}})^*, g(v_i \otimes \cdots \otimes v_{i + n - 1})
\rangle t^k 
\end{multline*}
for any $v_1,...,v_{m + n -1} \in V$.

If for arbitrary $v' \in V'$, $v_1,...,v_{m + n -1} \in V$, the formal
Laurent series in $t^{1/2}$
\[ \langle v', (f \; _i*_0 \; g)_{t^{1/2}}(v_1 \otimes \cdots \otimes
v_{m + n - 1}) \rangle \]
is absolutely convergent when $t^{1/2} = 1$, then $(f \; _i*_0
\; g)_1$ is well defined as an element of $\mathcal{SF}_V(m + n -1)$, and
we define  the {\it contraction $(f \; _i*_0 \; g)$ in $\mathcal{SF}_V(m
+ n -1)$ of $f$ and $g$} by  
\[(f \; _i*_0 \; g) = (f \;  _i*_0 \; g)_1 .\]
The following associativity of the $t^{1/2}$-contraction follows
immediately {}from the definition.    

\begin{prop}\label{contraction associativity}
Let $l \in \Z$, and $m,n \in \mathbb{N}$ such that $l + m - 1 \in \Z$,
and let $f_1 \in \mathcal{SF}_V(l)$, $f_2 \in \mathcal{SF}_V(m)$, $f_3 \in 
\mathcal{SF}_V(n)$, and $i,j \in \Z$ such that $1 \leq i \leq l$, and 
$1 \leq j \leq l + m -1$.  Then one of the following holds:

(i) $j<i$ and as formal power series in $t_1^{1/2}$ and
$t_2^{1/2}$  
\[((f_1 \; _i*_0 \; f_2)_{t_1^{1/2}} \; _j*_0 \; f_3)_{t_2^{1/2}} =
((f_1 \; _j*_0 \; f_3)_{t_2^{1/2}} \; _{i + n - 1}*_0 \; f_2
)_{t_1^{1/2}}; \]  

(ii) $j \geq i + m$ and as formal power series in $t_1^{1/2}$
and $t_2^{1/2}$ 
\[((f_1 \; _i*_0 \; f_2)_{t_1^{1/2}} \; _j*_0 \; f_3)_{t_2^{1/2}} =
((f_1 \; _{j - m + 1}*_0 \; f_3)_{t_2^{1/2}} \; _i*_0 \;
f_2)_{t_1^{1/2}} ; \]  

(iii) $i \leq j < i + m$ and as formal power series in
$t_1^{1/2}$ and $t_2^{1/2}$  
\[((f_1 \; _i*_0 \; f_2)_{t_1^{1/2}} \; _j*_0 \; f_3)_{t_2^{1/2}} =
(f_1 \; _i*_0 \; (f_2 \; _{j - i + 1}*_0 \;
f_3)_{t_2^{1/2}})_{t_1^{1/2}} .\]  
\end{prop}  

\begin{rema} {\em Note the similarity between Proposition
\ref{contraction associativity} and Proposition \ref{sewing associativity}. 
This similarity is due to the fact that both the sewing operation on the
moduli space of $N=1$ superspheres with tubes and the 
$t^{1/2}$-contraction define the formal substitution maps for a
partial operad structure (cf. \cite{M}, \cite{HL1},\cite{HL2}) on the  
moduli space of $N=1$ superspheres with tubes and on the $n$-ary 
functions {}from $V$ to $V[[t^{1/2}, t^{-1/2}]]$, respectively.  The fact
that the associativity of the contraction itself relies on a  
certain series being doubly absolutely convergent means that the 
contraction defines only a partial ``pseudo-operad" (cf. \cite{HL1},
\cite{HL2}) since operad associativity is not automatic.  

In addition, an operad has an action of the permutation group which we
define below for these $n$-ary functions. }    
\end{rema}

Let $(k \; l) \in S_{n}$ be the permutation on $n$ letters
which switches the $k$-th and $l$-th letters, for $k, l = 1, ..., n$,
with $k < l$.   We define the (right) action of the transposition $(k \;
l)$ on $V^{\otimes n}$ by 
\begin{multline*}
(v_1 \otimes \cdots \otimes v_{k-1} \otimes v_k \otimes v_{k+1} \otimes
\cdots \otimes v_{l-1} \otimes v_l  \otimes v_{l+1} \otimes \cdots \otimes
v_n) \cdot (k \; l)\\  
= (-1)^{\eta{(k \; l)}} (v_1 \otimes \cdots \otimes v_{k-1} \otimes v_l
\otimes  v_{k+1} \otimes \cdots \otimes v_{l-1} \otimes v_k  \otimes
v_{l+1} \otimes \cdots  \otimes v_n)
\end{multline*}
for $v_j$ of homogeneous sign in $V$, where 
\[\eta{(k \; l)} = \sum_{j = k + 1}^{l-1} \eta(v_j)(\eta(v_k) +
\eta(v_l)) + \eta(v_l)\eta(v_k) .\]
Let $\sigma \in S_n$ be a permutation on $n$ letters. Then $\sigma$ is
the product of transpositions $\sigma = \sigma_1 \cdots \sigma_m$,
$\sigma_i = (k_i \; l_i)$, for $k_i, l_i = 1,...,n$, with $k_i < l_i$, $i =
1,...,m$.   Thus we have a right action of $S_n$ on $V^{\otimes n}$
given by  
\begin{eqnarray*}
(v_1 \otimes \cdots \otimes v_n) \cdot \sigma &=& (v_1 \otimes \cdots 
\otimes  v_n) \cdot \sigma_1 \cdots
\sigma_m \cdot \\
&=& (-1)^{\eta(\sigma_1) + \cdots + \eta(\sigma_m)} v_{\sigma(1)} 
\otimes \cdots \otimes v_{\sigma(n)} \\ 
&=& (-1)^{\eta(\sigma)} v_{\sigma(1)} \otimes \cdots \otimes
v_{\sigma(n)}  .
\end{eqnarray*} 
This right action of $S_n$ on $n$-tuples in $V$, induces a left action
of $S_n$ on $\mathcal{SF}_V(n)$ given by 
\[(\sigma \cdot f) (v_1 \otimes \cdots \otimes v_n) = f((v_1
\otimes \cdots \otimes v_n) \cdot \sigma) , \]
for $f \in \mathcal{SF}_V(n)$.

Since $V$ is $\mathbb{Z}_2$-graded, End $V$ has 
a natural $\mathbb{Z}_2$-grading and a natural Lie superalgebra 
structure.  If $P \in \mbox{End} \; V$, the  
corresponding adjoint operator on $V'$, if it exists, is denoted by $P'$.  
The condition for the existence of $P'$ is that the linear functional on $V$ 
defined by the right-hand side of  
\[ \langle P'v',v\rangle = \langle v',Pv \rangle, \; \mbox{for} \; v
\in V, \; v' \in V' \]
should lie in $V'$.  If there exists $n \in \frac{1}{2} \mathbb{Z}$ such
that $P$ maps $V_{(k)}$ to $V_{(n + k)}$ for any $k \in \frac{1}{2}
\mathbb{Z}$, we say that $P$ has {\it weight $n$}.  Note that $P$ has 
weight $n$ if and only if its adjoint $P'$ exists and has weight $-n$ as 
an operator on $V'$.  In addition, $P$ is even (resp., odd) if and only if 
its adjoint is even (resp., odd).  In the case that $V$ is a module for 
$\mathfrak{ns}$ graded by the eigenvalues of $L(0)$, the adjoint 
operator $L'(-n)$ for $n \in \mathbb{Z}$ corresponding to $L(-n)$ exists and 
is even with weight $n$, and the adjoint operator $G'(-n - 1/2)$ for 
$n \in \mathbb{Z}$ corresponding to $G(-n - 1/2)$ exists and is odd 
with weight $n + 1/2$.
 
\subsection{The notion of $N=1$ supergeometric vertex operator 
superalgebra}\label{geom VOSAs}

\begin{defn} {\em An} $N = 1$ supergeometric vertex
operator superalgebra over $\bigwedge_*$ ($N=1$ SG-VOSA) {\em is a
$\frac{1}{2} \mathbb{Z}$-graded (by weight) $\bigwedge_\infty$-module
which is also  $\mathbb{Z}_2$-graded (by sign)  
\[V = \coprod_{k \in \frac{1}{2} \mathbb{Z}} V_{(k)} = \coprod_{k \in
 \frac{1}{2} \mathbb{Z}} V_{(k)}^0 \oplus \coprod_{k \in  \frac{1}{2}
\mathbb{Z}} V_{(k)}^1 = V^0 \oplus V^1 \]  
such that only $\bigwedge_* \subseteq \bigwedge_\infty$ acts nontrivially
on $V$,  
\[\dim V_{(k)} < \infty \quad \mbox{for} \quad k \in \frac{1}{2}
\mathbb{Z} ,\] 
and, for any $n \in \mathbb{N}$, a map
\[ \nu_n : SK(n) \rightarrow S\mathcal{F}_V (n) \]
satisfying the following axioms: 

(1) Positive energy axiom: 
\[V_{(k)} = 0 \quad \mbox{for $k$ sufficiently small.} \]

(2) Grading axiom:  Let $v' \in V'$, $v \in V_{(k)}$, and $a \in
(\bigwedge_\infty^0)^\times$.  Then 
\[ \langle v', \nu_1(\mathbf{0},(a, \mathbf{0}))(v)\rangle =
a^{-2k} \langle v',v \rangle .\]

(3) Supermeromorphicity axiom: For any $n \in \Z$, $v' \in V'$,
and $v_1,...,v_n \in V$, the function
\[Q \mapsto \langle v', \nu_n(Q) (v_1 \otimes \cdots \otimes v_n)
\rangle \]
on $SK(n)$ is a canonical supermeromorphic superfunction (in the
sense of (\ref{meromorphic})), and if $(z_i, \theta_i)$ and $(z_j,
\theta_j)$, for $i,j \in \{1,...,n\}$, $i \neq j$, are the $i$-th and
$j$-th punctures of $Q \in SK(n)$, respectively,  then for any $v_i$
and $v_j$ in $V$, there exists $N(v_i,v_j) \in \Z$ such that for any
$v' \in V'$ and $v_l \in V$, $l \neq i,j$, the order of the pole
$(z_i, \theta_i) = (z_j, \theta_j)$ of $\langle v', \nu_n(Q) (v_1
\otimes \cdots \otimes v_n) \rangle$ is less then $N(v_i,v_j)$.  

(4) Permutation axiom:  Let $\sigma \in S_n$.  Then for any $Q
\in SK(n)$
\[\sigma \cdot (\nu_n (Q)) = \nu_n(\sigma \cdot Q) . \]

(5) Sewing axiom:  There exists a unique complex number $c$ (the}
central charge {\em or} rank{\em) such that if $Q_1 \in SK(m)$ and
$Q_2 \in SK(n)$ are given by 
\begin{multline*}
Q_1 = ((z_1, \theta_1),...,(z_{m-1}, \theta_{m-1}); (A^{(0)},
M^{(0)}), (a^{(1)}, A^{(1)}, M^{(1)}),\\
...,(a^{(m)}, A^{(m)}, M^{(m)})), 
\end{multline*} 
and
\[Q_2 = ((z_1', \theta_1'),...,(z_{n-1}', \theta_{n-1}'); (B^{(0)},
N^{(0)}), (b^{(1)}, B^{(1)}, N^{(1)}),...,(b^{(n)}, B^{(n)},
N^{(n)})),\] 
and if the $i$-th tube of $Q_1$ can be sewn with the 0-th tube of $Q_2$,
then for any $v' \in V'$ and $v_1,...,v_{m + n - 1} \in V$, 
\[\langle v', (\nu_m (Q_1) \; _i*_0 \; \nu_n (Q_2))_{t^{1/2}} (v_1
\otimes \cdots \otimes v_{m + n - 1}) \rangle \]
is absolutely convergent when $t^{1/2} = 1$, and
\[\nu_{m + n - 1} (Q_1 \; _i\infty_0 \; Q_2) = (\nu_m (Q_1) \;
_i*_0  \; \nu_n (Q_2)) \; e^{-\Gamma (a^{(i)}, A^{(i)}, M^{(i)},
B^{(0)}, N^{(0)})c} . \] }
\end{defn}

Note that in this definition, we assume that Conjecture \ref{Gamma
converges} is true.

We denote the $N=1$ SG-VOSA defined above by 
\[(V, \nu = \{\nu_n \}_{n \in \mathbb{N}})\]
or just by $(V, \nu)$.   

We can extend this definition to Grassmann algebras over an arbitrary field
$\mathbb F$ of characteristic zero.  But  in this case, we have to first
introduce the notion of formal superspheres  with tubes over the Grassmann
algebra over $\mathbb F$, and the sewing  axiom will become a statement for
certain formal series in even and odd  variables with coefficients in this
Grassmann algebra.

Let $(V_1, \nu)$ and $(V_2, \mu)$ be two $N=1$ SG-VOSAs over
$\bigwedge_{*_1}$ and $\bigwedge_{*_2}$, respectively.  A {\it
homomorphism of $N=1$ SG-VOSAs} is a doubly graded
$\bigwedge_\infty$-module map $\gamma : V_1 \rightarrow V_2$ such that if
$\bar{\gamma} : \bar{V_1} \rightarrow \bar{V_2}$ is the unique  extension
of $\gamma$, then for any $n \in \mathbb{N}$ and any $Q \in  SK(n)$ 
\[\bar{\gamma} \circ \nu (Q) = \mu (Q) \circ \gamma^{\otimes
n} . \]

Note that an $N=1$ SG-VOSA over some subalgebra of the underlying
Grassmann algebra $\bigwedge_\infty$ is isomorphic to the same $N=1$
SG-VOSA over a different subalgebra.

\section{(Algebraic) $N=1$ Neveu-Schwarz vertex operator superalgebras}

In this section, we recall the notion of {\it $N = 1$ Neveu-Schwarz 
vertex operator superalgebra over a Grassmann algebra and with odd 
formal variables ($N=1$ NS-VOSA)} given in \cite{B thesis} and \cite{B
vosas} and  recall some of the consequences of this notion.  In addition,
in this section we give the formulations of supercommutativity and
associativity for an $N=1$ NS-VOSA with odd formal variables.

\subsection{Delta functions with odd formal variables}

Let $x$, $x_0$, $x_1$ and $x_2$ be even formal variables, and
let $\varphi$, $\varphi_0$, $\varphi_1$ and $\varphi_2$ be odd
formal variables.  For any formal Laurent series $f(x) \in  
\bigwedge_\infty[[x,x^{-1}]]$, we can define 
\begin{equation}\label{expansion}
f(x + \varphi_1 \varphi_2) = f(x) + \varphi_1
\varphi_2 f'(x) \; \in \mbox{$\bigwedge_\infty$} [[x,x^{-1}]] [\varphi_1]
[\varphi_2]. 
\end{equation}

Recall (cf. \cite{FLM}) the {\it formal $\delta$-function at $x=1$} given
by $\delta(x) = \sum_{n \in \mathbb{Z}} x^n$. 
As developed in \cite{B vosas}, we have the following $\delta$-function of 
expressions involving three even formal variables and two odd formal 
variables
\begin{eqnarray*}
\delta \biggl( \frac{x_1 - x_2 - \varphi_1 \varphi_2}{x_0} \biggr) &=&
\sum_{n \in \mathbb{Z}} (x_1 - x_2 - \varphi_1 \varphi_2)^n x_0^{-n} \\
&=& \delta \biggl( \frac{x_1 - x_2}{x_0} \biggr)  -
\varphi_1 \varphi_2 x_0^{-1} \delta' \biggl( \frac{x_1 - x_2}{x_0}
\biggr)
\end{eqnarray*} 
where $\delta'(x) = d/dx \; \delta (x) = \sum_{n \in \mathbb{Z}} n
x^{n-1}$, and we use the conventions that a function of even and odd
variables should be expanded about the even variables and any expression in
two even variables (such as $(x_1 - x_2)^n$, for $n \in \mathbb{Z}$) should
be expanded in positive powers of the second variable, (in this case $x_2$). 
        
{}From \cite{B vosas}, we have the following two $\delta$-function
identities
\begin{equation}\label{delta 2 terms with phis}
x_1^{-1} \delta \biggl( \frac{x_2 + x_0 + \varphi_1 \varphi_2}{x_1}
\biggr) = x_2^{-1} \delta \biggl( \frac{x_1 - x_0 - \varphi_1
\varphi_2}{x_2} \biggr)  
\end{equation}
\begin{multline}\label{delta 3 terms with phis}
x_0^{-1} \delta \biggl( \frac{x_1 - x_2 - \varphi_1 \varphi_2}{x_0}
\biggr) - x_0^{-1} 
\delta \biggl( \frac{x_2 - x_1 + \varphi_1 \varphi_2}{-x_0} \biggr) \\
= x_2^{-1} \delta \biggl( \frac{x_1 - x_0 - \varphi_1 
\varphi_2}{x_2} \biggr) . 
\end{multline}
The expressions on both sides of (\ref{delta 2 terms with phis}) and 
the three terms occurring in (\ref{delta 3 terms with phis}) all 
correspond to the same formal substitution $x_0 = x_1 - x_2  - \varphi_1 
\varphi_2$.  Letting $V$ be any $\bigwedge_\infty$-module and for
\[X(x_0,\varphi_0,x_1,\varphi_1,x_2,\varphi_2) \in (\mbox{End} \;V)
[[x_0,x_0^{-1},x_1,x_1^{-1},x_2,x_2^{-1}]] [\varphi_1, \varphi_2] ,\] 
a formal substitution corresponding to $x_0 = x_1 - x_2 - \varphi_1
\varphi_2$ can be made as long as the resulting expression is
well defined, e.g., if 
\[X(x_1,\varphi_1,x_2,\varphi_2) \in (\mbox{End} \;V)
[[x_1,x_1^{-1}]]((x_2)) [\varphi_1, \varphi_2] ,\] 
then  
\begin{multline}\label{delta substitute}
\delta \biggl( \frac{x_2 + x_0 + \varphi_1 \varphi_2}{x_1} \biggr)
X(x_1,\varphi_1,x_2,\varphi_2) \\
= \delta \biggl( \frac{x_2 + x_0 + \varphi_1
\varphi_2}{x_1} \biggr) X(x_2 + x_0 + \varphi_1 \varphi_2,
\varphi_1,x_2,\varphi_2) ,
\end{multline}
and (\ref{delta substitute}) holds for more general $X(x_1, \varphi_1,
x_2, \varphi_2)$ as long as the expression $X(x_2 + x_0 + \varphi_1
\varphi_2, \varphi_1,x_2,\varphi_2)$ is well defined.  Any such
substitution corresponding to $x_0 = x_1 - x_2 - \varphi_1 \varphi_2$ can
be thought of as the even part of a superconformal shift of $(x_1,
\varphi_1)$ by
$(x_2,\varphi_2)$.

\subsection{The notion of $N=1$ NS-VOSA over $\bigwedge_*$ and with
odd formal variables}

\begin{defn}\label{VOSA definition}
{\em An} $N = 1$ NS-VOSA over $\bigwedge_*$ and with odd variables
{\em is a $\frac{1}{2} \mathbb{Z}$-graded (by weight)
$\bigwedge_\infty$-module which is  also $\mathbb{Z}_2$-graded (by sign)  
\begin{equation}\label{vosa1}
V = \coprod_{k \in \frac{1}{2} \mathbb{Z}} V_{(k)} = \coprod_{k
\in \frac{1}{2}\mathbb{Z}} V_{(k)}^0 \oplus \coprod_{k \in
\frac{1}{2}\mathbb{Z}} V_{(k)}^1 = V^0 \oplus V^1 
\end{equation}  
such that only $\bigwedge_* \subseteq \bigwedge_\infty$ acts
nontrivially, 
\begin{equation}\label{vosa2}
\dim V_{(k)} < \infty \quad \mbox{for} \quad k \in \frac{1}{2}
\mathbb{Z} , 
\end{equation}
\begin{equation}\label{positive energy}
V_{(k)} = 0 \quad \mbox{for $k$ sufficiently small} , 
\end{equation}
equipped with a linear map $V \otimes V \longrightarrow V[[x,x^{-1}]]
[\varphi]$, or equivalently,
\begin{eqnarray*} 
V &\longrightarrow&  (\mbox{End} \; V)[[x,x^{-1}]][\varphi] \\
v  &\mapsto&  Y(v,(x,\varphi)) = \sum_{n \in \mathbb{Z}} v_n x^{-n-1} +
\varphi \sum_{n \in \mathbb{Z}} v_{n - \frac{1}{2}} x^{-n-1}
\end{eqnarray*}
where $v_n \in (\mbox{End} \; V)^{\eta(v)}$ and $v_{n - 1/2} \in 
(\mbox{End} \; V)^{(\eta(v) + 1) \mbox{\begin{footnotesize} mod 
\end{footnotesize}} 2}$ for $v$ of homogeneous sign in $V$, $x$ is an 
even formal variable, and $\varphi$ is an odd formal variable, and 
where $Y(v,(x,\varphi))$ denotes the} vertex operator associated with 
$v$, {\em and equipped also with two distinguished homogeneous vectors 
$\mathbf{1} \in V_{(0)}^0$ (the {\em vacuum}) and $\tau \in V_{(3/2)}^1$ 
(the {\em Neveu-Schwarz element}).  The following conditions are assumed 
for $u,v \in V$:   
\begin{equation}\label{truncation}
u_n v = 0 \quad \mbox{for $n \in \frac{1}{2} \mathbb{Z}$ sufficiently
large;} 
\end{equation}
\begin{equation}\label{vacuum identity}
Y(\mathbf{1}, (x, \varphi)) = 1 \quad \mbox{(1 on the right being
the identity operator);} 
\end{equation}
the} creation property {\em holds:
\[Y(v,(x,\varphi)) \mathbf{1} \in V[[x]][\varphi] \qquad \mbox{and}
\qquad \lim_{(x,\varphi) \rightarrow 0} Y(v,(x,\varphi)) \mathbf{1} =
v ; \] 
the} Jacobi identity {\em holds:  
\begin{multline*}
x_0^{-1} \delta \biggl( \frac{x_1 - x_2 - \varphi_1 \varphi_2}{x_0}
\biggr) Y(u,(x_1, \varphi_1))Y(v,(x_2, \varphi_2)) \\
- (-1)^{\eta(u)\eta(v)} x_0^{-1} \delta \biggl( \frac{x_2 - x_1 + 
\varphi_1 \varphi_2}{-x_0} \biggr)Y(v,(x_2, \varphi_2))Y(u,(x_1,
\varphi_1)) \\
= \; x_2^{-1} \delta \biggl( \frac{x_1 - x_0 - \varphi_1
\varphi_2}{x_2} \biggr) Y(Y(u,(x_0, \varphi_1 - \varphi_2))v,(x_2,
\varphi_2)) , 
\end{multline*}
for $u,v$ of homogeneous sign in $V$; the $N=1$ Neveu-Schwarz algebra 
relations hold:
\begin{eqnarray*}
\left[L(m),L(n) \right] \! &=& \!(m - n)L(m + n) + \frac{1}{12} (m^3 - m)
\delta_{m + n , 0} (\mbox{rank} \; V) , \\ \label{V9}
\biggl[ G(m + \frac{1}{2}),L(n) \biggr] \! &=& \! (m - \frac{n - 1}{2} ) G(m
+ n + \frac{1}{2}) ,\\ \label{V10}
\biggl[ G (m + \frac{1}{2} ) , G(n - \frac{1}{2} ) \biggr] \! &=& \! 2L(m +
n) + \frac{1}{3} (m^2 + m) \delta_{m + n , 0} (\mbox{rank} \; V) , \label{V11}
\end{eqnarray*}
for $m,n \in \mathbb{Z}$, where 
\[G(n + \frac{1}{2}) = \tau_{n + 1}, \qquad \mbox{and} \qquad 2L(n) =
\tau_{n + \frac{1}{2}} \qquad \mbox{for} \; n \in \mathbb{Z} , \]
i.e., 
\begin{equation}\label{stress tensor}
Y(\tau,(x,\varphi)) = \sum_{n \in \mathbb{Z}} G (n + \frac{1}{2}) x^{-
n - \frac{1}{2} - \frac{3}{2}} \; + \; 2 \varphi \sum_{n \in \mathbb{Z}}
L(n) x^{- n - 2} ,
\end{equation}
and $\mbox{rank} \; V \in \mathbb{C}$; 
\begin{equation}\label{grading for vosa with}
L(0)v = kv \quad \mbox{for} \quad k \in \frac{1}{2} \mathbb{Z} \quad
\mbox{and} \quad v \in V_{(k)}; 
\end{equation}
and the} $G(-1/2)$-derivative property {\em holds:
\begin{equation}\label{G(-1/2)-derivative}
\biggl( \frac{\partial}{\partial \varphi} + \varphi
\frac{\partial}{\partial x} \biggr) Y(v,(x,\varphi)) =  Y(G(- 
\frac{1}{2})v,(x,\varphi)) . 
\end{equation} }
\end{defn}

The $N=1$ NS-VOSA just defined is denoted
by $(V,Y(\cdot,(x,\varphi)),\mathbf{1},\tau)$, or for simplicity by $V$.

We have the following consequences of the definition of $N=1$ NS-VOSA with
odd formal variables.  We have 
\[ L(n) \mathbf{1} = G(n + \frac{1}{2}) \mathbf{1} = 0, \quad \mbox{for
$n \geq -1$}.\]
There exists $\omega = (1/2)G(-1/2) \tau \in V_{(2)}$ such that 
\begin{equation}
Y(\omega,(x,\varphi)) = \sum_{n \in \mathbb{Z}} L(n) x^{-n-2} -
\frac{\varphi}{2} \sum_{n \in \mathbb{Z}} (n + 1) G(n -
\frac{1}{2}) x^{-n-2} .
\end{equation}
The supercommutator formula is given by 
\begin{multline}\label{bracket relation for a vosa}
[ Y(u, (x_1,\varphi_1)), Y(v,(x_2,\varphi_2))] \\
= \mbox{Res}_{x_0} x_2^{-1} \delta \biggl( \frac{x_1
- x_0 - \varphi_1 \varphi_2}{x_2} \biggr) Y(Y(u,(x_0, \varphi_1 -
\varphi_2))v,(x_2, \varphi_2))  
\end{multline}
where $\mbox{Res}_{x_0}$ of a power series in $x_0$ is the coefficient of 
$x_0^{-1}$. 
We have
\begin{eqnarray}
x_0^{2L(0)} Y(v, (x, \varphi)) x_0^{-2 L(0)} = Y(x_0^{2L(0)}v, (x_0^2 x,
x_0 \varphi)), \label{conjugate by L(0)}
\end{eqnarray}
\begin{eqnarray}
Y(e^{x_0 L(-1) + \varphi_0 G(-\frac{1}{2})}v, (x,\varphi)) &=& e^{x_0
\frac{\partial}{\partial x} + \varphi_0 \left(
\frac{\partial}{\partial \varphi} + \varphi \frac{\partial}{\partial 
x} \right)} Y(v,(x,\varphi)) \\ 
&=& Y(v,(x_0 + x + \varphi_0 \varphi, \varphi_0 + \varphi)) 
\label{exponential L(-1) and G(-1/2) property}
\end{eqnarray}
and
\begin{eqnarray}
& & \hspace{-2.2in} e^{x_0 L(-1) + \varphi_0 G(-\frac{1}{2})} Y(v,(x,\varphi)) e^{- x_0
L(-1) - \varphi_0 G(-\frac{1}{2})} =  \nonumber\\
\hspace{1.9in} &=& Y(e^{x_0 L(-1) + \varphi_0
G(-\frac{1}{2}) - 2\varphi_0 \varphi L(-1)}v,(x,\varphi)) 
\label{L(-1) and G(-1/2) exp 1}\\
&=& Y(v,(x +x_0  + \varphi \varphi_0, \varphi + \varphi_0)). 
\label{L(-1) and G(-1/2) exp 2}
\end{eqnarray} 
And finally, we have the property of {\it skew-supersymmetry}: for $u,v$
of homogeneous sign in $V$
\begin{equation}\label{skew-supersymmetry}
e^{x L(-1) + \varphi G(-\frac{1}{2})} Y(v,(-x,-\varphi))u \; = \;
(-1)^{\eta(v) \eta(u)} Y(u,(x,\varphi))v.
\end{equation}

Let $(V_1, Y_1(\cdot,(x,\varphi)),\mathbf{1}_1,\tau_1)$ and $(V_2,
Y_2(\cdot,(x,\varphi)),\mathbf{1}_2,\tau_2)$ be two $N=1$ NS-VOSAs over
$\bigwedge_{*_1}$ and $\bigwedge_{*_2}$, respectively.  A {\it
homomorphism of $N=1$ NS-VOSAs} is a doubly graded
$\bigwedge_\infty$-module homomorphism $\gamma : V_1 \longrightarrow V_2
\;$ such that 
\[\gamma (Y_1(u,(x,\varphi))v) = Y_2(\gamma(u),(x,\varphi))\gamma(v)
\quad \mbox{for} \quad u,v \in V_1 ,\]
$\gamma(\mathbf{1}_1) = \mathbf{1}_2$, and $\gamma(\tau_1) =
\tau_2$.  

Of course an $N=1$ NS-VOSA over
$\bigwedge_*$ is isomorphic to the same $N=1$ NS-VOSA with coefficients
in any other $\bigwedge_{*'}$. 

\begin{rema}{\em
Let $(V,Y(\cdot,(x,\varphi)), \mathbf{1}, \tau)$ be an $N=1$ NS-VOSA.  In
\cite{B vosas}, we also study the notion of $N=1$ NS-VOSA over a
Grassmann algebra without odd formal variables and show that
$(V,Y(\cdot,(x,0)), \mathbf{1}, \tau)$ is such an algebra.  Conversely,
given an $N=1$ NS-VOSA without odd formal variables $(V,Y(\cdot,x),
\mathbf{1}, \tau)$,  we can define $\tilde{Y}(v, (x,\varphi)) = Y(v,x) +
\varphi Y(G(-1/2) v,x)$,  and then $(V,\tilde{Y}(\cdot,(x,\varphi)),
\mathbf{1}, \tau)$ is an $N=1$ NS-VOSA with odd formal variables.  Using
this correspondence, in \cite{B vosas} we prove that the category of
$N=1$ NS-VOSAS over $\bigwedge_*$ with odd formal variables and with
central charge, i.e., rank, $c \in \mathbb{C}$, is isomorphic to the
category  of $N=1$ NS-VOSAs over $\bigwedge_*$ without odd formal 
variables and with central charge $c \in \mathbb{C}$.   However, in
including the odd formal variables the correspondence with the geometry
and the role of the operator $G(-1/2)$ as in the $G(-1/2)$-derivative
property is more  explicit.  }
\end{rema}

\subsection{Expansions of rational superfunctions}\label{iota}

Let $\bigwedge_\infty[x_1,x_2,...,x_n]_S$ be the ring of rational
functions obtained by inverting (localizing with respect to) the set 
\[S = \biggl\{\sum_{i = 1}^{n} a_i x_i \; | \; a_i \in
\mbox{$\bigwedge_\infty^0$},
\; \mbox{not all} \; (a_i)_B = 0\biggr\} . \] 
Recall the map $\iota_{i_1 ... i_n} : \mathbb{F}[x_1,...,x_n]_S
\longrightarrow \mathbb{F}[[x_1, x_1^{-1},...,  x_n, x_n^{-1}]]$ defined
in \cite{FLM} where coefficients of elements in $S$ are restricted to the
field $\mathbb{F}$.  We extend this map to
$\bigwedge_\infty[x_1,x_2,...,x_n]_S[\varphi_1,\varphi_2,...,\varphi_n] =
\bigwedge_\infty[x_1,\varphi_1,x_2,\varphi_2,...,x_n,\varphi_n]_{S}$ in
the obvious way obtaining
\[\iota_{i_1 ... i_n} : \mbox{$\bigwedge_\infty$}
[x_1,\varphi_1,...,x_n,\varphi_n]_S
\longrightarrow \mbox{$\bigwedge_\infty$}[[x_1, x_1^{-1},..., x_n,
x_n^{-1}]][\varphi_1,...,\varphi_n] .\]
Let $\bigwedge_\infty[x_1, \varphi_1, x_2, \varphi_2,...,x_n,
\varphi_n]_{S'}$ be the ring of rational functions obtained by inverting
the set   
\[S' = \biggl\{\sum_{\stackrel{i,j = 1}{i<j}}^{n} (a_i x_i + a_{ij} \varphi_i
\varphi_j) \; | \; a_i, a_{ij} \in \mbox{$\bigwedge_\infty^0$}, \;
\mbox{not all} \; (a_i)_B = 0\biggr\}. \]
Since we use the convention that a function of even and odd
variables should be expanded about the even variables, we have 
\[\frac{1}{\sum_{\stackrel{i,j = 1}{i<j}}^{n} (a_i x_i + a_{ij} \varphi_i
\varphi_j)} = \frac{1}{\sum_{i = 1}^{n} a_i x_i} -
\frac{\sum_{\stackrel{i,j = 1}{i<j}}^{n} a_{ij} \varphi_i
\varphi_j}{(\sum_{i = 1}^{n} a_i x_i)^2} .\]
Thus $\bigwedge_\infty[x_1,\varphi_1,x_2,\varphi_2,...,x_n,\varphi_n]_{S'}
\subseteq
\bigwedge_\infty[x_1,\varphi_1,x_2,\varphi_2,...,x_n,\varphi_n]_S$,    and
$\iota_{i_1 ... i_n}$ is well defined on
$\bigwedge_\infty[x_1,\varphi_1,x_2,\varphi_2,...,x_n,\varphi_n]_{S'}$.

In the case $n = 2$, the map $\iota_{12} : \bigwedge_\infty
[x_1,\varphi_1,x_2,\varphi_2]_{S'} \longrightarrow \bigwedge_\infty
[[x_1,x_2]][\varphi_1,\varphi_2]$ is given by first expanding an
element of $\bigwedge_\infty [x_1,\varphi_1,x_2,\varphi_2]_{S'}$ as a
formal series in $\bigwedge_\infty[x_1,\varphi_1,x_2,\varphi_2]_S$ and
then expanding each term as a series in $\bigwedge_\infty [[x_1,x_2]]
[\varphi_1,\varphi_2]$ containing at most finitely many negative powers
of $x_2$ (using binomial expansions for negative powers of linear
polynomials involving both $x_1$ and $x_2$).

\subsection{Duality for vertex operator superalgebras}\label{Duality}

In Section \ref{geom VOSAs}, we formulated the notion of $N=1$ SG-VOSA.  In
Section \ref{alg {}from geom} we will show that any such object defines an
$N=1$ NS-VOSA with odd formal variables.  To show that the alleged $N=1$
NS-VOSA satisfies the Jacobi identity, we need the notions of associativity
and supercommutativity for an $N=1$ NS-VOSA with odd formal variables as
formulated in \cite{B thesis} and \cite{B vosas}.  Together, these notions
of associativity and (super)commutativity are known as ``duality'', a term
which arose {}from physics. 

Let $(V, Y(\cdot,(x,\varphi)), \mathbf{1}, \tau)$ be an $N=1$ NS-VOSA. 
As in Section \ref{Linear algebra}, let  $V'$ denote the graded dual
space of $V$, let $\bar{V}$ denote the algebraic completion of $V$, and
let $\langle \cdot , \cdot \rangle$ denote the natural pairing between
$V'$ and $\bar{V}$.  

\begin{prop}\label{supercommutativity} (\cite{B vosas})
{\bf (a) (rationality of products)}  For $u, v, w \in V$, with $u$, and $v$
of homogeneous sign, and $v' \in V'$, the formal series 
\[ \langle v', Y(u,(x_1,\varphi_1)) Y(v,(x_2,\varphi_2)) w \rangle, \]
which involves
only finitely many negative powers of $x_2$ and only finitely many
positive powers of $x_1$, lies in the image of the map $\iota_{12}$: 
\[\langle v', Y(u,(x_1,\varphi_1)) Y(v,(x_2,\varphi_2)) w \rangle
= \iota_{12} f(x_1,\varphi_1,x_2,\varphi_2) , \]
where the (uniquely determined) element $f \in \bigwedge_\infty
[x_1,\varphi_1,x_2,\varphi_2]_{S'}$ is of the form  
\[f(x_1,\varphi_1,x_2,\varphi_2) =
\frac{g(x_1,\varphi_1,x_2,\varphi_2)}{x_1^r x_2^s (x_1 - x_2 - 
\varphi_1 \varphi_2)^t} \]
for some $g \in \bigwedge_\infty [x_1,\varphi_1,x_2,\varphi_2]$ and $r, s,
t \in \mathbb{Z}$.  

{\bf (b) (supercommutativity)}  We also have
\[\langle v', Y(v,(x_2,\varphi_2))Y(u,(x_1,\varphi_1)) w \rangle =
(-1)^{\eta(u)\eta(v)} \iota_{21} f(x_1,\varphi_1,x_2,\varphi_2)
, \] 
i.e,
\begin{multline*}
\iota_{12}^{-1} \langle v', Y(u,(x_1,\varphi_1))
Y(v,(x_2,\varphi_2)) w \rangle \\
= (-1)^{\eta(u)\eta(v)} \iota_{21}^{-1} \langle v',
Y(v,(x_2,\varphi_2)) Y(u,(x_1,\varphi_1)) w \rangle .
\end{multline*}  
\end{prop}

\begin{prop} \label{iterates}(\cite{B vosas})
{\bf (a) (rationality of iterates)} For $u, v, w \in V$, and $v' \in
V'$, the formal series $\langle v', Y(Y(u,(x_0,\varphi_1 -
\varphi_2))v,(x_2,\varphi_2)) w \rangle$, which involves only finitely
many negative powers of $x_0$ and only finitely many positive powers
of $x_2$, lies in the image of the map $\iota_{20}$: 
\[\langle v', Y(Y(u,(x_0,\varphi_1 - \varphi_2))v,(x_2,\varphi_2))w
\rangle  = \iota_{20} h(x_0,\varphi_1 - \varphi_2,x_2,\varphi_2) , \]
where the (uniquely determined) element $h \in \bigwedge_\infty
[x_0,\varphi_1,x_2,\varphi_2]_{S'}$ is of the form      
\[h(x_0,\varphi_1 - \varphi_2,x_2,\varphi_2) = \frac{k(x_0,\varphi_1 -
\varphi_2,x_2,\varphi_2)}{x_0^r x_2^s (x_0 + x_2 - \varphi_1
\varphi_2)^t} \]   
for some $k \in \bigwedge_\infty [x_0,\varphi_1,x_2,\varphi_2]$ and $r, s,
t \in \mathbb{Z}$.  

{\bf (b)}  The formal series $\langle v', Y(u,(x_0 + x_2  + \varphi_1
\varphi_2, \varphi_1)) Y(v,(x_2,\varphi_2)) w \rangle$, 
which involves only finitely many negative powers of $x_2$ and only
finitely many positive powers of $x_0$, lies in the image of
$\iota_{02}$, and in fact  
\[\langle v', Y(u,(x_0 + x_2  + \varphi_1 \varphi_2, \varphi_1))
Y(v,(x_2,\varphi_2)) w \rangle =  \iota_{02} h(x_0,\varphi_1 -
\varphi_2, x_2, \varphi_2) . \]    
\end{prop}

\begin{prop}\label{associativity} (\cite{B vosas})
{\bf (associativity)}  We have the following equality of rational
functions: 
\begin{multline*}
\iota^{-1}_{12} \langle v', Y(u,(x_1,\varphi_1))
Y(v,(x_2,\varphi_2)) w \rangle \\
= \left. \left( \iota^{-1}_{20} \langle v',
Y(Y(u,(x_0,\varphi_1 - \varphi_2))v,(x_2,\varphi_2)) w \rangle \right)
\right|_{x_0 = x_1 - x_2 - \varphi_1 \varphi_2} 
\end{multline*}
\end{prop}

\begin{prop}\label{duality}(\cite{B vosas})
In the presence of the other axioms in the definition of $N=1$ NS-VOSA
with odd variables, the Jacobi identity follows {}from the rationality of
products and iterates and supercommutativity and associativity.  In
particular, the Jacobi identity may be replaced by these properties. 
\end{prop}

\section{$N=1$ NS-VOSAs {}from $N=1$ SG-VOSAs}\label{alg {}from
geom}

In this section, given an $N=1$ supergeometric VOSA, we construct an
(algebraic) $N=1$ NS-VOSA.  Since $N=1$ NS- (resp., SG-) VOSAs are
isomorphic if taken over different subalgebras of the underlying
Grassmann algebra $\bigwedge_\infty$, we suppress the details of which
subalgebras are actually acting non-trivially.

Let $(V,\nu)$ be an $N=1$ SG-VOSA, we define
the vacuum $\mathbf{1}_\nu \in \bar{V}$ by 
\[\mathbf{1}_\nu = \nu_0 (\mathbf{0}) ; \]
an element $\tau_\nu \in \bar{V}$ by
\[\tau_\nu = - \frac{\partial}{\partial \epsilon} \nu_0 (\mathbf{0},
M(\epsilon, \frac{3}{2})) ,\]
(recall (\ref{define tau series})); and the vertex operator  
\[Y_\nu (v_1, (x, \varphi)) = \sum_{n \in \mathbb{Z}} (v_1)_n x^{-n-1} +
\varphi \sum_{n \in \mathbb{Z}} (v_1)_{n - 1/2} x^{-n-1}\] 
associated with $v_1 \in V$ by
\begin{equation}\label{define Y}
(v_1)_n v_2 + \theta (v_1)_{n - \frac{1}{2}} v_2 = \mbox{Res}_z
\left( z^n \nu_2 ((z,\theta); \mathbf{0}, (1, \mathbf{0}), (1,
\mathbf{0}) ) (v' \otimes v_1 \otimes v_2) \right) , 
\end{equation}
where $\mbox{Res}_z$ means taking the residue at the singularity $z =
0$, i.e., taking the coefficient of $z^{-1}$.

\begin{prop}\label{get a vosa}  
The elements $\mathbf{1}_\nu$ and $\tau_\nu$ of $\bar{V}$ are in
fact in $V_{(0)}$ and $V_{(\frac{3}{2})}$, respectively.  If the rank
of $(V,\nu)$ is $c$, then $(V, Y_\nu( \cdot, (x,\varphi)),
\mathbf{1}_\nu, \tau_\nu)$ is an $N=1$ NS-VOSA with odd formal variables and
with rank $c$ . 
\end{prop}

\begin{proof} To see that $\mathbf{1}_\nu \in
V_{(0)}$, we use the grading and sewing axioms for $(V,\nu)$ to note
that for all $a \in (\bigwedge_\infty^0)^\times$
\begin{eqnarray*}
\langle v', \mathbf{1}_\nu \rangle &=& \langle v', \nu_0 (\mbox{\bf
0}) \rangle\\  
&=& \langle v', \nu_0 ((\mathbf{0}, (a, \mathbf{0})) \; _1\infty_0
\; (\mathbf{0})) \rangle \\ 
&=& \langle v', (\nu_1 (\mathbf{0}, (a, \mathbf{0})) \; _1*_0 \;
\nu_0(\mathbf{0})) \rangle e^{- \Gamma(a,\mathbf{0},\mathbf{0})
c} \\ 
&=& \sum_{k \in \frac{1}{2} \mathbb{Z}} \langle v', \nu_1 (\mathbf{0},
(a, \mathbf{0})) (P(k) (\mathbf{1}_\nu)) \rangle \\
&=& \sum_{k \in \frac{1}{2} \mathbb{Z}} a^{-2k} \langle v', P(k)
(\mathbf{1}_\nu) \rangle . 
\end{eqnarray*}
Thus $\mathbf{1}_\nu \in V_{(0)}$. Similarly for $\tau_\nu$, for all $a
\in (\bigwedge_\infty^0)^\times$,
\begin{eqnarray*}
a^{-3} \langle v', \tau_\nu \rangle &=& - \langle v' ,
\frac{\partial}{\partial \epsilon} \nu_0 (\mathbf{0},
M(a^{-3} \epsilon, \frac{3}{2})) \rangle \\
&=& - \langle v', \frac{\partial}{\partial \epsilon} \nu_0 ((\mbox{\bf
0}, (a ,\mathbf{0})) \; _1\infty_0 \; (\mathbf{0}, M(\epsilon,
\frac{3}{2}))) \rangle \\
&=& - \langle v', \frac{\partial}{\partial \epsilon} (\nu_1 (\mbox{\bf
0}, (a ,\mathbf{0})) \; _1*_0 \; \nu_0(\mathbf{0}, M(\epsilon,
\frac{3}{2}))) \rangle e^{- \Gamma(a,\mathbf{0},\mathbf{0},
\mathbf{0},M(\epsilon, \frac{3}{2}) )c} \\ 
&=& - \sum_{k \in \frac{1}{2} \mathbb{Z}} \langle v', \nu_1 (\mathbf{0},
(a ,\mathbf{0})) (P(k) (\frac{\partial}{\partial \epsilon} \nu_0
(\mathbf{0}, M(\epsilon, \frac{3}{2})))) \rangle \\
&=& \sum_{k \in \frac{1}{2} \mathbb{Z}} a^{-2k} \langle v', P(k)
(\tau_\nu) \rangle .
\end{eqnarray*}
Thus $\tau_\nu \in V_{(\frac{3}{2})}$.

The axioms (\ref{vosa1}), (\ref{vosa2}), and (\ref{positive energy})
are obvious.   Since for any $v' \in V'$ and $v_1, v_2 \in V$, the
superfunction $\langle v', \nu_2 ((z,\theta); \mathbf{0}, (1, \mathbf{0}), 
(1,\mathbf{0})) (v_1 \otimes v_2) \rangle$ is a canonical
supermeromorphic function in $(z,\theta)$ with poles only at $0$ and
$\infty$, there exists formal Laurent series $f_0(x)$ and $f_1(x)$
such that  
\[\left. (f_0(x) + \varphi f_1(x)) \right|_{(x,\varphi) = (z,\theta)}
= \langle v', \nu_2 ((z,\theta); \mathbf{0}, (1, \mathbf{0}), (1,
\mathbf{0})) (v_1 \otimes v_2) \rangle . \] {}From this equality and 
(\ref{define Y}), we see that 
\[Y_\nu (v_1, (x,\varphi)) \in (\mbox{End} \; V) [[x,x^{-1}]] [\varphi]
.\] 
By the supermeromorphicity axiom, $\langle v', \nu_2 ((z,\theta);
\mathbf{0}, (1, \mathbf{0}), (1, \mathbf{0})) (u \otimes v)
\rangle$ is supermeromorphic in $(z,\theta)$, and there exists $N(u,v)
\in \Z$ such that for any $v' \in V'$ the order of the pole at $0$ is
less than $N(u,v)$.  Thus by (\ref{define Y}), we must have $u_nv = 0$
for $n \geq N(u,v)$, and $n \in \frac{1}{2} \mathbb{Z}$. 

For any $(z,\theta) \in (\bigwedge_\infty^0)^\times \oplus
\bigwedge_\infty^1$, using the sewing and grading axioms,
\begin{eqnarray*}
\left. \langle v', Y_\nu (\mathbf{1}_\nu, (x,\varphi))v \rangle
\right|_{(x,\varphi) = (z,\theta)} 
&=& \sum_{k \in \frac{1}{2} \mathbb{Z}} \left. \langle v', Y_\nu
(P(k) (\mathbf{1}_\nu), (x,\varphi))v \rangle \right|_{(x,\varphi) = 
(z,\theta)} \\
&=& \langle v', (\nu_2 ((z,\theta);\mathbf{0}, (1, \mathbf{0}),
(1, \mathbf{0})) \; _1*_0 \; \nu_0(\mathbf{0})) (v) \rangle \\
&=& \langle v', \nu_1 (((z,\theta);\mathbf{0}, (1, \mathbf{0}),
(1, \mathbf{0})) \; _1\infty_0 \; (\mathbf{0})) (v) \rangle \\
&=& \langle v', \nu_1 (\mathbf{0}, (1, \mathbf{0})) (v) \rangle \\
&=& \langle v',v \rangle .
\end{eqnarray*}
Thus $Y_\nu (\mathbf{1}_\nu, (x,\varphi)) = 1$.  Similarly,
\begin{eqnarray*}
\left. \langle v', Y_\nu (v, (x,\varphi)) \mathbf{1}_\nu \rangle
\right|_{(x,\varphi) = (z,\theta)} 
&=& \sum_{k \in \frac{1}{2} \mathbb{Z}} \left. \langle v', Y_\nu
(v, (x,\varphi)) P(k) (\mathbf{1}_\nu) \rangle \right|_{(x,\varphi) =
(z,\theta)} \\
&=& \langle v', (\nu_2 ((z,\theta);\mathbf{0}, (1, \mathbf{0}),
(1, \mathbf{0})) \; _2*_0 \; \nu_0(\mathbf{0})) (v) \rangle \\
&=& \langle v', \nu_1 (((z,\theta);\mathbf{0}, (1, \mathbf{0}),
(1, \mathbf{0})) \; _2\infty_0 \; (\mathbf{0})) (v) \rangle \\
&=& \langle v', \nu_1 ((A(-z,1), M(-\theta, \frac{1}{2})), (1, \mbox{\bf
0})) (v) \rangle .
\end{eqnarray*}
By the supermeromorphicity axiom, we know that $\langle v', \nu_1 (
\cdot ) (v) \rangle$ is a canonical supermeromorphic function on
$SK(2)$.  Thus $\langle v', \nu_1 ((A(-z,1), M(-\theta, 1/2)),
(1, \mathbf{0})) (v) \rangle$ is a polynomial in $z$ and $\theta$,
and by the grading axiom, 
\begin{eqnarray*}
\lim_{(z,\theta) \rightarrow 0} \langle v', \nu_1 ((A(-z,1), M(-\theta,
\frac{1}{2})), (1, \mathbf{0})) (v) \rangle &=& \langle v', \nu_1, (
\mathbf{0} , (1, \mathbf{0})) (v) \rangle \\
&=& \langle v',v \rangle .
\end{eqnarray*}
Thus $\lim_{(x,\varphi) \rightarrow 0} Y_\nu (v, (x,\varphi)) \mathbf{1}_\nu 
= v$ which proves the creation property.

We next prove the axioms involving the $N=1$ Neveu-Schwarz algebra.  Using
the notation $AM^{(01)} = ((A^{(0)}, M^{(0)}), (a^{(1)}, A^{(1)},
M^{(1)})) \in SK(1)$, we define operators $L(n)$ and $G(n -1/2)$, for 
$n \in \mathbb{Z}$, {}from $V$ to $(V')^*$ by 
\begin{eqnarray}
\biggl. - \frac{\partial}{\partial A_{-n}^{(0)}} \langle v',
\nu_1(AM^{(01)}) (v) \rangle \biggr|_{AM^{(01)} = e} &=& \langle v',
L(n) v \rangle , \quad n < 0 , \\
\biggl. - \frac{1}{2} \frac{\partial}{\partial a^{(1)}} \langle v',
\nu_1(AM^{(01)}) (v) \rangle \biggr|_{AM^{(01)} = e} &=& \langle v', 
L(0) v \rangle ,  \\
\biggl. - \frac{\partial}{\partial A_{n}^{(1)}} \langle v',
\nu_1(AM^{(01)}) (v) \rangle \biggr|_{AM^{(01)} = e} &=& \langle v',
L(n) v \rangle , \quad n > 0 , \\
\hspace{.2in} \biggl. - \frac{\partial}{\partial M_{-n -
\frac{1}{2}}^{(0)}} \langle
 v', \nu_1(AM^{(01)}) (v) \rangle \biggr|_{AM^{(01)} = e} &=&
\langle v', G(n + \frac{1}{2}) v \rangle , \quad n < 0 , \label{defineG}
\\
\biggl. - \frac{\partial}{\partial M_{n - \frac{1}{2}}^{(1)}} \langle
v', \nu_1(AM^{(01)}) (v) \rangle \biggr|_{AM^{(01)} = e} &=& \langle
v', G(n - \frac{1}{2}) v \rangle , \quad n > 0 .
\end{eqnarray} 
{}From the grading axiom, when $v \in V_{(k)}$, for $k \in \frac{1}{2}
\mathbb{Z}$,
\begin{eqnarray*}
\langle v', L(0)v \rangle &=& \biggl. - \frac{1}{2}
\frac{\partial}{\partial a^{(1)}} \langle v', \nu_1(AM^{(01)}) (v)
\rangle \biggr|_{AM^{(01)} = e} \\
&=& \biggl. - \frac{1}{2} \frac{\partial}{\partial a^{(1)}} \langle v',
\nu_1(\mathbf{0}, (a^{(1)},\mathbf{0})) (v) \rangle
\biggr|_{a^{(1)}=1} \\
&=& \biggl. - \frac{1}{2} \frac{\partial}{\partial a^{(1)}}
((a^{(1)})^{-2k} \langle v', v \rangle )\biggr|_{a^{(1)}=1} \\
&=& k \langle v',v \rangle 
\end{eqnarray*}
for all $v' \in V'$.  Thus $L(0)v = kv$ for $v \in V_{(k)}$.

By the sewing axiom 
\begin{multline*}
\langle v', \nu_1 (AM^{(01)} \; _1\infty_0 \; BN^{(01)}) (v)
\rangle \\
=\langle v', ( \nu_1 (AM^{(01)}) \; _1*_0 \;
\nu_1(BN^{(01)})) (v) \rangle  e^{- \Gamma (a^{(1)}, A^{(1)},
M^{(1)}, B^{(0)}, N^{(0)})c} . 
\end{multline*}
Using this, Proposition \ref{NS bracket}, the bracket operation on
$\hat{T}_eSK(1)$ (the subspace of the supermeromorphic tangent space
at the identity of $SK(1)$ defined in Section \ref{tangent section}), the
definitions for $L(n)$ and $G(n - 1/2)$ given above, and Proposition
\ref{central charge sewing}, we obtain
\begin{multline}\label{L relations}
\sum_{k \in \frac{1}{2} \mathbb{Z}} \langle v', L(m) P(k) (L(n)v)
\rangle - \sum_{k \in \frac{1}{2} \mathbb{Z}} \langle v', L(n) P(k)
(L(m)v) \rangle \\
= \langle v', (m-n)L(m + n)v \rangle + \langle v',
\frac{m^3 - m}{12} c \delta_{m+n,0} v \rangle , 
\end{multline}
\begin{multline}\label{GL relations}
\sum_{k \in \frac{1}{2} \mathbb{Z}} \langle v', G(m + \frac{1}{2})
P(k) (L(n)v) \rangle - \sum_{k \in \frac{1}{2} \mathbb{Z}} \langle v',
L(n) P(k) (G(m + \frac{1}{2})v) \rangle \\
 = \langle v', (m- \frac{n - 1}{2})G(m + n + \frac{1}{2})v \rangle , 
\end{multline}
\begin{multline}\label{G relations}
\sum_{k \in \frac{1}{2} \mathbb{Z}} \langle v', G(m + \frac{1}{2})
P(k) (G(n - \frac{1}{2})v) \rangle - \sum_{k \in \frac{1}{2} \mathbb{Z}}
\langle v', G(n - \frac{1}{2}) P(k) (G(m + \frac{1}{2})v) \rangle\\
= \langle v', 2L(m + n)v \rangle + \langle v',
\frac{m^2 + m}{3} c \delta_{m+n,0} v \rangle , 
\end{multline}
for $m,n \in \mathbb{Z}$.  Taking $m = 0$, $n \neq 0$ and $v \in V_{(j)}$ 
in (\ref{L relations}), we obtain
\[ \langle v', L(0)L(n)v \rangle = \langle v', (j-n)L(n)v \rangle
. \]
Thus $L(0)L(n)v = (j-n)L(n)v$ which implies that $L(n)$ is an operator
mapping $V_{(j)}$ to $V_{(j - n)}$.  In particular, $L(n)$ is an
operator {}from $V$ to itself.
Taking $n = 0$, and $v \in V_{(j)}$ in (\ref{GL relations}), we obtain
\[ \langle v', L(0)G(m + \frac{1}{2})v \rangle = \langle v',
(j- (m + \frac{1}{2}))G(m + \frac{1}{2})v \rangle . \]
Thus $L(0)G(m + 1/2)v = (j- (m + 1/2))G(m +1/2)v$ which implies that 
$G(m +1/2)$ is an operator mapping $V_{(j)}$ to $V_{(j - (m + 1/2))}$.  In 
particular, $G(m + 1/2)$ is an operator {}from $V$ to itself.  Thus (\ref{L
relations}), (\ref{GL relations}) and (\ref{G relations}) are
equivalent to the $\mathfrak{ns}$ relations.  

By the definition of $\mathcal{G}_e (z, \theta)$ (see (\ref{define G})),
the definition of $\tau_\nu$, and the sewing axiom, 
\begin{eqnarray*}
\lefteqn{\mathcal{G}_e (z, \theta) ( \langle v', \nu_1 ( \cdot )
(v) \rangle ) } \\
&=& \Bigl. \frac{\partial}{\partial \epsilon} \langle v', \nu_1
(((z,\theta); \mathbf{0}, (1, \mathbf{0}), (1, \mathbf{0})) \;
_1\infty_0 \; (\mathbf{0}, M(\epsilon, \frac{3}{2}))) (v) \rangle
\Bigr|_{\epsilon = 0} \\
&=& \frac{\partial}{\partial \epsilon} \langle
v', (\nu_2 ((z,\theta); \mathbf{0}, (1, \mathbf{0}), (1,
\mathbf{0})) \; _1*_0 \; \nu_0(\mathbf{0}, M(\epsilon, \frac{3}{2})))
(v) \rangle \cdot \\
& & \hspace{2.5in} \Bigl. e^{-\Gamma(1,\mathbf{0},\mathbf{0}, \mathbf{0}, 
M(\epsilon, \frac{3}{2}) )c} \Bigr|_{\epsilon = 0} \\
&=& \Bigl. (-1)^{\eta(v')} \langle v', (\nu_2
((z,\theta); \mathbf{0}, (1, \mathbf{0}), (1, \mathbf{0})) \;
_1*_0 \; \frac{\partial}{\partial \epsilon} \nu_0(\mathbf{0},
M(\epsilon, \frac{3}{2}))\Bigr|_{\epsilon = 0}) (v) \rangle \\
&=& - \biggl. (-1)^{\eta(v')} \sum_{k \in \frac{1}{2} \mathbb{Z}} \langle
v', Y_\nu (P(k) (\tau_\nu),(x,\varphi)) v \rangle \biggr|_{(x,\varphi)
= (z,\theta)} \\
&=& - \biggl. (-1)^{\eta(v')} \sum_{k \in \frac{1}{2} \mathbb{Z}} \langle
v', Y_\nu (\tau_\nu,(x,\varphi)) v \rangle \biggr|_{(x,\varphi) =
(z,\theta)} .
\end{eqnarray*}
On the other hand, by Proposition \ref{the linear functional G}
\begin{eqnarray*}
\lefteqn{\mathcal{G}_e (z, \theta) ( \langle v', \nu_1 ( \cdot ) (v) \rangle
) }\\ 
&=& \Biggl( \sum_{k = 0}^1 \sum_{j \in \Z} z^{-(2k - 1)j
- 2 + k} \biggl. \frac{\partial}{\partial M_{j - \frac{1}{2}}^{(k)}}
\biggr|_e \Biggr. \\ 
& & \Biggl.\qquad + \; 2 \theta \biggl( z^{-2}
\biggl. \frac{\partial}{\partial a_0^{(1)}} \biggr|_e + \sum_{k = 0}^1
\sum_{j \in \Z} z^{-(2k - 1)j - 2} \biggl. \frac{\partial}{\partial
A_j^{(k)}} \biggr|_e \biggr) \Biggr) (\langle v', \nu_1 ( \cdot ) (v)
\rangle ) \\ 
&=& - (-1)^{\eta(v')} \sum_{k = 0}^1 \sum_{j \in \Z}
\langle v', z^{-(2k - 1)j - 2 + k} G((2k-1)(j - \frac{1}{2}))v \rangle
\\ 
& & \qquad  + \; 2 \theta \langle v', z^{-2} L(0)v
\rangle + 2 \theta \sum_{k = 0}^1 \sum_{j \in \Z} \langle v', z^{-(2k
- 1)j - 2} L((2k-1)j)v \rangle \\ 
&=& - \biggl. (-1)^{\eta(v')} \langle v', \biggl( \sum_{n
\in \mathbb{Z}} G(n + \frac{1}{2}) x^{-n-2} + 2
\varphi \sum_{n \in \mathbb{Z}} L(n) x^{-n-2} \biggr) v \rangle
\biggr|_{(x,\varphi) = (z,\theta)} .
\end{eqnarray*}
Thus $Y_\nu (\tau_\nu, (x,\varphi)) = \sum_{n \in \mathbb{Z}} G(n +
1/2) x^{-n- 1/2 - 3/2} + 2 \varphi \sum_{n \in
\mathbb{Z}} L(n) x^{-n-2}$.

By Proposition \ref{for derivatives} and the sewing axiom, 
\begin{eqnarray}
\lefteqn{\left. \langle v', Y_\nu (v_1,(x,\varphi))v_2 \rangle
\right|_{(x,\varphi) = (z_0 + z + \theta_0 \theta, \theta + \theta_0)}} \nonumber\\
&=& \! \! \langle v', \nu_2 ((z_0 + z + \theta_0 \theta, \theta +
\theta_0); \mathbf{0} , (1, \mathbf{0}), (1, \mathbf{0})) (v_1 \otimes
v_2) \rangle \nonumber \\
\hspace{.2in} &=& \! \! \langle v', (\nu_2 ((z, \theta); \mathbf{0} , (1,
\mathbf{0}), (1, \mathbf{0})) \; _1*_0 \; \nu_1 ((A(-z_0,1), M(-\theta_0,
\frac{1}{2})), (1, \mathbf{0})))  \label{derivative sewing} \\
& &\hspace{1.7in} (v_1 \otimes v_2) \cdot e^{-\Gamma(1, \mathbf{0},
\mathbf{0},A(-z_0,1), M(-\theta_0, \frac{1}{2})) c} \rangle . \nonumber
\end{eqnarray}
By taking $\partial / \partial \theta_0 + \theta_0
\partial / \partial z_0$ of (\ref{derivative sewing}) and
setting $(z_0,\theta_0) = 0$, we obtain 
\begin{eqnarray*}
\lefteqn{\biggl. \biggl( \Bigl. \Bigl( \frac{\partial}{\partial
\theta_0} + \theta_0 \frac{\partial}{\partial z_0} \Bigr)  \langle v',
Y_\nu (v_1, (x, \varphi))v_2 \rangle \Bigr|_{(x,\varphi) = (z_0 + z +
\theta_0 \theta, \theta + \theta_0)} \biggr) \biggr|_{(z_0,\theta_0) = 
0} } \\
&=& (-1)^{\eta(v')} \langle v', (\nu_2 ((z, \theta);
\mathbf{0} , (1, \mathbf{0}), (1, \mathbf{0})) \; _1*_0 \; \\
& & \qquad \biggl. \Bigl( \frac{\partial}{\partial
\theta_0} + \theta_0 \frac{\partial}{\partial z_0} \Bigr) \nu_1
((A(-z_0,1), M(-\theta_0, \frac{1}{2})), (1, \mathbf{0}))
\biggr|_{(z_0, \theta_0) = 0} ) (v_1 \otimes v_2) \rangle  \\  
&=& (-1)^{\eta(v')} \langle  v', (\nu_2 ((z, \theta); \mathbf{0} ,
(1, \mathbf{0}), (1, \mathbf{0})) \; _1*_0 \; \\
& & \hspace{1.6in} \left. \frac{\partial}{\partial \theta_0} \nu_1
((\mathbf{0}, M(-\theta_0, \frac{1}{2})), (1, \mathbf{0}))
\right|_{\theta_0 = 0} ) (v_1 \otimes v_2) \rangle .  
\end{eqnarray*}
Since for any $f(x) + \varphi g(x) \in \bigwedge_\infty [[x,x^{-1}]]
[\varphi]$,  
\begin{eqnarray*}
\lefteqn{ \left. \Bigl( \frac{\partial}{\partial \varphi} +
\varphi \frac{\partial}{\partial x} \Bigr) (f(x) + \varphi g(x))
\right|_{(x,\varphi) = (z,\theta)} }\\ 
&=& \theta f'(z) + g(z) \\
&=& \left. (\theta f'(z_0 + z) + g(z_0 + z)) \right|_{z_0 = 0} \\
&=& \left. \frac{\partial}{\partial \theta_0} (f(z_0 + z) + \theta_0
\theta f'(z_0 + z) + (\theta + \theta_0) g(z_0 + z) \right|_{(z_0,
\theta_0) = 0} \\
&=& \left. \Bigl( \frac{\partial}{\partial \theta_0} +
\theta_0 \frac{\partial}{\partial z} \Bigr) (f(z_0 + z + \theta_0
\theta) + (\theta + \theta_0) g(z_0 + z + \theta_0 \theta))
\right|_{(z_0, \theta_0) = 0},
\end{eqnarray*}
by the definition of $G(- 1/2)$ given by (\ref{defineG}),
\begin{eqnarray*}
\lefteqn{\Bigl. \langle v', \Bigl( \frac{\partial}{\partial
\varphi} + \varphi \frac{\partial}{\partial x} \Bigr) Y_\nu (v_1,(x,
\varphi)) v_2 \rangle \Bigr|_{(x,\varphi) = (z, \theta)}} \\
&=& \! \!  (-1)^{\eta(v')} \Bigl.  \Bigl( \frac{\partial}{\partial
\varphi} + \varphi \frac{\partial}{\partial x} \Bigr) \langle v',
Y_\nu (v_1, (x, \varphi))v_2) \rangle  \Bigr|_{(x,\varphi) =
(z,\theta)} \\  
&=& \! \! \Bigl. (-1)^{\eta(v')} \Bigl( \frac{\partial}{\partial
\theta_0} + \theta_0 \frac{\partial}{\partial z_0} \Bigr) \! \left(
\left. \langle v', Y_\nu (v_1, (x, \varphi))v_2) \rangle
\right|_{(x,\varphi) = (z_0 + z + \theta_0 \theta, \theta + \theta_0)}
\right)  \Bigr|_{(z_0, \theta_0) = 0} \\
&=& \! \! \langle \Bigl. v', (\nu_2 ((z, \theta); \mathbf{0}, (1,
\mathbf{0}), (1, \mathbf{0})) \; _1*_0 \; \frac{\partial}{\partial
\theta_0} \nu_1 ( (\mathbf{0}, M(-\theta_0,  \frac{1}{2})), (1,
\mathbf{0})) \Bigr|_{\theta_0 = 0} ) (v_1 \otimes v_2) \rangle \\ 
&=& \! \! \Bigl. \sum_{k \in \frac{1}{2} \mathbb{Z}} \langle v', Y_\nu
(P(k) (G(-\frac{1}{2})v_1), (x, \varphi)) v_2 \rangle
\Bigr|_{(x,\varphi) = (z, \theta)} \\ 
&=& \! \! \Bigl. \langle v', Y_\nu (G(-\frac{1}{2})v_1, (x, \varphi)) v_2
\rangle \Bigr|_{(x,\varphi) = (z, \theta)} ,
\end{eqnarray*}
or equivalently
\[\Bigl( \frac{\partial}{\partial \varphi} + \varphi
\frac{\partial}{\partial x} \Bigr) Y_\nu (v_1,(x, \varphi)) =  Y_\nu
(G(-\frac{1}{2})v_1, (x, \varphi)) .\]

Finally, we prove the Jacobi identity by proving duality.  Using
the sewing axiom and Proposition \ref{for comm and assoc} and assuming
that $|(z_1)_B| > |(z_2)_B| > 0$, we have 
\begin{eqnarray*}
\lefteqn{\left. \langle v', Y_\nu (v_1,(x_1,\varphi_1))
Y_\nu (v_2,(x_2,\varphi_2)) v \rangle \right|_{(x_i,\varphi_i) =
(z_i,\theta_i)}} \\
&=& \sum_{k \in \frac{1}{2} \mathbb{Z}} \left. \langle v', Y_\nu
(v_1,(x_1,\varphi_1)) P(k) (Y_\nu (v_2,(x_2,\varphi_2))v) \rangle
\right|_{(x_i,\varphi_i) = (z_i,\theta_i)} \\
&=& \langle v', (\nu_2 ((z_1,\theta_1); \mathbf{0}, (1, \mbox{\bf
0}), (1, \mathbf{0})) \; _2*_0 \; \nu_2 ((z_2,\theta_2); \mbox{\bf
0}, (1, \mathbf{0}), (1, \mathbf{0}))) (v_1 \otimes v_2 \otimes v)
\rangle \\
&=& \langle v', \nu_3 (((z_1,\theta_1); \mathbf{0}, (1, \mbox{\bf
0}), (1, \mathbf{0})) \; _2\infty_0 \; ((z_2,\theta_2); \mbox{\bf
0}, (1, \mathbf{0}), (1, \mathbf{0}))) (v_1 \otimes v_2 \otimes v)
\rangle \\
&=& \langle v', \nu_3 ((z_1,\theta_1),(z_2,\theta_2); \mathbf{0},
(1, \mathbf{0}), (1, \mathbf{0}), (1, \mathbf{0})) (v_1 \otimes
v_2 \otimes v) \rangle .
\end{eqnarray*} 
Similarly, when $0 < |(z_1)_B| < |(z_2)_B|$,
\begin{multline*}
\left. \langle v', Y_\nu (v_2,(x_2,\varphi_2))
Y_\nu (v_1,(x_1,\varphi_1)) v \rangle \right|_{(x_i,\varphi_i) =
(z_i,\theta_i)} \\
= \langle v', \nu_3 ((z_2,\theta_2),(z_1,\theta_1); \mathbf{0},
(1, \mathbf{0}), (1, \mathbf{0}), (1, \mathbf{0})) (v_2 \otimes
v_1 \otimes v) \rangle . 
\end{multline*}
By the permutation axiom,
\begin{multline*}
\langle v', \nu_3 ((z_1,\theta_1),(z_2,\theta_2);
\mathbf{0}, (1, \mathbf{0}), (1, \mathbf{0}), (1, \mathbf{0}))
(v_1 \otimes v_2 \otimes v) \rangle \\
= (-1)^{\eta(v_1)\eta(v_2)} \langle v', \nu_3
((z_2,\theta_2),(z_1,\theta_1); \mathbf{0}, (1, \mathbf{0}), (1,
\mathbf{0}), (1, \mathbf{0})) (v_2 \otimes v_1 \otimes v) \rangle
\end{multline*}
for any $(z_1,\theta_1),(z_2,\theta_2) \in \bigwedge_\infty$ satisfying
$(z_1)_B \neq (z_2)_B$, and $(z_1)_B , (z_2)_B \neq 0$.  By the
supermeromorphicity axiom, both the right-hand side and the left-hand
side of the above equality are rational superfunctions of the form
\[ \frac{g((z_1,\theta_1),(z_2,\theta_2))}{z_1^{s_1} z_2^{s_2} (z_1 -
z_2 - \theta_1 \theta_2)^{s_{12}}} \]
where $g((z_1,\theta_1),(z_2,\theta_2)) \in \bigwedge_\infty [z_1,
\theta_1, z_2, \theta_2]$ and $s_1, s_2, s_{12} \in \mathbb{N}$.  Thus we
have proved the rationality and supercommutativity.

To prove associativity, we use the sewing axiom and Proposition
\ref{for comm and assoc} to see that for $|(z_1)_B| > |(z_2)_B| >
|(z_1)_B - (z_2)_B| > 0$
\begin{eqnarray*}
\lefteqn{\Bigl. \langle v', Y_\nu (Y_\nu (v_1, (x_0,
\varphi_1 - \varphi_2)) v_2, (x_2,\varphi_2)) v_3 \rangle
\Bigr|_{ \begin{scriptsize} \begin{array}{rl}
(x_0, \varphi_1 - \varphi_2) \! \! \! &= (z_1 - z_2 - \theta_1 \theta_2
, \theta_1 - \theta_2), \\
(x_2,\varphi_2) \! \! \! &= (z_2, \theta_2) \end{array} \end{scriptsize}} }\\
&=& \sum_{k \in \frac{1}{2} \mathbb{Z}} \langle v', Y_\nu (P(k) 
(Y_\nu (v_1, (x_0, \varphi_1 - \varphi_2)) v_2), \\
& & \hspace{1.6in} \Bigl.  (x_2,\varphi_2)) v_3 
\rangle \Bigr|_{\begin{scriptsize} \begin{array}{rl}
(x_0, \varphi_1 - \varphi_2) \! \! \! &= (z_1 - z_2 -
\theta_1 \theta_2 , \theta_1 - \theta_2), \\
(x_2,\varphi_2) \! \! \! &= (z_2, \theta_2) \end{array} \end{scriptsize} } \\  
&=& \langle v', (\nu_2((z_2,\theta_2); \mathbf{0}, (1, \mbox{\bf
0}), (1, \mathbf{0})) \; _1*_0 \; \nu_2((z_1 - z_2 - \theta_1
\theta_2, \theta_1 - \theta_2); \\
& & \hspace{2.5in} \mathbf{0}, (1, \mathbf{0}), (1,
\mathbf{0}))) (v_1 \otimes v_2 \otimes v_3) \rangle \\
&=& \langle v', \nu_3(((z_2,\theta_2); \mathbf{0}, (1, \mbox{\bf
0}), (1, \mathbf{0})) \; _1\infty_0 \; ((z_1 - z_2 - \theta_1
\theta_2, \theta_1 - \theta_2); \\
& & \hspace{2.5in} \mathbf{0}, (1, \mathbf{0}), (1,\mathbf{0}))) 
(v_1 \otimes v_2 \otimes v_3) \rangle \\
&=& \langle v', \nu_3((z_1, \theta_1),(z_2,\theta_2); \mathbf{0},
(1, \mathbf{0}), (1, \mathbf{0}), (1, \mathbf{0})) (v_1 \otimes
v_2 \otimes v_3) \rangle \\ 
&=& \left. \langle v', Y_\nu (v_1, (x_1, \varphi_1)) Y_\nu (v_2,
(x_2,\varphi_2)) v_3 \rangle \right|_{(x_i, \varphi_i) = (z_i,
\theta_i)} . 
\end{eqnarray*}
{}From the rationality, supercommutativity, associativity and
Proposition \ref{duality}, we obtain the Jacobi identity.

This finishes the proof that $(V, Y_\nu, \mathbf{1}_\nu, \tau_\nu)$
is an $N=1$ NS-VOSA.  
\end{proof}

\section{$N=1$ SG-VOSAs {}from $N=1$ NS-VOSAs}\label{geom {}from
alg}

In this section, given an $N=1$ NS-VOSA, we construct an $N=1$ SG-VOSA,
assuming that Conjecture \ref{Gamma converges} is true. Let
$(V, Y(\cdot,(x,\varphi)),\bf{1}, \tau)$, be an $N=1$ NS-VOSA with
odd formal variables and rank $c \in \mathbb{C}$.  We define a sequence
of maps: 
\begin{eqnarray*}
\nu_n^Y : SK(n) & \rightarrow & SF_V(n) \\
Q & \mapsto & \nu_n^Y (Q)
\end{eqnarray*}
by
\begin{multline*}
\langle v', \nu_n^Y \left( (z_1, \theta_1), ...,(z_{n-1},
\theta_{n-1}); (A^{(0)}, M^{(0)}), (a^{(1)}, A^{(1)}, M^{(1)}), 
... \right. \\
\left. ... (a^{(n)}, A^{(n)}, M^{(n)}) \right) (v_1 
\otimes \cdots \otimes v_n) \rangle 
\end{multline*} 
\begin{multline}\label{define nu} 
= \iota^{-1}_{1 \cdots n-1} \langle e^{- \sum_{j \in \Z} \left(
A^{(0)}_j L'(j) + M^{(0)}_{j - \frac{1}{2}} G'(j - \frac{1}{2})
\right)} v',\\
Y(e^{- \sum_{j \in \Z} \left(A^{(1)}_j L(j) + M^{(1)}_{j -
\frac{1}{2}} G(j - \frac{1}{2}) \right)} \cdot (a^{(1)})^{-2L(0)}
\cdot v_1, (x_1, \varphi_1)) \cdots\\
Y(e^{- \sum_{j \in \Z} \left(A^{(n-1)}_j L(j) + M^{(n-1)}_{j -
\frac{1}{2}} G(j - \frac{1}{2}) \right)} \cdot (a^{(n-1)})^{-2L(0)} 
\cdot v_{n-1}, (x_{n-1}, \varphi_{n-1})) \cdot \\
\left. e^{- \sum_{j \in \Z} \left(A^{(n)}_j L(j) +
M^{(n)}_{j - \frac{1}{2}} G(j - \frac{1}{2}) \right)} \cdot 
(a^{(n)})^{-2L(0)} \cdot v_n \rangle \right|_{(x_i,\varphi_i) =
(z_i,\theta_i)} 
\end{multline}
for $n \in \Z$, and
\[\langle v', \nu_0^Y (A^{(0)}, M^{(0)})\rangle = \langle e^{- \sum_{j
\in \Z} \left( A^{(0)}_j L'(j) + M^{(0)}_{j - \frac{1}{2}} G'(j -
\frac{1}{2}) \right)} v', \mathbf{1} \rangle \]
where $\iota^{-1}_{1 \cdots n-1}$ is defined in Section \ref{iota}.

To prove that $(V,\nu^Y)$ is an $N=1$ SG-VOSA, and in particular, that it
satisfies the sewing axiom, we will need the following proposition.

\begin{prop}\label{doubly convergent}
Let $l \in \Z$, and $m,n \in \mathbb{N}$ such that $l + m - 1 \in \Z$,
and let $Q_1 \in SK(l)$, $Q_2 \in SK(m)$, $Q_3 \in SK(n)$, and $i,j \in
\Z$ such that $1 \leq i \leq l$, and $1 \leq j \leq l + m -1$.  If the 
sewings $(Q_1 \; _i\infty_0 \; Q_2) \; _j\infty_0 \; Q_3$ exist, then 
for any $v' \in V'$, $v_1,...,v_{l + m + n - 2} \in V$, the series 
\begin{equation}\label{doubly-convergent1}
\langle v', ((\nu^Y_l(Q_1) \; _i*_0 \; \nu_m^Y(Q_2))_{t_1^{1/2}}
\; _j*_0 \; \nu^Y_n(Q_3) )_{t_2^{1/2}} (v_1 \otimes \cdots \otimes
v_{l + m + n - 2}) \rangle
\end{equation}
is doubly absolutely convergent at $t_1^{1/2} = t_2^{1/2} = 1$. 

Similarly, for $1 \leq i \leq l$, $1 \leq j \leq m$, if the sewings 
$Q_1 \; _i\infty_0 \; (Q_2 \; _j\infty_0 \; Q_3)$ exist, then for any
$v' \in V'$, $v_1,...,v_{l + m + n - 2} \in V$, the series 
\begin{equation}\label{doubly-convergent2}
\langle v', (\nu^Y_l(Q_1) \; _i*_0 \; (\nu_m^Y(Q_2) \; _j*_0 \;
\nu^Y_n(Q_3))_{t_1^{1/2}} )_{t_2^{1/2}} (v_1 \otimes \cdots \otimes v_{l +
m + n - 2}) \rangle
\end{equation}  
is doubly absolutely convergent at $t_1^{1/2} = t_2^{1/2} = 1$.  
\end{prop}

\begin{proof}
We prove the first statement.  The second statement is proved similarly.  
Let $t_1^{1/2}, t_2^{1/2} \in \mathbb{C}^\times$.  Let 
$Q_2(t_1^{1/2})$ and $Q_3(t_2^{1/2})$ be the elements of 
$SK(m)$ and $SK(n)$, respectively, obtained {}from $Q_2$ and $Q_3$, 
respectively, by multiplying the local coordinate maps at the $i$-th 
puncture of $Q_2$ and the $j$-th puncture of $Q_3$ by $t_1$ and $t_2$, 
respectively, in the even coordinate, and by $t_1^{1/2}$ and 
$t_2^{1/2}$, respectively, in the odd coordinate.  Then there exist 
two neighborhoods $\Delta_1$ and $\Delta_2$ of $1 \in \mathbb{C}$ such 
that for any $(t_1^{1/2}, t_2^{1/2}) \in \Delta_1 \times
\Delta_2$, the sewings $(Q_1 \; _i\infty_0 \; Q_2(t_1^{1/2})) \;
_j \infty_0 \; Q_3(t_2^{1/2})$ exist.   

This gives a two-parameter family of superspheres
\[\{ (Q_1 \; _i\infty_0 \; Q_2(t_1^\frac{1}{2})) \; _j \infty_0 \;
Q_3(t_2^\frac{1}{2}) \; | \; (t_1^\frac{1}{2}, t_2^\frac{1}{2}) \in \Delta_1
\times \Delta_2 \} .\] 
Then as in the proof of Proposition 4.19 in \cite{B memoirs}, we can take 
the fiber bundle structure of these superspheres and define global sections
to obtain a complex analytic family of compact complex manifolds.
Then using the definition of $\nu^Y$ and the contraction $*$, the proof
that (\ref{doubly-convergent1}) is doubly absolutely convergent at
$t_1^{1/2} = t_2^{1/2} = 1$, is completely analogous to the proof of the
convergence of the series $\Psi_j(t^{1/2})$ in Proposition 4.19 in
\cite{B vosas}, and similarly for (\ref{doubly-convergent2}).
\end{proof}

\begin{prop}\label{get a sgvosa}
The pair $(V,\nu^Y)$ is an $N=1$ SG-VOSA.    
\end{prop}

\begin{proof} 
We must check that $(V,\nu^Y)$ satisfies the five axioms in the
definition of an $N=1$  SG-VOSA. 

(1) Positive energy axiom: This is obvious.

(2) Grading axiom:  For any $v' \in V'$, $v \in V_{(k)}$, for $k \in
\frac{1}{2} \mathbb{Z}$, and $a \in (\bigwedge_\infty^0)^\times$,
\[\langle v', \nu_1^Y (\mathbf{0}, (a,\mathbf{0}))(v) \rangle =
\langle v', a^{-2L(0)} v \rangle = a^{-2k} \langle v',v \rangle . \]

(3) Supermeromorphicity axiom:  The fact that for any $n \in \Z$, $v'
\in V'$, $v_1,...,v_n \in V$, and $Q \in SK(n)$, the function $Q
\mapsto \langle v', \nu_n^Y (Q) (v_1 \otimes \cdots \otimes v_n)
\rangle$ is a canonical supermeromorphic function follows {}from the
rationality of $(V,Y, \bf{1}, \tau)$ and the fact that $V$ is a
positive energy module for $\mathfrak{ns}$.  To prove
the rest of the supermeromorphicity axiom, without loss of generality,
we only discuss those $Q$ of the form
$((z_1,\theta_1),...,(z_{n-1},\theta_{n-1}); \mathbf{0}, (1,
\mathbf{0}),...,(1,\mathbf{0}))$.

Take $N(v_i,v_j) \in \frac{1}{2} \Z$ such that $(v_i)_k v_j = 0$ 
for any $k \geq N(v_i,v_j)$. If $i<j$, then by
supercommutativity and associativity, we have
\begin{eqnarray*}
\lefteqn{\langle v', \nu_1^Y(Q) (v_1 \otimes \cdots \otimes v_n) \rangle}\\
&=& \left. \iota^{-1}_{1 \cdots n - 1} \left( \langle v', Y(v_1,
(x_1, \varphi_1)) \cdots Y(v_{n-1}, x_{n-1}, \varphi_{n-1})) v_n
\rangle \right) \right|_{(x_l, \varphi_l) = (z_l , \varphi_l)} \\
&=& (-1)^{\sum_{m = 1}^{j-i-1} \eta(v_{j - m})\eta(v_j)}
\iota^{-1}_{1 \cdots (i-1)ij(i+1) \cdots (j-1)(j+1) \cdots n-1}
\left( \langle v', Y(v_1, (x_1, \varphi_1)) \cdots \right. \\
& &Y(v_{i-1}, (x_{i-1}, \varphi_{i-1})) Y(v_i, (x_i,\varphi_i))
Y(v_j, (x_j, \varphi_j)) Y(v_{i+1}, (x_{i+1}, \varphi_{i+1}))  \cdots
\\ 
& &Y(v_{j-1}, (x_{j-1}, \varphi_{j-1})) Y(v_{j+1},
(x_{j+1}, \varphi_{j+1})) \cdots \\
& & \left. \left. \hspace{2in} Y(v_{n-1}, (x_{n-1}, \varphi_{n-1}))
v_n \rangle \right) \right|_{(x_l, \varphi_l) = (z_l , \theta_l)} \\
&=& (-1)^{\sum_{m = 1}^{j-i-1} \eta(v_{j - m})\eta(v_j)}
\iota^{-1}_{1 \cdots (i-1)j(i+1) \cdots (j-1)(j+1) \cdots (n-1)0}
\left( \langle v', Y(v_1, (x_1, \varphi_1)) \cdots \right. \\
& &Y(v_{i-1}, (x_{i-1}, \varphi_{i-1})) Y(Y(v_i, (x_0, \varphi_i -
\varphi_j))v_j, (x_j, \varphi_j)) Y(v_{i+1}, (x_{i+1}, \varphi_{i+1}))\\
& & \cdots Y(v_{j-1}, (x_{j-1}, \varphi_{j-1})) Y(v_{j+1},
(x_{j+1}, \varphi_{j+1})) \cdots \\
& & \Bigl. \left. \hspace{1.5in}  Y(v_{n-1}, (x_{n-1}, \varphi_{n-1}))
v_n \rangle \right) \Bigr|_{\begin{scriptsize} \begin{array}{c} (x_l,
\varphi_l) = (z_l , \varphi_l) ,
\; \mbox{for} \; l \neq i,\\ 
x_0 = z_i - z_j - \varphi_i \varphi_j . \end{array} \end{scriptsize}}
\end{eqnarray*}   
By the definition of $N(v_i,v_j)$, we see that the order of the pole
$(z_i, \theta_i) = (z_j, \theta_j)$ of $\langle v', \nu_n^Y (Q) (v_1
\otimes \cdots \otimes v_n) \rangle$ is less than $N(v_i,v_j)$.  The
result for $i > j$ is proved similarly.

(4) Permutation axiom:  If $\sigma \in S_{n-1} \subset S_n$ is a
permutation on the first $n-1$ letters, then {}from supercommutativity, we
have
\begin{equation*}
\langle v',  \nu_n^Y (\sigma \cdot Q)(v_1 \otimes \cdots \otimes
v_n) \rangle \hspace{3.1in}
\end{equation*}  
\begin{multline*} 
= \langle v', \nu_n^Y \left( \sigma \cdot ((z_1, \theta_1), ...,(z_{n-1},
\theta_{n-1}); (A^{(0)}, M^{(0)}),(a^{(1)}, A^{(1)}, M^{(1)}), ..., 
\right. \\
\left. (a^{(n)}, A^{(n)}, M^{(n)})) \right) (v_1 \otimes \cdots \otimes
v_n) \rangle 
\end{multline*}
\begin{multline*}
= \langle v', \nu_n^Y ((z_{\sigma^{-1}(1)}, \theta_{\sigma^{-1}(1)}),
...,(z_{\sigma^{-1}(n-1)}, \theta_{\sigma^{-1}(n-1)}); (A^{(0)},
M^{(0)}), \\
(a^{(\sigma^{-1}(1))}, A^{(\sigma^{-1}(1))}, 
M^{(\sigma^{-1}(1))}), ...,(a^{(\sigma^{-1}(n-1))},  
A^{(\sigma^{-1}(n-1))}, M^{(\sigma^{-1}(n-1))}) ,\\ 
(a^{(n)}, A^{(n)}, M^{(n)})) (v_1 \otimes \cdots 
\otimes v_n) \rangle 
\end{multline*}
\begin{multline*}
= \iota^{-1}_{\sigma^{-1}(1) \cdots \sigma^{-1}(n-1)} \langle  e^{-
\sum_{j \in \Z} \left(A^{(0)}_j L'(j) + M^{(0)}_{j - \frac{1}{2}} G'(j
- \frac{1}{2}) \right)} v', \\
Y(e^{- \sum_{j \in \Z} \left(A^{(\sigma^{-1}(1))}_j
L(j) + M^{(\sigma^{-1}(1))}_{j - \frac{1}{2}} G(j - \frac{1}{2})
\right)} \cdot (a^{(\sigma^{-1}(1))})^{-2L(0)} \cdot v_1,\\  
(x_{\sigma^{-1}(1)}, \varphi_{\sigma^{-1}(1)})) \cdots Y(e^{- \sum_{j \in \Z}
\left(A^{(\sigma^{-1}(n-1))}_j L(j) + M^{(\sigma^{-1}(n-1))}_{j - \frac{1}{2}} 
G(j - \frac{1}{2}) \right)}
\cdot \\ 
(a^{(\sigma^{-1}(n-1))})^{-2L(0)} \cdot v_{n-1}, (x_{\sigma^{-1}(n-1)},
\varphi_{\sigma^{-1}(n-1)})) \cdot \\
e^{- \sum_{j \in \Z} \left(A^{(n)}_j L(j) + M^{(n)}_{j -
\frac{1}{2}} G(j - \frac{1}{2}) \right)} \cdot \left. (a^{(n)})^{-2L(0)} 
\cdot v_n \rangle \right|_{(x_i,\varphi_i) = (z_i,\theta_i)} 
\end{multline*} 
\begin{multline*}
= (-1)^{\eta(\sigma)} \iota^{-1}_{1 \cdots n-1} \langle
e^{- \sum_{j \in \Z} \left( A^{(0)}_j L'(j) + M^{(0)}_{j -
\frac{1}{2}} G'(j - \frac{1}{2}) \right)} v', \\
Y(e^{- \sum_{j \in \Z} \left(A^{(1)}_j L(j) + 
M^{(1)}_{j - \frac{1}{2}} G(j - \frac{1}{2}) \right)} \cdot
(a^{(1)})^{-2L(0)} \cdot v_{\sigma(1)}, (x_1, \varphi_1)) \cdots \\
Y(e^{- \sum_{j \in \Z} \left(A^{(n-1)}_j L(j) + M^{(n-1)}_{j
- \frac{1}{2}} G(j - \frac{1}{2}) \right)} \cdot (a^{(n-1)})^{-2L(0)}
\cdot v_{\sigma(n - 1)}, \\
(x_{n-1}, \varphi_{n-1})) e^{- \sum_{j \in \Z} \left(A^{(n)}_j L(j)
+ M^{(n)}_{j - \frac{1}{2}} G(j - \frac{1}{2}) \right)} \cdot \\
\left.  (a^{(n)})^{-2L(0)} \cdot v_n \rangle \right|_{(x_i,\varphi_i) =
(z_i,\theta_i)} 
\end{multline*}
\begin{multline*}
= (-1)^{\eta(\sigma)} \langle v', \nu_n^Y ( (z_1, \theta_1),
...,(z_{n-1}, \theta_{n-1}); (A^{(0)}, M^{(0)}), (a^{(1)}, A^{(1)},
M^{(1)}),...,\\
(a^{(n)}, A^{(n)}, M^{(n)})) (v_{\sigma(1)} \otimes
\cdots  \otimes v_{\sigma(n-1)}\otimes v_n) \rangle 
\end{multline*}
\begin{multline*}
= \langle v', \nu_n^Y ((z_1, \theta_1), ...,(z_{n-1}, 
\theta_{n-1}); (A^{(0)}, M^{(0)}), (a^{(1)}, A^{(1)}, M^{(1)}),
..., \\
(a^{(n)}, A^{(n)}, M^{(n)}) ) ((v_1 \otimes \cdots \otimes
v_n) \cdot \sigma) \rangle 
\end{multline*}
\begin{equation*}  
= \langle v',  \sigma \cdot (\nu_n^Y (Q))
(v_1 \otimes \cdots \otimes v_n) \rangle \hspace{2.7in}
\end{equation*}
which proves the permutation axiom for $\sigma \in S_{n-1} \subset
S_n$.  For the transposition switching $n$ and $n - 1$,  i.e., for
$(n \; n - 1) \in S_n$, using skew-supersymmetry
(\ref{skew-supersymmetry}),
\[Y(u,(x,\varphi))v = (-1)^{\eta(u) \eta(v)} e^{xL(-1) + \varphi
G(-\frac{1}{2})} Y(v,(-x, - \varphi))u ,\]
we have
\begin{eqnarray*}
\lefteqn{ \langle v', (n \; n - 1) \cdot (\nu_n^Y (Q))(v_1 \otimes
\cdots \otimes v_n) \rangle }\\ 
&=& (-1)^{\eta(v_n) \eta(v_{n-1})} \langle v',  \nu_n^Y (Q)(v_1 \otimes
\cdots \otimes v_{n-2} \otimes v_n \otimes v_{n-1} ) \rangle
\hspace{1.1in}
\end{eqnarray*}
\begin{multline*}
= (-1)^{\eta(v_n) \eta(v_{n-1})} \langle v',\nu_n^Y  \left(((z_1,
\theta_1), ...,(z_{n-1}, \theta_{n-1}); (A^{(0)},
M^{(0)}), \right. \\
\left. (a^{(1)}, A^{(1)}, M^{(1)}), ..., (a^{(n)}, A^{(n)}, M^{(n)}))
\right) (v_1 \otimes \cdots \otimes v_{n-2} \otimes v_n \otimes
v_{n-1})
\rangle 
\end{multline*}
\begin{multline*}
= (-1)^{\eta(v_n) \eta(v_{n-1})} \iota^{-1}_{1 \cdots n-1} \langle
e^{- \sum_{j \in \Z} \left( A^{(0)}_j L'(j) + M^{(0)}_{j -
\frac{1}{2}} G'(j - \frac{1}{2}) \right)} v', \\
Y(e^{- \sum_{j \in \Z} \left(A^{(1)}_j L(j) + M^{(1)}_{j -
\frac{1}{2}} G(j - \frac{1}{2}) \right)} \cdot (a^{(1)})^{-2L(0)}
\cdot v_1, (x_1, \varphi_1)) \cdots \\ 
Y(e^{- \sum_{j \in \Z} \left(A^{(n-1)}_j L(j) + M^{(n-1)}_{j
- \frac{1}{2}} G(j - \frac{1}{2}) \right)} \cdot (a^{(n-1)})^{-2L(0)}
\cdot v_n, (x_{n-1}, \varphi_{n-1})) \cdot \\
e^{- \sum_{j \in \Z} \left(A^{(n)}_j L(j)
+ M^{(n)}_{j - \frac{1}{2}} G(j - \frac{1}{2}) \right)} \cdot \left.
(a^{(n)})^{-2L(0)} \cdot v_{n-1} \rangle \right|_{(x_i,\varphi_i) =
(z_i,\theta_i)} 
\end{multline*}
\begin{multline}\label{permute1}
= \iota^{-1}_{1 \cdots n-1} \langle e^{- \sum_{j \in \Z} \left(
A^{(0)}_j L'(j) + M^{(0)}_{j - \frac{1}{2}} G'(j - \frac{1}{2})
\right)} v', \\
Y(e^{- \sum_{j \in \Z} \left(A^{(1)}_j L(j) + M^{(1)}_{j
- \frac{1}{2}} G(j - \frac{1}{2}) \right)} \cdot (a^{(1)})^{-2L(0)}
\cdot v_1, (x_1, \varphi_1)) \cdots \\
e^{x_{n-1}L(-1) + \varphi_{n-1}
G(-\frac{1}{2})} Y(e^{- \sum_{j \in \Z} \left(A^{(n)}_j L(j) + M^{(n)}_{j -
\frac{1}{2}} G(j - \frac{1}{2}) \right)} \cdot (a^{(n)})^{-2L(0)}
\cdot \\ 
v_{n-1}, (-x_{n-1}, -\varphi_{n-1})) \cdot e^{- \sum_{j \in \Z} \left(A^{(n-1)}_j
L(j) + M^{(n-1)}_{j - \frac{1}{2}} G(j - \frac{1}{2}) \right)} \cdot \\
\left. (a^{(n-1)})^{-2L(0)} \cdot v_n \rangle \right|_{(x_i,\varphi_i) =
(z_i,\theta_i)} .
\end{multline}
{}Using equations (\ref{L(-1) and G(-1/2) exp 1}) and 
(\ref{L(-1) and G(-1/2) exp 2}), the right-hand side of
(\ref{permute1}) is equal to  
\begin {multline*}
\iota^{-1}_{1 \cdots n-1} \langle e^{- \sum_{j \in \Z} \left(
A^{(0)}_j L'(j) + M^{(0)}_{j - \frac{1}{2}} G'(j - \frac{1}{2})
\right)} v', e^{x_{n-1}L(-1) + \varphi_{n-1} G(-\frac{1}{2})} \cdot \\
Y(e^{- \sum_{j \in \Z} \left(A^{(1)}_j L(j) + M^{(1)}_{j
- \frac{1}{2}} G(j - \frac{1}{2}) \right)} \cdot (a^{(1)})^{-2L(0)}
\cdot v_1, s_{(x_{n-1}, \varphi_{n-1})}(x_1, \varphi_1)) \cdot\\
\cdots Y(e^{- \sum_{j \in \Z} \left(A^{(n-2)}_j L(j) + M^{(n-2)}_{j -
\frac{1}{2}} G(j - \frac{1}{2}) \right)} \cdot  (a^{(n-2)})^{-2L(0)}
v_{n-2},\\ s_{(x_{n-1}, \varphi_{n-1})}(x_{n-2}, \varphi_{n-2})) \cdot\\ 
Y(e^{- \sum_{j \in \Z} \left(A^{(n)}_j L(j) + M^{(n)}_{j -
\frac{1}{2}} G(j - \frac{1}{2}) \right)} \cdot (a^{(n)})^{-2L(0)}
\cdot v_{n-1}, (-x_{n-1}, -\varphi_{n-1})) \cdot \\ 
 e^{- \sum_{j \in \Z} \left(A^{(n-1)}_j
L(j) + M^{(n-1)}_{j - \frac{1}{2}} G(j - \frac{1}{2}) \right)} \cdot
\left.(a^{(n-1)})^{-2L(0)} \cdot v_n \rangle \right|_{(x_i,\varphi_i) =
(z_i,\theta_i)} 
\end{multline*}
\begin{multline}\label{permute2}
= \iota^{-1}_{1 \cdots n-1} \langle e^{x_{n-1}L'(1) +
\varphi_{n-1} G'(\frac{1}{2})} \cdot e^{- \sum_{j \in \Z} \left(
A^{(0)}_j L'(j) + M^{(0)}_{j - \frac{1}{2}} G'(j - \frac{1}{2})
\right)} v', \\
Y(e^{- \sum_{j \in \Z} \left(A^{(1)}_j L(j) + M^{(1)}_{j
- \frac{1}{2}} G(j - \frac{1}{2}) \right)} \cdot (a^{(1)})^{-2L(0)}
\cdot v_1, s_{(x_{n-1}, \varphi_{n-1})}(x_1, \varphi_1)) \cdot \\ 
\cdots Y(e^{- \sum_{j \in \Z} \left(A^{(n-2)}_j 
L(j) + M^{(n-2)}_{j - \frac{1}{2}} G(j - \frac{1}{2}) \right)} \cdot 
(a^{(n-2)})^{-2L(0)} v_{n-2},   \\
s_{(x_{n-1},\varphi_{n-1})}(x_{n-2},\varphi_{n-2})) \cdot\\
Y(e^{- \sum_{j \in \Z} \left(A^{(n)}_j L(j) + M^{(n)}_{j -
\frac{1}{2}} G(j - \frac{1}{2}) \right)} \cdot (a^{(n)})^{-2L(0)}
\cdot v_{n-1}, (-x_{n-1}, -\varphi_{n-1})) \cdot \\
e^{- \sum_{j \in \Z} \left(A^{(n-1)}_j L(j) +
M^{(n-1)}_{j - \frac{1}{2}} G(j - \frac{1}{2}) \right)} \cdot 
 \left. (a^{(n-1)})^{-2L(0)} \cdot v_n \rangle \right|_{(x_i,\varphi_i) =
(z_i,\theta_i)} .
\end{multline}
{}From (\ref{Sn action}) and the definition of $\nu^Y$, the right-hand
side of (\ref{permute2}) is equal to 
\begin{multline*} 
\iota^{-1}_{1 \cdots n-1} \langle e^{- \sum_{j \in \Z} \left(
\tilde{A}^{(0)}_j L'(j) + \tilde{M}^{(0)}_{j - \frac{1}{2}} G'(j -
\frac{1}{2}) \right)} v', Y(e^{- \sum_{j \in \Z} \left(A^{(1)}_j L(j) +
M^{(1)}_{j - \frac{1}{2}} G(j - \frac{1}{2}) \right)} \cdot \\
(a^{(1)})^{-2L(0)} \cdot v_1, s_{(x_{n-1}, \varphi_{n-1})}(x_1, \varphi_1))
\cdots Y(e^{- \sum_{j \in \Z} \left(A^{(n-2)}_j L(j) + M^{(n-2)}_{j
- \frac{1}{2}} G(j - \frac{1}{2}) \right)} \cdot \\
\cdot (a^{(n-2)})^{-2L(0)} \cdot v_{n-2}, 
s_{(x_{n-1}, \varphi_{n-1})}(x_{n-2},\varphi_{n-2})) \cdot\\
Y(e^{- \sum_{j \in \Z} \left(A^{(n)}_j L(j) + M^{(n)}_{j -
\frac{1}{2}} G(j - \frac{1}{2}) \right)} \cdot (a^{(n)})^{-2L(0)}
\cdot v_{n-1}, (-x_{n-1}, -\varphi_{n-1})) \cdot \\ 
e^{- \sum_{j \in \Z} \left(A^{(n-1)}_j L(j) +
M^{(n-1)}_{j - \frac{1}{2}} G(j - \frac{1}{2}) \right)} \cdot 
\left. (a^{(n-1)})^{-2L(0)} \cdot v_n \rangle \right|_{(x_i,\varphi_i) =
(z_i,\theta_i)} 
\end{multline*}
\begin{multline*}
= \langle v', \nu_n^Y(s_{(z_{n-1}, \theta_{n-1})}(z_1, \theta_1),
s_{(z_{n-1}, \theta_{n-1})} (z_2, \theta_2),..., s_{(z_{n-1}, \theta_{n-1})}
(z_{n-2},\theta_{n-2}), \\
s_{(z_{n-1}, \theta_{n-1})}(0); (\tilde{A}^{(0)},\tilde{M}^{(0)}),  (a^{(1)},
A^{(1)},M^{(1)}),...,(a^{(n-2)}, A^{(n-2)}, M^{(n-2)}), \\ (a^{(n)}, A^{(n)},
M^{(n)}),(a^{(n-1)}, A^{(n-1)}, M^{(n-1)})) (v_1 \otimes \cdots \otimes
v_n) \rangle 
\end{multline*}
\begin{multline*}
= \langle v', \nu_n^Y \left( (n \; \; n - 1) \cdot ((z_1, \theta_1),
...,(z_{n-1}, \theta_{n-1}); (A^{(0)}, M^{(0)}), (a^{(1)}, A^{(1)},
M^{(1)}), \right.\\
\left. ...,(a^{(n)}, A^{(n)}, M^{(n)})) \right)
(v_1 \otimes \cdots \otimes v_n) \rangle  
\end{multline*}
\begin{equation}\label{permute3}
= \langle v', \nu_n^Y ( (n \; \; n - 1) \cdot Q )
(v_1 \otimes \cdots \otimes v_n) \rangle .  \hspace{2in}
\end{equation}
Combining (\ref{permute1}), (\ref{permute2}), and (\ref{permute3}), we
see that the permutation axiom is true for the transposition $(n \; n
- 1) \in S_n$.  Since $S_n$ is generated by $S_{n-1}$ and $(n \; n -
1)$, the permutation axiom holds for arbitrary $\sigma \in S_n$.

(5) Sewing axiom: Suppose $Q_1 \in SK(m)$ and $Q_2 \in SK(n)$ are
given by 
\begin{multline*}
Q_1 = ((z_1, \theta_1),...,(z_{m-1}, \theta_{m-1});
(A^{(0)},M^{(0)}), (a^{(1)}, A^{(1)}, M^{(1)}),\\
...,(a^{(m)}, A^{(m)}, M^{(m)})) , 
\end{multline*}
\begin{multline*}
Q_2 = ((z_1', \theta_1'),...,(z_{n-1}', \theta_{n-1}');
(B^{(0)},N^{(0)}), (b^{(1)}, B^{(1)}, N^{(1)}),\\
...,(b^{(n)}, B^{(n)}, N^{(n)})) 
\end{multline*}
and the $i$-th tube of $Q_1$, for some $1 \leq i \leq m$, can be sewn with
the 0-th tube of $Q_2$.  We show by induction on $m$ and $n$ that the
sewing axiom holds for $Q_1$ and $Q_2$.  The proof is in eight steps:

(a)  The case $m = 1$, $n = 0,1$. 

(b) - (f) Several special cases needed for the induction.  

(g) Induction on $n$.

(h) Induction on $m$.

We start with step (a): We prove the case $m = 1$, $n = 0,1$ by
considering
\[Q_1 = ((A^{(0)},M^{(0)}), (a^{(1)}, A^{(1)}, M^{(1)})) \in
SK(1)\]
\[Q_2 = ((B^{(0)},N^{(0)}), (b^{(1)}, B^{(1)}, N^{(1)})) \in
SK(1) \]
where $Q_1 \; _1\infty_0 \; Q_2$ exists.  Using the bases (\ref{V
basis}) and (\ref{dual basis}), we have 
\[\langle v', (\nu_1^Y (Q_1) \; _1*_0 \; \nu_1^Y (Q_2))_{t^{1/2}} (v) 
\rangle \hspace{3.4in} \]
\begin{multline*}
= \sum_{k \in \frac{1}{2}\mathbb{Z}} \sum_{i^{(k)} = 1}^{\dim
V_{(k)}} \langle e^{- \sum_{j \in \Z} \left(A^{(0)}_j L'(j) +
M^{(0)}_{j - \frac{1}{2}} G'(j - \frac{1}{2}) \right)} v', \\
e^{- \sum_{j \in \Z} \left(A^{(1)}_j L(j) + M^{(1)}_{j
- \frac{1}{2}} G(j - \frac{1}{2}) \right)} \cdot 
(a^{(1)})^{-2L(0)} \cdot e_{i^{(k)}}^{(k)} \rangle \cdot \\
\langle e^{- \sum_{j \in \Z} \left(B^{(0)}_j L'(j) + N^{(0)}_{j - \frac{1}{2}} 
G'(j - \frac{1}{2}) \right)} (e_{i^{(k)}}^{(k)})^*, \\
e^{- \sum_{j \in \Z} \left(B^{(1)}_j L(j) + N^{(1)}_{j - \frac{1}{2}} 
G(j - \frac{1}{2}) \right)} \cdot (b^{(1)})^{-2L(0)} \cdot v \rangle t^k 
\end{multline*}
\begin{multline*}
= \sum_{k \in \frac{1}{2}\mathbb{Z}} \sum_{i^{(k)} = 1}^{\dim V_{(k)}}
\langle e^{- \sum_{j \in \Z} \left( A^{(0)}_j L'(j) + M^{(0)}_{j -
\frac{1}{2}} G'(j - \frac{1}{2}) \right)} v', \\
e^{- \sum_{j \in \Z} \left(A^{(1)}_j L(j) + M^{(1)}_{j - \frac{1}{2}} 
G(j - \frac{1}{2}) \right)} \cdot (a^{(1)})^{-2L(0)} \cdot e_{i^{(k)}}^{(k)} 
\rangle \cdot \\ 
\langle e^{- \sum_{j \in \Z} \left( B^{(0)}_j L'(j) + N^{(0)}_{j -
\frac{1}{2}} G'(j - \frac{1}{2}) \right)} \cdot t^{L'(0)}
(e_{i^{(k)}}^{(k)})^*, \\
e^{- \sum_{j \in \Z} \left(B^{(1)}_j L(j) +
N^{(1)}_{j - \frac{1}{2}} G(j - \frac{1}{2}) \right)} \cdot
(b^{(1)})^{-2L(0)} \cdot v \rangle
\end{multline*}
\begin{multline*}
= \sum_{k \in \frac{1}{2}\mathbb{Z}} \sum_{i^{(k)} = 1}^{\dim V_{(k)}}
\langle e^{- \sum_{j \in \Z} \left( A^{(0)}_j L'(j) + M^{(0)}_{j -
\frac{1}{2}} G'(j - \frac{1}{2}) \right)} v', \\
e^{- \sum_{j \in \Z} \left(A^{(1)}_j L(j) + M^{(1)}_{j - \frac{1}{2}} 
G(j - \frac{1}{2}) \right)} \cdot (a^{(1)})^{-2L(0)} \cdot e_{i^{(k)}}^{(k)} 
\rangle \cdot \\ 
\langle (e_{i^{(k)}}^{(k)})^*,t^{L(0)} \cdot e^{- \sum_{j \in \Z}
\left( B^{(0)}_j L(-j) + N^{(0)}_{j - \frac{1}{2}} G(- j +
\frac{1}{2}) \right)} \cdot \\
e^{- \sum_{j \in \Z} \left(B^{(1)}_j L(j) + N^{(1)}_{j - \frac{1}{2}} 
G(j - \frac{1}{2}) \right)} \cdot (b^{(1)})^{-2L(0)} \cdot v \rangle
\end{multline*} 
\begin{multline}\label{sew1}
 = \langle e^{- \sum_{j \in \Z} \left(A^{(0)}_j L'(j) + M^{(0)}_{j -
\frac{1}{2}} G'(j - \frac{1}{2}) \right)} v', 
e^{- \sum_{j \in \Z} \left(A^{(1)}_j L(j) + M^{(1)}_{j - \frac{1}{2}} 
G(j - \frac{1}{2}) \right)} \cdot \\
(t^\frac{1}{2} a^{(1)})^{-2L(0)} \cdot 
e^{- \sum_{j \in \Z} \left( B^{(0)}_j L(-j) + N^{(0)}_{j - \frac{1}{2}} 
G(- j + \frac{1}{2}) \right)} \cdot \\
e^{- \sum_{j \in \Z} \left(B^{(1)}_j L(j) 
+ N^{(1)}_{j - \frac{1}{2}} G(j - \frac{1}{2}) \right)} \cdot (b^{(1)})^{-2L(0)} 
\cdot v \rangle .
\end{multline}

Let $\Psi_j(t^{1/2}) = \Psi_j(t^{-1/2} a^{(1)}, A^{(1)},
M^{(1)}, B^{(0)}, N^{(0)})$, for $j \in \frac{1}{2} \mathbb{Z}$, and let 
$\Gamma(t^{1/2}) = \Gamma(t^{-1/2} a^{(1)}, A^{(1)}, M^{(1)},
B^{(0)},N^{(0)})$. {}From Corollary \ref{normal order in End}, equation
(\ref{sew1}) is equal to   
\begin{multline*}
\langle e^{- \sum_{j \in \Z} \left(A^{(0)}_j L'(j) + M^{(0)}_{j -
\frac{1}{2}} G'(j - \frac{1}{2}) \right)} v', e^{\sum_{j \in
\Z} \left(\Psi_{-j}(t^\frac{1}{2}) L(-j) + \Psi_{- j +
\frac{1}{2}}(t^\frac{1}{2}) G(- j + \frac{1}{2}) \right)}\\
e^{\sum_{j \in \Z} \left(\Psi_j(t^\frac{1}{2}) L(j) + \Psi_{j -
\frac{1}{2}}(t^\frac{1}{2}) G(j - \frac{1}{2}) \right)} \;
e^{2\Psi_0(t^\frac{1}{2}) L(0)} 
(t^{-\frac{1}{2}} a^{(1)})^{-2L(0)} \; e^{\Gamma(t^\frac{1}{2})c}
\cdot \\ 
e^{- \sum_{j \in \Z} \left(B^{(1)}_j L(j) + N^{(1)}_{j
- \frac{1}{2}} G(j - \frac{1}{2}) \right)} \cdot (b^{(1)})^{-2L(0)}
\cdot v \rangle 
\end{multline*}
\begin{multline*}  
= \langle e^{\sum_{j \in
\Z} \left(\Psi_{-j}(t^\frac{1}{2}) L'(j) + \Psi_{- j +
\frac{1}{2}}(t^\frac{1}{2}) G'(j - \frac{1}{2}) \right)} \cdot e^{-
\sum_{j \in \Z} \left(A^{(0)}_j L'(j) + M^{(0)}_{j - \frac{1}{2}} G'(j
- \frac{1}{2}) \right)} v', \\
e^{\sum_{j \in \Z} \left(\Psi_j(t^\frac{1}{2}) L(j) + \Psi_{j -
\frac{1}{2}}(t^\frac{1}{2}) G(j - \frac{1}{2}) \right)} \;
e^{2\Psi_0(t^\frac{1}{2}) L(0)} (t^{-\frac{1}{2}}a^{(1)})^{-2L(0)}
\cdot \\ 
e^{- \sum_{j \in \Z} \left(B^{(1)}_j L(j) + N^{(1)}_{j
- \frac{1}{2}} G(j - \frac{1}{2}) \right)} \cdot (b^{(1)})^{-2L(0)}
\cdot v \rangle e^{\Gamma(t^\frac{1}{2})c} .
\end{multline*}
Thus 
\[\langle v', (\nu_1^Y (Q_1) \; _1*_0 \; \nu_1^Y (Q_2))_{t^{1/2}} (v)
\rangle e^{-\Gamma(t^\frac{1}{2})c}  \hspace{2.9in}\]
\begin{multline}\label{sew2}
= \langle e^{\sum_{j \in
\Z} \left(\Psi_{-j}(t^\frac{1}{2}) L'(j) + \Psi_{- j +
\frac{1}{2}}(t^\frac{1}{2}) G'(j - \frac{1}{2}) \right)} \cdot e^{-
\sum_{j \in \Z} \left(A^{(0)}_j L'(j) + M^{(0)}_{j - \frac{1}{2}} G'(j
- \frac{1}{2}) \right)} v', \\ 
e^{\sum_{j \in \Z} \left(\Psi_j(t^\frac{1}{2}) L(j) + \Psi_{j
-\frac{1}{2}}(t^\frac{1}{2}) G(j - \frac{1}{2}) \right)} \;
e^{2\Psi_0(t^\frac{1}{2}) L(0)} (t^{-\frac{1}{2}}a^{(1)})^{-2L(0)}
\cdot \\ 
e^{- \sum_{j \in \Z} \left(B^{(1)}_j L(j) + N^{(1)}_{j
- \frac{1}{2}} G(j - \frac{1}{2}) \right)} \cdot (b^{(1)})^{-2L(0)}
\cdot v \rangle .
\end{multline}
Since $V$ is a positive energy representation of $\mathfrak{ns}$, the
right-hand side of (\ref{sew2}) is a polynomial in
$\Psi_j(t^{1/2})$, for $j \in \frac{1}{2} \mathbb{Z}$, $j \neq 0$, and 
$e^{\pm2\Psi_0(t^{1/2})}$ with coefficients in
$\bigwedge_\infty[t^{1/2}, t^{-1/2}]$.  By Theorem \ref{actual sewing}, the
series $\Psi_j(t^{1/2})$ and $e^{\pm 2 \tilde{\Psi}_0(t^{1/2})}$ converge
to $\tilde{\Psi}_j(a^{(1)},A^{(1)},M^{(1)},B^{(0)},N^{(0)})$ and
$e^{\pm 2 \tilde{\Psi}_0(a^{(1)}, A^{(1)},M^{(1)},B^{(0)},N^{(0)})}$,
respectively, at $t^{1/2} = 1$.  Let $\Psi_j = \Psi_j(a^{(1)}, A^{(1)},
M^{(1)}, B^{(0)},N^{(0)})$ for $j \in \frac{1}{2} \mathbb{Z}$.
{}From equations (\ref{F1}) and (\ref{F2}), Proposition 3.21 of
\cite{B memoirs}, and Theorem \ref{actual sewing}, the
local coordinate at $\infty$ of the canonical supersphere with $1+1$
tubes representing $Q_1 \; _1\infty_0 \; Q_2$, in terms of the
superderivation representation of $\mathfrak{ns}$,
(\ref{L notation}) and (\ref{G notation}), is given by
\begin{multline*} 
\exp \Bigl( -\sum_{j \in \Z} \left(\Psi_{-j} L_{-j}(w,\rho) + 
\Psi_{-j + \frac{1}{2}} G_{-j +\frac{1}{2}}(w,\rho) \right)\Bigr) \cdot  \\
\exp \Bigl(\sum_{j \in \Z} \left(A^{(0)}_j L_{-j}(w,\rho) +
M^{(0)}_{j - \frac{1}{2}} G_{-j + \frac{1}{2}}(w,\rho) \right) \Bigr)\cdot 
\Bigl( \frac{1}{w}, \frac{i\rho}{w} \Bigr), 
\end{multline*}
and the local coordinate at zero is given by
\begin{multline*}
\exp \Bigl(- \sum_{j \in \Z} \left(\Psi_j L_j(w,\rho) + \Psi_{j -
\frac{1}{2}} G_{j - \frac{1}{2}}(w,\rho) \right) \Bigr) \cdot 
(a^{(1)})^{-L_0(w,\rho)} \cdot e^{\Psi_0 L_0(w,\rho)} \cdot \\
\exp \Bigl( -\sum_{j \in \Z} \left(
B^{(1)}_j L_j(w,\rho) + N^{(1)}_{j - \frac{1}{2}} G_{j - \frac{1}{2}}(w,\rho) 
\right)\Bigr) \cdot
(b^{(1)})^{-L_0(w,\rho)} \cdot(w,\rho) .
\end{multline*}
By Proposition 3.17 of \cite{B memoirs} and by Proposition \ref{above},
there  exist series  $(C^{(0)}, O^{(0)})$, $(C^{(1)},O^{(1)}) \in
\bigwedge_\infty^\infty$ and $c^{(1)} \in (\bigwedge_\infty^0)^\times$
such that these local coordinates can be written as 
\[\exp \Bigl(\sum_{j \in \Z} \left(C^{(0)}_j L_{-j}(w,\rho) +
O^{(0)}_{j - \frac{1}{2}} G_{-j + \frac{1}{2}}(w,\rho) \right) \Bigr) \cdot 
\Bigl( \frac{1}{w},\frac{i\rho}{w} \Bigr) \]
and
\[\exp \Bigl( -\sum_{j \in \Z} \left(C^{(1)}_j L_j(w,\rho) + O^{(1)}_{j
- \frac{1}{2}} G_{j - \frac{1}{2}}(w,\rho) \right)\Bigr) \cdot (c^{(1)})^{-L_0(w,\rho)}
\cdot (w,\rho) ,\]
respectively.  Since $V$ is a positive-energy representation of
$\mathfrak{ns}$ and by Proposition \ref{Gamma converges},  when $t^{1/2} =
1$ the series $\Gamma(t^{-1/2} a^{(1)}, A^{(1)}, M^{(1)}, B^{(0)},
N^{(0)})$ converges to $\Gamma(a^{(1)}, A^{(1)}, M^{(1)},
B^{(0)},N^{(0)})$, and we have
\[\langle v', (\nu_1^Y(Q_1) \; _1*_0 \; \nu_2^Y(Q_1)) (v) \rangle
e^{-\Gamma(a^{(1)}, A^{(1)},M^{(1)},B^{(0)},N^{(0)})c}
\hspace{1.9in} \]
\begin{eqnarray*}
&=& \langle e^{- \sum_{j \in \Z} \left(C^{(0)}_j L'(j) + O^{(0)}_{j -
\frac{1}{2}} G'(j - \frac{1}{2}) \right)} v', \\
& & \hspace{1.5in} e^{- \sum_{j \in \Z}
\left(C^{(1)}_j L(j) + O^{(1)}_{j - \frac{1}{2}} G(j - \frac{1}{2})
\right)} \cdot (c^{(1)})^{-2L(0)} v \rangle \\
&=& \langle v', \nu_1^Y(Q_1 \; _1\infty_0 \; Q_2) (v)
\rangle 
\end{eqnarray*} 
which proves the sewing axiom for $n = 1$, $m = 1$.  The case $m = 1$,
$n = 0$ is then given by letting $(b^{(1)},B^{(1)},N^{(1)}) = (1,
\mathbf{0})$ and $v = \mathbf{1}$ above.

(b) We prove the sewing axiom for the special case $Q_1 \; _2\infty_0
\; Q_2$ when 
\[Q_1 = ((z,\theta); \mathbf{0}, (a^{(1)}, A^{(1)}, M^{(1)}), (1,
\mathbf{0})) \in SK(2)\]  
and 
\[Q_2 = ((B^{(0)}, N^{(0)}),(b^{(1)}, B^{(1)}, N^{(1)})) \in SK(1)
. \]  

In this case,
\[\langle v', (\nu_2^Y(Q_1) \; _2*_0 \; \nu_1^Y(Q_2))_{t^{1/2}} (v_1 \otimes
v_2) \rangle \hspace{2.5in} \]
\begin{multline*}
= \sum_{k \in \frac{1}{2}\mathbb{Z}} \sum_{i^{(k)} =
1}^{\dim V_{(k)}} \langle v', Y(e^{- \sum_{j \in \Z} \left(A^{(1)}_j
L(j) +  M^{(1)}_{j - \frac{1}{2}} G(j - \frac{1}{2}) \right)} \cdot
(a^{(1)})^{-2L(0)} \cdot v_1, \\
(x, \varphi)) e_{i^{(k)}}^{(k)} \rangle \cdot \langle e^{- \sum_{j \in \Z} 
\left(B^{(0)}_j L'(j) + N^{(0)}_{j - \frac{1}{2}} G'(j - \frac{1}{2}) \right)} 
\cdot (e^{(k)}_{i^{(k)}})^*, \\
e^{- \sum_{j \in \Z} \left(B^{(1)}_j L(j) +
N^{(1)}_{j - \frac{1}{2}} G(j - \frac{1}{2}) \right)} \cdot 
\left. (b^{(1)})^{-2L(0)} \cdot v_2 \rangle t^k
\right|_{(x,\varphi) = (z,\theta)} 
\end{multline*}
\begin{multline*}
= \sum_{k \in \frac{1}{2}\mathbb{Z}} \sum_{i^{(k)} =
1}^{\dim V_{(k)}} \langle v', Y(e^{- \sum_{j \in \Z} \left(A^{(1)}_j
L(j) +  M^{(1)}_{j - \frac{1}{2}} G(j - \frac{1}{2}) \right)} \cdot
(a^{(1)})^{-2L(0)} \cdot v_1, \\
(x, \varphi)) e_{i^{(k)}}^{(k)} \rangle \cdot \langle e^{- \sum_{j \in \Z} 
\left(B^{(0)}_j L'(j) + N^{(0)}_{j - \frac{1}{2}} G'(j - \frac{1}{2}) \right)} 
\cdot t^{L'(0)} (e^{(k)}_{i^{(k)}})^*, \\
e^{- \sum_{j \in \Z} \left(B^{(1)}_j L(j) + N^{(1)}_{j - \frac{1}{2}} 
G(j - \frac{1}{2}) \right)} \cdot \left. (b^{(1)})^{-2L(0)} \cdot v_2 \rangle 
\right|_{(x,\varphi) = (z,\theta)} 
\end{multline*} 
\begin{multline*}
= \sum_{k \in \frac{1}{2}\mathbb{Z}} \sum_{i^{(k)} =
1}^{\dim V_{(k)}} \langle v', Y(e^{- \sum_{j \in \Z} \left(A^{(1)}_j
L(j) +  M^{(1)}_{j - \frac{1}{2}} G(j - \frac{1}{2}) \right)} \cdot
(a^{(1)})^{-2L(0)} \cdot v_1, \\
(x, \varphi)) e_{i^{(k)}}^{(k)} \rangle \cdot \langle t^{L'(0)} \cdot 
e^{- \sum_{j \in \Z} \left(t^j B^{(0)}_j L'(j) + t^{j - \frac{1}{2}} 
N^{(0)}_{j - \frac{1}{2}} G'(j - \frac{1}{2}) \right)} \cdot (e^{(k)}_{i^{(k)}})^*, \\
e^{- \sum_{j \in \Z} \left(B^{(1)}_j L(j) + N^{(1)}_{j - \frac{1}{2}} G(j -
\frac{1}{2}) \right)} \cdot \left. (b^{(1)})^{-2L(0)} \cdot v_2 \rangle 
\right|_{(x,\varphi) = (z,\theta)} 
\end{multline*} 
\begin{multline*}
= \sum_{k \in \frac{1}{2}\mathbb{Z}} \sum_{i^{(k)} = 1}^{\dim V_{(k)}}
\langle v', Y(e^{- \sum_{j \in \Z} \left(A^{(1)}_j L(j) +  M^{(1)}_{j
- \frac{1}{2}} G(j - \frac{1}{2}) \right)} \cdot (a^{(1)})^{-2L(0)}
\cdot v_1, \\
(x, \varphi)) e_{i^{(k)}}^{(k)} \rangle \cdot \langle (e^{(k)}_{i^{(k)}})^*, 
e^{- \sum_{j \in \Z} \left(t^j B^{(0)}_j L(-j) + t^{j - \frac{1}{2}}N^{(0)}_{j -
\frac{1}{2}} G(- j + \frac{1}{2}) \right)} \cdot t^{L(0)} \cdot \\  
e^{- \sum_{j \in \Z} \left(B^{(1)}_j L(j) +
N^{(1)}_{j - \frac{1}{2}} G(j - \frac{1}{2}) \right)} \cdot \left.
(b^{(1)})^{-2L(0)} \cdot v_2 \rangle \right|_{(x,\varphi) =
(z,\theta)} 
\end{multline*}
\begin{multline}\label{for case d}
= \langle v', Y(e^{- \sum_{j \in \Z} \left(A^{(1)}_j L(j) +
M^{(1)}_{j - \frac{1}{2}} G(j - \frac{1}{2}) \right)} \cdot
(a^{(1)})^{-2L(0)} \cdot v_1, (x, \varphi)) \\
e^{- \sum_{j \in \Z} \left(t^j B^{(0)}_j L(-j) + t^{j - \frac{1}{2}}
N^{(0)}_{j - \frac{1}{2}} G(- j + \frac{1}{2}) \right)} \cdot t^{L(0)}
\cdot \\
e^{- \sum_{j \in \Z} \left(B^{(1)}_j L(j) + N^{(1)}_{j -
\frac{1}{2}} G(j - \frac{1}{2}) \right)} \cdot  
\left. (b^{(1)})^{-2L(0)} \cdot v_2 \rangle
\right|_{(x,\varphi) = (z,\theta)}  
\end{multline}
\begin{multline}\label{case b}
= \langle e^{- \sum_{j \in \Z} \left(t^j B^{(0)}_j L'(j) + t^{j -
\frac{1}{2}} N^{(0)}_{j - \frac{1}{2}} G'(j - \frac{1}{2}) \right)}
\cdot v', \\
e^{\sum_{j \in \Z} \left(t^j B^{(0)}_j L(-j) + t^{j - \frac{1}{2}}
N^{(0)}_{j - \frac{1}{2}} G(- j + \frac{1}{2}) \right)} \cdot \\
Y(e^{- \sum_{j \in \Z} \left(A^{(1)}_j L(j) + M^{(1)}_{j -
\frac{1}{2}} G(j - \frac{1}{2}) \right)} \cdot (a^{(1)})^{-2L(0)}
\cdot v_1, (x, \varphi)) \\
e^{- \sum_{j \in \Z} \left(t^j B^{(0)}_j L(-j) + t^{j - \frac{1}{2}}
N^{(0)}_{j - \frac{1}{2}} G(- j + \frac{1}{2}) \right)} \cdot t^{L(0)}
\cdot \\
e^{- \sum_{j \in \Z} \left(B^{(1)}_j L(j) + N^{(1)}_{j -
\frac{1}{2}} G(j - \frac{1}{2}) \right)} \cdot  
\left. (b^{(1)})^{-2L(0)} \cdot v_2 \rangle
\right|_{(x,\varphi) = (z,\theta)} . 
\end{multline}

The right-hand side of (\ref{case b}) is a Laurent polynomial in
$t^{1/2}$.  If we can prove
\begin{multline*}
e^{\sum_{j \in \Z} \left(t^j B^{(0)}_j L(-j) + t^{j - \frac{1}{2}}
N^{(0)}_{j - \frac{1}{2}} G(- j + \frac{1}{2}) \right)} Y(u,
(x,\varphi)) \\
e^{- \sum_{j \in \Z} \left(t^j B^{(0)}_j L(-j) + t^{j
- \frac{1}{2}} N^{(0)}_{j - \frac{1}{2}} G(- j + \frac{1}{2}) \right)}
\end{multline*}
\begin{equation}\label{new sewing coordinates}
= \; Y\Bigl(e^{-\sum_{j \in \Z} \left(\Theta^{(2)}_j L(j) + \Theta^{(2)}_{j -
\frac{1}{2}} G(j - \frac{1}{2}) \right)} e^{- 2\Theta^{(2)}_0 L(0)} u,
H_t^{-1} \circ I (x,\varphi) \Bigr) 
\end{equation}
where
\begin{equation}\label{Theta in case b}
\Theta^{(2)}_j = \Theta^{(2)}_j(\{t^k B^{(0)}_k, t^{k - \frac{1}{2}}
N^{(0)}_{k - \frac{1}{2}} \}_{k \in \Z},(x,\varphi)) 
\end{equation}
is given by (\ref{define second Theta}) and
\begin{multline}\label{Ht in case b}
H_{t^{1/2}}^{-1} (x,\varphi) = \exp \Bigl(\sum_{j \in \Z} \left(
t^j B^{(0)}_j L_{-j}(x,\varphi) +  t^{j - \frac{1}{2}} N^{(0)}_{j 
- \frac{1}{2}} G_{-j + \frac{1}{2}}(x,\varphi) \right) \Bigr) \cdot
\Bigl(\frac{1}{x},\frac{i\varphi}{x}\Bigr) 
\end{multline}   
then when $t^{1/2} = 1$, the right-hand side of (\ref{case b})
is equal to  
\[\langle v', \nu^Y_2(Q_1 \; _2\infty_0 \; Q_2) (v_1 \otimes v_2) \rangle.
\] 
Thus in this case, the proof of the sewing axiom is reduced to the
proof of (\ref{new sewing coordinates}).  Letting 
\[H_{t^{1/2}}^{-1} \circ I (x,\varphi) = (\tilde{x} (t^\frac{1}{2}),
\tilde{\varphi}(t^\frac{1}{2})) ,\] 
by (\ref{exponential L(-1) and G(-1/2) property}) and Corollary 
\ref{second Theta identity in End}, 
\[Y \Bigl(e^{-\sum_{j \in \Z} \left(\Theta^{(2)}_j L(j) + \Theta^{(2)}_{j -
\frac{1}{2}} G(j - \frac{1}{2}) \right)} e^{- 2\Theta^{(2)}_0 L(0)} u,
H_{t^{1/2}}^{-1} \circ I (x,\varphi) \Bigr) \hspace{1in}\]
\begin{multline*}
= \; Y \Bigl( e^{(\tilde{x}(t^\frac{1}{2}) - x) L(-1) +
(\tilde{\varphi} (t^\frac{1}{2}) - \varphi) G(-\frac{1}{2})} \cdot
e^{- \sum_{j \in \Z} \left( \Theta^{(2)}_j L(j) + \Theta^{(2)}_{j - \frac{1}{2}}
G(j - \frac{1}{2}) \right)} \cdot \\
e^{- 2\Theta^{(2)}_0 L(0)} u, (x,\varphi) \Bigr) 
\end{multline*}
\begin{multline*}
= \; Y \Biggl(\exp \biggl( \sum_{m = -1}^{\infty} \sum_{j \in \Z} \binom{
-j + 1}{ m + 1 } x^{-j - m} \biggl( \Bigl( t^j B_j^{(0)} + 2\varphi t^{j -
\frac{1}{2}} N_{j - \frac{1}{2}}^{(0)} \Bigr)L(m) \\  
+ \; \Bigl( t^{j - \frac{1}{2}}  N_{j - \frac{1}{2}}^{(0)} + \varphi 
x^{-1} \frac{(-j-m)}{2} t^j B_j^{(0)} \Bigr) G(m + \frac{1}{2}) \biggr) 
 \biggr)  u, (x,\varphi) \Biggr) .
\end{multline*} 
By this equality, the proof of (\ref{new sewing coordinates}) is
reduced to the proof of the following bracket formula
\begin{equation}\label{special case b bracket}
\Bigl[ \sum_{j \in \Z} \Bigl( t^j B_j^{(0)} L(-j) + t^{j -
\frac{1}{2}} N_{j -\frac{1}{2}}^{(0)} G( - j + \frac{1}{2}) \Bigr) ,
Y(u,(x,\varphi)) \Bigr] \hspace{.6in}
\end{equation}
\begin{multline*}
= Y \biggl( \sum_{m = -1}^{\infty} \sum_{j \in \Z} \binom{-j + 1}{m + 1}
x^{-j - m} \biggl( \Bigl( t^j B_j^{(0)} + 2\varphi t^{j - \frac{1}{2}} 
N_{j - \frac{1}{2}}^{(0)} \Bigr) L(m)  \\
+ \; \Bigl( t^{j - \frac{1}{2}}  N_{j - \frac{1}{2}}^{(0)} + \varphi
x^{-1} \frac{(-j-m)}{2} t^j B_j^{(0)} \Bigr) G(m + \frac{1}{2})
\biggr) u, (x,\varphi) \biggr) .
\end{multline*} 
{}From (\ref{bracket relation for a vosa}), we have
\begin{multline*} 
\left[ Y(\tau, (x_1,\varphi_1)), Y(u,(x,\varphi)) \right] \\
= \mbox{Res}_{x_0} x^{-1} \delta \left( \frac{x_1 - x_0 -
\varphi_1 \varphi}{x} \right) Y(Y(\tau,(x_0, \varphi_1 -
\varphi))u,(x, \varphi)) .
\end{multline*} 
Let 
\[l_t(x_1) = \sum_{j \in \Z} t^j B_j^{(0)} x_1^{-j + 1} \quad
\mbox{and} \quad g_t(x_1) = \sum_{j \in \Z} t^{j - \frac{1}{2}} N_{j -
\frac{1}{2}}^{(0)} x_1^{-j + 1} .\]
Then 
\begin{eqnarray*}
\varphi_1 \mbox{Res}_{x_1} l_t(x_1) Y(\omega, (x_1,\varphi_1)) &=&
\varphi_1 \sum_{j \in \Z} t^j B_j^{(0)} L(-j) \\ 
\varphi_1 \mbox{Res}_{x_1} g_t(x_1) Y(\tau, (x_1,\varphi_1)) &=&
\varphi_1 \sum_{j \in \Z} t^{j - \frac{1}{2}} N_{j -
\frac{1}{2}}^{(0)} G(- j + \frac{1}{2}) ,
\end{eqnarray*}
and by the $G(-1/2)$-derivative property
(\ref{G(-1/2)-derivative}), we have 
\begin{eqnarray*}
Y(\omega, (x_1,\varphi_1)) &=& Y( \frac{1}{2} G(-\frac{1}{2}) \tau,
(x_1, \varphi_1)) \\
&=& \frac{1}{2} \Bigl( \frac{\partial}{\partial \varphi_1} + \varphi_1
\frac{\partial}{\partial x_1} \Bigr) Y(\tau, (x_1,\varphi_1))  .
\end{eqnarray*}
Thus 
\begin{eqnarray*}
& & \hspace{-.5in} \varphi_1 \biggl[ \sum_{j \in \Z} \left( t^j B_j^{(0)}
L(-j) + t^{j - \frac{1}{2}} N_{j -\frac{1}{2}}^{(0)} G( - j +
\frac{1}{2}) \right) , Y(u,(x,\varphi)) \biggr] \\
&=& \varphi_1 \mbox{Res}_{x_1} \left[ \Bigl( \frac{1}{2}l_t(x_1)
\Bigl( \frac{\partial}{\partial \varphi_1} + \varphi_1
\frac{\partial}{\partial x_1} \Bigr) + g_t(x_1) \Bigr) Y(\tau,
(x_1,\varphi_1)), Y(u,(x,\varphi)) \right] \\
&=& \varphi_1 \mbox{Res}_{x_1} \left( \Bigl( \frac{1}{2}l_t(x_1) \Bigl(
\frac{\partial}{\partial \varphi_1} + \varphi_1
\frac{\partial}{\partial x_1} \Bigr) + g_t(x_1) \Bigr) \right. \\
& & \hspace{.7in} \left. \mbox{Res}_{x_0} x^{-1} \delta \Bigl(
\frac{x_1 - x_0 - \varphi_1 \varphi}{x} \Bigr) Y(Y(\tau,(x_0,
\varphi_1 - \varphi))u,(x, \varphi)) \right)\\
&=& Y\left(\varphi_1 \mbox{Res}_{x_1}\mbox{Res}_{x_0} \Bigl(
\frac{1}{2}l_t(x_1) \Bigl( \frac{\partial}{\partial \varphi_1} +
\varphi_1 \frac{\partial}{\partial x_1} \Bigr) + g_t(x_1) \Bigr)
\right. \\
& & \left. \hspace{1in} x^{-1} \delta \Bigl( \frac{x_1 - x_0 -
\varphi_1 \varphi}{x} \Bigr) Y(\tau,(x_0, \varphi_1 - \varphi))u,(x,
\varphi) \right) .  
\end{eqnarray*}
But
\begin{multline*}
\varphi_1 \mbox{Res}_{x_1}\mbox{Res}_{x_0} \left(
\frac{1}{2}l_t(x_1) \Bigl( \frac{\partial}{\partial \varphi_1} +
\varphi_1 \frac{\partial}{\partial x_1} \Bigr) + g_t(x_1) \right) \\
x^{-1} \delta \Bigl( \frac{x_1 - x_0 - \varphi_1
\varphi}{x} \Bigr) Y(\tau,(x_0, \varphi_1 - \varphi)) 
\end{multline*}
\begin{multline*}
=  \varphi_1 \mbox{Res}_{x_1}\mbox{Res}_{x_0} \left( \frac{1}{2}
l_t(x_1) \Bigl( \frac{\partial}{\partial \varphi_1} + \varphi_1
\frac{\partial}{\partial x_1} \Bigr) + g_t(x_1) \right) \\
x_1^{-1} \delta \Bigl( \frac{x + x_0 + \varphi_1
\varphi}{x_1} \Bigr) Y(\tau,(x_0, \varphi_1 - \varphi)) 
\end{multline*}
\begin{multline*}
= \varphi_1 \mbox{Res}_{x_1}\mbox{Res}_{x_0} \left(
\frac{1}{2}l_t(x_1) \Bigl( \frac{\partial}{\partial \varphi_1} 
x_1^{-1} \delta \Bigl( \frac{x + x_0 + \varphi_1 \varphi}{x_1} \Bigr)
\Bigr) Y(\tau,(x_0,- \varphi)) \right.\\ 
 + \; \frac{1}{2}l_t(x_1) x_1^{-1} \delta \Bigl(
\frac{x + x_0}{x_1} \Bigr) \frac{\partial}{\partial \varphi_1}
Y(\tau, (x_0,\varphi_1)) \\
\left. + \; g_t(x_1) x_1^{-1} \delta \Bigl( \frac{x
+ x_0}{x_1} \Bigr) Y(\tau,(x_0,- \varphi)) \right) 
\end{multline*}
\begin{multline*}
= \varphi_1 \mbox{Res}_{x_1} \mbox{Res}_{x_0} \biggl( 
\biggl(\sum_{n \in \mathbb{Z}} n (x + x_0)^{n-1} \varphi x_1^{-n-1}
\biggr) \biggl( \frac{1}{2} \sum_{j \in \Z} t^j B^{(0)}_j x_1^{-j + 1}
\biggr) \biggr.\\
\biggl( \sum_{k \in \mathbb{Z}} G(k + \frac{1}{2})
x_0^{-k - 2} \biggr) 
+ \biggl(\sum_{n \in \mathbb{Z}}(x + x_0)^n x_1^{-n-1} \biggr) 
\biggl( \frac{1}{2} \sum_{j \in \Z} t^j B^{(0)}_j x_1^{-j + 1} \biggr)\\
\biggl( 2\sum_{k \in \mathbb{Z}} L(k) x_0^{-k - 2} \biggr)
+ \biggl(\sum_{n \in \mathbb{Z}}(x + x_0)^n x_1^{-n-1} \biggr)
\biggl( \sum_{j \in \Z} t^{j-\frac{1}{2}} N^{(0)}_{j - \frac{1}{2}}
x_1^{-j + 1} \biggr) \\
\biggl. \biggl( \sum_{k \in \mathbb{Z}} G(k + \frac{1}{2})
x_0^{-k - 2} - 2\varphi L(k) x_0^{-k - 2} \biggr) \biggr)
\end{multline*}
\begin{multline*}
= \varphi_1 \mbox{Res}_{x_0} \biggl( \sum_{l \in \mathbb{N}} \sum_{k \in
\mathbb{Z}} \sum_{j \in \Z} \biggl( (-j + 1) \binom{-j}{l} x^{-j -l} x_0^l
\frac{\varphi}{2} t^j B^{(0)}_j G(k + \frac{1}{2}) x_0^{-k - 2}
\biggr. \biggr. \\
+ \; \binom{-j + 1}{l} x^{-j-l + 1} x_0^l t^j B^{(0)}_j
L(k) x_0^{-k - 2} \\
\biggl. \biggr. + \; \binom{-j + 1}{l} x^{-j-l + 1} x_0^l
t^{j-\frac{1}{2}} N^{(0)}_{j - \frac{1}{2}} \Bigl(G(k + \frac{1}{2}) -
2\varphi L(k) \Bigr) x_0^{-k - 2} \biggr) \biggr) 
\end{multline*}
\begin{multline*}
= \varphi_1 \sum_{l \in \mathbb{N}} \sum_{j \in \Z} \left( \frac{(-j +
1)}{2} \binom{-j}{l} x^{-j -l}  \varphi t^j
B^{(0)}_j G(l - \frac{1}{2}) \right. \\
+ \; \binom{-j + 1}{l} x^{-j-l + 1} t^j B^{(0)}_j
L(l - 1) \\
\left. + \; \binom{-j + 1}{l} x^{-j-l + 1} t^{j-\frac{1}{2}}
N^{(0)}_{j - \frac{1}{2}} (G(l - \frac{1}{2}) - 
2\varphi L(l - 1)) \right) 
\end{multline*}
\begin{multline*}
= \varphi_1 \sum_{m = -1}^\infty \sum_{j \in \Z} \left( \binom{-j + 1}{m+1} 
\left( t^j B^{(0)}_j L(m) + t^{j-\frac{1}{2}} N^{(0)}_{j - \frac{1}{2}} 
(G(m + \frac{1}{2}) \right) x^{-j-m} \right. \\
\left. + \; \varphi \binom{-j + 1}{m+1}\left( \frac{(-j -m )}{2} x^{-1}
t^j B^{(0)}_j G(m + \frac{1}{2}) + 2 t^{j-\frac{1}{2}} N^{(0)}_{j -
\frac{1}{2}} L(m) \right) x^{-j-m} \right)
\end{multline*}
\begin{multline*}
= \varphi_1  \sum_{m = -1}^{\infty} \sum_{j \in \Z} \binom{-j + 1}{m + 1}
x^{-j - m} \biggl( \Bigl( t^j B_j^{(0)} + 2\varphi t^{j - \frac{1}{2}} 
N_{j - \frac{1}{2}}^{(0)} \Bigr) L(m)  \\
+ \; \Bigl( t^{j - \frac{1}{2}}  N_{j - \frac{1}{2}}^{(0)} + \varphi
x^{-1} \frac{(-j-m)}{2} t^j B_j^{(0)} \Bigr) G(m + \frac{1}{2})
\biggr) .
\end{multline*}
This finishes the proof of (\ref{special case b bracket}) thus
concluding the proof of the sewing axiom for case (b).

(c) We prove the sewing axiom for the special case $Q_1 \; _2\infty_0
\; Q_2$ when 
\[Q_1 = ((z,\theta); \mathbf{0}, (a^{(1)}, A^{(1)}, M^{(1)}),
(a^{(2)}, A^{(2)}, M^{(2)})) \in SK(2)\] 
and 
\[Q_2 = ((B^{(0)}, N^{(0)}),(b^{(1)}, B^{(1)}, N^{(1)})) \in SK(1)
. \] 

Using steps (a) and (b) and Propositions \ref{sewing associativity},
\ref{contraction associativity} and \ref{doubly convergent}, we have 
\begin{eqnarray*}
\lefteqn{(\nu_2^Y(Q_1) \; _2*_0 \; \nu_1^Y(Q_2))
e^{-\Gamma(a^{(2)}, A^{(2)}, M^{(2)},B^{(0)}, N^{(0)})c}}\\
&=& \! (\nu_2^Y(((z,\theta); \mathbf{0}, (a^{(1)}, A^{(1)}, M^{(1)}),
(1,\mathbf{0})) \; _2\infty_0 \; (\mathbf{0},(a^{(2)}, A^{(2)},
M^{(2)}))) \; _2*_0 \; \\
& & \hspace{2.2in} \nu_1^Y(Q_2))  \; e^{-\Gamma(a^{(2)}, A^{(2)}, M^{(2)},B^{(0)},
N^{(0)})c} \\
&=& \! ((\nu_2^Y((z,\theta); \mathbf{0}, (a^{(1)}, A^{(1)}, M^{(1)}),
(1,\mathbf{0})) \; _2*_0 \; \nu_1^Y (\mathbf{0},(a^{(2)}, A^{(2)},
M^{(2)}))) \; _2*_0 \; \\
& & \hspace{2.2in} \nu_1^Y(Q_2)) \; e^{-\Gamma(a^{(2)}, A^{(2)}, M^{(2)},B^{(0)},
N^{(0)})c} \\
&=& \! (\nu_2^Y((z,\theta); \mathbf{0}, (a^{(1)}, A^{(1)}, M^{(1)}),
(1,\mathbf{0})) \; _2*_0 \; (\nu_1^Y (\mathbf{0},(a^{(2)}, A^{(2)}, 
M^{(2)})) \; _1*_0 \; \\
& & \hspace{2.2in} \nu_1^Y(Q_2))) \; e^{-\Gamma(a^{(2)}, A^{(2)}, M^{(2)},B^{(0)},
N^{(0)})c} \\
&=& \! \nu_2^Y((z,\theta); \mathbf{0}, (a^{(1)}, A^{(1)}, M^{(1)}),
(1,\mathbf{0})) \; _2*_0 \; \nu_1^Y ((\mathbf{0},(a^{(2)}, A^{(2)}, 
M^{(2)})) \; _1\infty_0 \; Q_2) \\
&=& \! \nu_2^Y(((z,\theta); \mathbf{0}, (a^{(1)}, A^{(1)}, M^{(1)}),
(1,\mathbf{0})) \; _2\infty_0 \; ((\mathbf{0},(a^{(2)}, A^{(2)}, 
M^{(2)})) \; _1\infty_0 \; Q_2)) \\
&=& \! \nu_2^Y((((z,\theta); \mathbf{0}, (a^{(1)}, A^{(1)}, M^{(1)}),
(1,\mathbf{0})) \; _2\infty_0 \; (\mathbf{0},(a^{(2)}, A^{(2)}, 
M^{(2)}))) \; _2\infty_0 \; Q_2) \\
&=& \! \nu_2^Y( Q_1 \; _2\infty_0 \; Q_2) 
\end{eqnarray*}

(d) We prove the sewing axiom for the special case  $Q_1 \; _1\infty_0
\; Q_2$ when 
\[Q_1 = ((z,\theta); \mathbf{0}, (1,\mathbf{0}), (1,\mathbf{0})) \in
SK(2)\]   
and 
\[Q_2 ((B^{(0)}, N^{(0)}),(b^{(1)}, B^{(1)}, N^{(1)})) \in SK(1)
. \]

In this case,
\begin{eqnarray*}
\lefteqn{\langle v', (\nu_2^Y(Q_1) \; _1*_0 \;
\nu_1^Y(Q_2))_{t^{1/2}} (v_1 \otimes v_2) \rangle} \\
&=&  \sum_{k \in \frac{1}{2}\mathbb{Z}} \sum_{i^{(k)} = 1}^{\dim
V_{(k)}} \langle v', Y(e_{i^{(k)}}^{(k)}, (x, \varphi)) v_2 \rangle
\cdot \langle e^{- \sum_{j \in \Z} \left(B^{(0)}_j L'(j) + N^{(0)}_{j
- \frac{1}{2}} G'(j - \frac{1}{2}) \right)} \cdot \\
& & \qquad (e^{(k)}_{i^{(k)}})^*, e^{- \sum_{j \in \Z} \left(B^{(1)}_j L(j) +
N^{(1)}_{j - \frac{1}{2}} G(j - \frac{1}{2}) \right)} \cdot \left.
(b^{(1)})^{-2L(0)} \cdot v_1 \rangle t^k \right|_{(x,\varphi) =
(z,\theta)} \\ 
&=&  \sum_{k \in \frac{1}{2}\mathbb{Z}} \sum_{i^{(k)} = 1}^{\dim
V_{(k)}} \langle v', Y(e_{i^{(k)}}^{(k)}, (x, \varphi))v_2 \rangle
\cdot \langle (e^{(k)}_{i^{(k)}})^*,  \\ 
& & \qquad e^{- \sum_{j \in \Z} \left(t^j B^{(0)}_j L(-j) + t^{j -
\frac{1}{2}} N^{(0)}_{j - \frac{1}{2}} G(- j + \frac{1}{2}) \right)}
t^{L(0)} \cdot \\
& & \hspace{1in} e^{- \sum_{j \in \Z} \left(B^{(1)}_j  L(j) +
N^{(1)}_{j - \frac{1}{2}} G(j - \frac{1}{2}) \right)} \cdot 
\left. (b^{(1)})^{-2L(0)} \cdot v_1 \rangle
\right|_{(x,\varphi) = (z,\theta)} \\ 
&=& \langle v', Y(e^{- \sum_{j \in \Z} \left( t^j B^{(0)}_j L(-j) +
t^{j - \frac{1}{2}} N^{(0)}_{j - \frac{1}{2}} G(- j + \frac{1}{2})
\right)} t^{L(0)} \cdot \\
& & \qquad e^{- \sum_{j \in \Z} \left(B^{(1)}_j L(j) +
N^{(1)}_{j - \frac{1}{2}} G(j - \frac{1}{2}) \right)} \cdot \left.
(b^{(1)})^{-2L(0)} \cdot v_1, (x, \varphi)) \cdot v_2 \rangle
\right|_{(x,\varphi) = (z,\theta)}   \\
&=& (-1)^{\eta(v_1)\eta(v_2)} \langle v', e^{xL(-1) + \varphi
G(-\frac{1}{2})} Y(v_2, (-x,-\varphi)) \\
& & \qquad e^{- \sum_{j \in \Z} \left(t^j B^{(0)}_j L(-j) +
t^{j - \frac{1}{2}} N^{(0)}_{j - \frac{1}{2}} G(- j + \frac{1}{2})
\right)} t^{L(0)} \cdot \\
& & \hspace{1in} e^{- \sum_{j \in \Z} \left(B^{(1)}_j L(j) +
N^{(1)}_{j - \frac{1}{2}} G(j - \frac{1}{2}) \right)} \cdot 
\left. (b^{(1)})^{-2L(0)} \cdot v_1 \rangle
\right|_{(x,\varphi) = (z, \theta)} 
\end{eqnarray*}
\begin{multline}\label{in case d} 
= (-1)^{\eta(v_1)\eta(v_2)} \langle e^{xL'(1) + \varphi
G'(\frac{1}{2})} v', Y(v_2, (-x,-\varphi)) \\
e^{- \sum_{j \in \Z} \left(t^j B^{(0)}_j L(-j) +
t^{j - \frac{1}{2}} N^{(0)}_{j - \frac{1}{2}} G(- j + \frac{1}{2})
\right)} t^{L(0)} \cdot \\
e^{- \sum_{j \in \Z} \left(B^{(1)}_j L(j) +
N^{(1)}_{j - \frac{1}{2}} G(j - \frac{1}{2}) \right)} \cdot 
\left. (b^{(1)})^{-2L(0)} \cdot v_1 \rangle
\right|_{(x,\varphi) = (z, \theta)} .
\end{multline}
By (\ref{for case d}) in case (b), the right-hand side of (\ref{in
case d}) is equal to
\begin{multline*}
(-1)^{\eta(v_1)\eta(v_2)} (\nu_2^Y ( (-z,-\theta) ; \mathbf{0},  (1,
\mathbf{0}), (1, \mathbf{0})) \\
\; _2*_0 \; \nu_1^Y ((B^{(0)}, N^{(0)}), (b^{(1)},
B^{(1)}, N^{(1)})) )_{t^{1/2}} ( e^{zL'(1) + \theta G'(\frac{1}{2})} v'
\otimes v_2 \otimes v_1 ) . 
\end{multline*} 
Let $\Theta^{(2)}_j$, for $j \in \frac{1}{2} \Z$, be defined by
(\ref{Theta in case b}), and let $H_{t^{1/2}} (x, \varphi) =
(H_{t^{1/2}}^0 (x, \varphi), H_{t^{1/2}}^1 (x, \varphi))$ be defined by
(\ref{Ht in case b}).  Then by case (b), we have  
\[\langle v', (\nu_2^Y(Q_1) \; _1*_0 \;
\nu_1^Y(Q_2))_{t^{1/2}} (v_1 \otimes v_2) \rangle \hspace{2.4in}\]
\begin{multline*}
= (-1)^{\eta(v_1)\eta(v_2)} \langle e^{- \sum_{j \in \Z} \left(t^j
B^{(0)}_j L'(j) + t^{j - \frac{1}{2}} N^{(0)}_{j - \frac{1}{2}} G'(j -
\frac{1}{2}) \right)} e^{-xL'(1) - \varphi G'(\frac{1}{2})} v', \\
Y(e^{- \sum_{j \in \Z} \left(\Theta^{(2)}_j L(j) +
\Theta^{(2)}_{j - \frac{1}{2}} G(j - \frac{1}{2}) \right)} e^{- 2\Theta^{(2)}_0
L(0)} v_2, H_{t^{1/2}}^{-1} \circ (x,\varphi) ) \\
t^{L(0)} e^{- \sum_{j \in \Z} \left(B^{(1)}_j L(j) +  N^{(1)}_{j - \frac{1}{2}} 
G(j - \frac{1}{2}) \right)} \left. (b^{(1)}_0)^{-L(0)} \cdot v_1 \rangle 
\right|_{(x,\varphi) = (- z, - \theta)} 
\end{multline*}
\begin{multline*}
= \langle e^{\left( H_{t^{1/2}}^0(x, \varphi) L'(1) + H_{t^{1/2}}^1
(x, \varphi) G'(\frac{1}{2}) \right)} e^{- \sum_{j \in \Z} \left(t^j
B^{(0)}_j L'(j) + t^{j - \frac{1}{2}} N^{(0)}_{j - \frac{1}{2}} G'(j - 
\frac{1}{2}) \right)} \\
e^{-xL'(1) - \varphi G'(\frac{1}{2})} v', 
Y(t^{L(0)} e^{- \sum_{j \in \Z} \left(B^{(1)}_j
L(j) +  N^{(1)}_{j - \frac{1}{2}} G(j - \frac{1}{2}) \right)}
(b^{(1)})^{-2L(0)} \cdot v_1, \\
- H_{t^{1/2}}^{-1} \circ I (x,\varphi)) 
e^{- \sum_{j \in \Z} \left(\Theta^{(2)}_j L(j) + \Theta^{(2)}_{j - \frac{1}{2}} 
G(j - \frac{1}{2}) \right)} \left. e^{- 2\Theta^{(2)}_0 L(0)} v_2 \rangle 
\right|_{(x,\varphi) = (- z, - \theta)} .
\end{multline*}
The right-hand side of the above expression converges when $t^{1/2} = 1$ and,
by the definition of sewing, is equal to 
\[\langle v', \nu_2^Y( Q_1 \; _1\infty_0 \; Q_2)(v_1 \otimes v_2)
\rangle.\] 

(e) We prove the sewing axiom for the special case  $Q_1 \; _1\infty_0
\; Q_2$ when 
\[Q_1 = ((z,\theta); \mathbf{0}, (a^{(1)}, A^{(1)}, M^{(1)}),
(a^{(2)}, A^{(2)}, M^{(2)})) \in SK(2)\] 
and 
\[Q_2 ((B^{(0)}, N^{(0)}),(b^{(1)}, B^{(1)}, N^{(1)})) \in SK(1) . \]

Using cases (a), (c) and (d) and Propositions \ref{sewing
associativity}, \ref{contraction associativity} and \ref{doubly
convergent}, and letting 
\[e^{- \Gamma(a^{(1)}, A^{(1)}, M^{(1)}, B^{(0)},
N^{(0)})c} = e^{- \Gamma c} , \]
we have 
\begin{eqnarray*}
\lefteqn{(\nu_2^Y(Q_1)  \; _1*_0 \; \nu_1^Y(Q_2) e^{-
\Gamma(a^{(1)}, A^{(1)}, M^{(1)}, B^{(0)}, N^{(0)})c} } \\
&=& (\nu_2^Y((((z,\theta); \mathbf{0}, (1, \mathbf{0}), (1,
\mathbf{0})) \; _1\infty_0 \; (\mathbf{0}, (a^{(1)}, A^{(1)},
M^{(1)}))) \; _2\infty_0  \\ 
& & \quad  (\mathbf{0}, (a^{(2)}, A^{(2)},
M^{(2)}))) \; _1*_0 \; \nu_1^Y((B^{(0)}, N^{(0)}),(b^{(1)},
B^{(1)}, N^{(1)})))  e^{- \Gamma c}  \\ 
&=& (((\nu_2^Y((z,\theta); \mathbf{0}, (1, \mathbf{0}), (1,
\mathbf{0})) \; _1*_0 \; \nu_1^Y(\mathbf{0}, (a^{(1)}, A^{(1)},
M^{(1)}))) \;  _2*_0 \\
& & \quad \nu_1^Y(\mathbf{0}, (a^{(2)}, A^{(2)},
M^{(2)}))) \; _1*_0 \; \nu_1^Y((B^{(0)}, N^{(0)}),(b^{(1)},
B^{(1)}, N^{(1)}))) e^{- \Gamma c}  \\ 
&=& ((\nu_2^Y((z,\theta); \mathbf{0}, (1, \mathbf{0}), (1,
\mathbf{0})) \; _1*_0 \; (\nu_1^Y(\mathbf{0}, (a^{(1)}, A^{(1)},
M^{(1)})) \; _1*_0 \\
& & \quad  \nu_1^Y((B^{(0)}, N^{(0)}),(b^{(1)},
B^{(1)}, N^{(1)})))) \; _2*_0 \; \nu_1^Y(\mathbf{0}, (a^{(2)},
A^{(2)}, M^{(2)})))  e^{- \Gamma c}  \\
&=& (\nu_2^Y((z,\theta); \mathbf{0}, (1, \mathbf{0}), (1, \mathbf{0}))
\; _1*_0 \; \nu_1^Y((\mathbf{0}, (a^{(1)}, A^{(1)}, M^{(1)})) \;
_1\infty_0 \;\\ 
& &  \quad  ((B^{(0)}, N^{(0)}),(b^{(1)}, B^{(1)}, N^{(1)})))) \; _2*_0 \;
\nu_1^Y(\mathbf{0}, (a^{(2)}, A^{(2)}, M^{(2)})) \\ 
&=& \nu_2^Y(((z,\theta); \mathbf{0}, (1, \mathbf{0}), (1, \mathbf{0}))
\; _1\infty_0 \; ((\mathbf{0}, (a^{(1)}, A^{(1)}, M^{(1)})) \; _1\infty_0
\\ 
& &  \quad ((B^{(0)}, N^{(0)}),(b^{(1)}, B^{(1)}, N^{(1)})))) \; _2*_0 \;
\nu_1^Y(\mathbf{0}, (a^{(2)}, A^{(2)}, M^{(2)})) \\ 
&=& \nu_2^Y((((z,\theta); \mathbf{0}, (1, \mathbf{0}), (1,
\mathbf{0})) \; _1\infty_0 \; ((\mathbf{0}, (a^{(1)}, A^{(1)},
M^{(1)})) \; _1\infty_0 \\ 
& & \quad ((B^{(0)}, N^{(0)}),(b^{(1)}, B^{(1)}, N^{(1)})))) \;_2\infty_0
\; (\mathbf{0}, (a^{(2)}, A^{(2)}, M^{(2)}))) \\ 
&=& \nu_2^Y(((((z,\theta); \mathbf{0}, (1, \mathbf{0}), (1,
\mathbf{0})) \; _1\infty_0 \; (\mathbf{0}, (a^{(1)}, A^{(1)},
M^{(1)}))) \; _2\infty_0  \\
& & \quad (\mathbf{0}, (a^{(2)}, A^{(2)},
M^{(2)}))) \; _1\infty_0 \; ((B^{(0)}, N^{(0)}),(b^{(1)},
B^{(1)}, N^{(1)})))  \\ 
&=& \nu_2^Y(Q_1 \; _1\infty_0 \; Q_2)  
\end{eqnarray*}

(f) We prove the sewing axiom for the special case $Q_1 \; _m\infty_0
\; Q_2$ when
\begin{multline*}
Q_1 = ((z_1, \theta_1),...,(z_{m-1}, \theta_{m-1});
(A^{(0)},M^{(0)}), (a^{(1)}, A^{(1)}, M^{(1)}),...,\\
(a^{(m)}, A^{(m)}, M^{(m)})) \in SK(m),  
\end{multline*} 
and
\[Q_2 = ((z, \theta); \mathbf{0}, (b^{(1)}, B^{(1)}, N^{(1)}),
(b^{(2)}, B^{(2)}, N^{(2)})) \in SK(2) .\]
In this case,
\[\langle v', (\nu_m^Y(Q_1) \; _m*_0 \; \nu_2^Y(Q_1))_{t^{1/2}} (v_1
\otimes \cdots v_{m +1}) \rangle \hspace{2.7in} \]
\begin{multline*}
= \sum_{k \in \frac{1}{2}\mathbb{Z}} \sum_{i^{(k)} = 1}^{\dim V_{(k)}}
\iota^{-1}_{1 \cdots m-1} \langle e^{- \sum_{j \in \Z} \left(
A^{(0)}_j L'(j) + M^{(0)}_{j - \frac{1}{2}} G'(j - \frac{1}{2})
\right)} v', \\
Y(e^{- \sum_{j \in \Z} \left(A^{(1)}_j L(j) +  M^{(1)}_{j
- \frac{1}{2}} G(j - \frac{1}{2}) \right)} \cdot (a^{(1)})^{-2L(0)}
\cdot v_1, (x_1, \varphi_1)) \cdots \\
Y(e^{- \sum_{j \in \Z} \left(A^{(m-1)}_j L(j) + M^{(m-1)}_{j -
\frac{1}{2}} G(j - \frac{1}{2}) \right)} \cdot (a^{(m-1)})^{-2L(0)}
\cdot v_{m-1}, \\
(x_{m-1}, \varphi_{m-1})) \cdot e^{- \sum_{j \in \Z} \left(A^{(m)}_j L(j) +
M^{(m)}_{j - \frac{1}{2}} G(j - \frac{1}{2}) \right)} \cdot
(a^{(m)})^{-2L(0)} \cdot  e_{i^{(k)}}^{(k)} \rangle \cdot \\
\langle (e^{(k)}_{i^{(k)}})^*, Y(e^{- \sum_{j \in \Z}
\left(B^{(1)}_j L(j) + N^{(1)}_{j - \frac{1}{2}} G(j - \frac{1}{2})
\right)} \cdot (b^{(1)})^{-2L(0)} \cdot v_m, (x,\varphi)) \cdot \\
e^{- \sum_{j \in \Z} \left(B^{(2)}_j L(j)
+ N^{(2)}_{j - \frac{1}{2}} G(j - \frac{1}{2}) \right)} \cdot \left.
(b^{(2)})^{-2L(0)} \cdot v_{m+1}\rangle t^k \right|_{\begin{scriptsize}
\begin{array}{c}
 (x_i,\varphi_i)
= (z_i,\theta_i),\\
(x,\varphi) = (z,\theta)  \end{array} \end{scriptsize}}
\end{multline*}
\begin{multline*} 
=  \sum_{k \in \frac{1}{2}\mathbb{Z}} \sum_{i^{(k)} = 1}^{\dim V_{(k)}}
\iota^{-1}_{1 \cdots m-1} \langle e^{- \sum_{j \in \Z} \left(
A^{(0)}_j L'(j) + M^{(0)}_{j - \frac{1}{2}} G'(j - \frac{1}{2})
\right)} v', \\
Y(e^{- \sum_{j \in \Z} \left(A^{(1)}_j L(j) +  M^{(1)}_{j -
\frac{1}{2}} G(j - \frac{1}{2}) \right)} \cdot (a^{(1)})^{-2L(0)}
\cdot v_1, (x_1, \varphi_1)) \cdots \\
Y(e^{- \sum_{j \in \Z} \left(A^{(m-1)}_j L(j) + M^{(m-1)}_{j -
\frac{1}{2}} G(j - \frac{1}{2}) \right)} \cdot (a^{(m-1)})^{-2L(0)}
\cdot v_{m-1}, \\
(x_{m-1}, \varphi_{m-1})) \cdot e^{- \sum_{j \in \Z} \left(A^{(m)}_j L(j) + 
M^{(m)}_{j - \frac{1}{2}} G(j - \frac{1}{2}) \right)} \cdot
(a^{(m)})^{-2L(0)} \cdot e_{i^{(k)}}^{(k)} \rangle \cdot \\ 
\langle t^{L'(0)} (e^{(k)}_{i^{(k)}})^*, Y(e^{- \sum_{j \in \Z}
\left(B^{(1)}_j L(j) + N^{(1)}_{j - \frac{1}{2}} G(j - \frac{1}{2})
\right)} \cdot (b^{(1)})^{-2L(0)} \cdot v_m,  \\
(x,\varphi)) \cdot e^{- \sum_{j \in \Z} \left(B^{(2)}_j L(j)
+ N^{(2)}_{j - \frac{1}{2}} G(j - \frac{1}{2}) \right)} \cdot \\
\left.(b^{(2)})^{-2L(0)} \cdot v_{m+1}\rangle \right|_{\begin{scriptsize}
\begin{array}{c}
 (x_i,\varphi_i)
= (z_i,\theta_i),\\
(x,\varphi) = (z,\theta)  \end{array} \end{scriptsize}}
\end{multline*}
\begin{multline*} 
=  \sum_{k \in \frac{1}{2}\mathbb{Z}} \sum_{i^{(k)} = 1}^{\dim V_{(k)}}
\iota^{-1}_{1 \cdots m-1} \langle e^{- \sum_{j \in \Z} \left(
A^{(0)}_j L'(j) + M^{(0)}_{j - \frac{1}{2}} G'(j - \frac{1}{2})
\right)} v', \\
Y(e^{- \sum_{j \in \Z} \left(A^{(1)}_j L(j) +  M^{(1)}_{j -
\frac{1}{2}} G(j - \frac{1}{2}) \right)} \cdot (a^{(1)})^{-2L(0)}
\cdot v_1, (x_1, \varphi_1)) \cdots \\
Y(e^{- \sum_{j \in \Z} \left(A^{(m-1)}_j L(j) + M^{(m-1)}_{j -
\frac{1}{2}} G(j - \frac{1}{2}) \right)} \cdot (a^{(m-1)})^{-2L(0)}
\cdot v_{m-1}, \\
(x_{m-1}, \varphi_{m-1})) \cdot e^{- \sum_{j \in \Z} \left(A^{(m)}_j 
L(j) + M^{(m)}_{j -\frac{1}{2}} G(j - \frac{1}{2}) \right)} \cdot 
(a^{(m)})^{-2L(0)}\cdot e_{i^{(k)}}^{(k)} \rangle \cdot \\
\langle (e^{(k)}_{i^{(k)}})^*,  t^{L(0)} Y(e^{- \sum_{j \in \Z}
\left(B^{(1)}_j L(j) + N^{(1)}_{j - \frac{1}{2}} G(j - \frac{1}{2})
\right)} \cdot (b^{(1)})^{-2L(0)} \cdot v_m, \\
(x,\varphi)) \cdot e^{- \sum_{j \in \Z} \left(B^{(2)}_j L(j)
+ N^{(2)}_{j - \frac{1}{2}} G(j - \frac{1}{2}) \right)} \cdot \\
\left. (b^{(2)})^{-2L(0)} \cdot v_{m+1}\rangle \right|_{\begin{scriptsize}
\begin{array}{c}
 (x_i,\varphi_i) = (z_i,\theta_i),\\
(x,\varphi) = (z,\theta)  \end{array} \end{scriptsize}}
\end{multline*}
\begin{multline}\label{case f} 
=  \iota^{-1}_{1 \cdots m} \langle e^{- \sum_{j \in \Z} \left(
A^{(0)}_j L'(j) + M^{(0)}_{j - \frac{1}{2}} G'(j - \frac{1}{2})
\right)} v', \\
Y(e^{- \sum_{j \in \Z} \left(A^{(1)}_j L(j) +  M^{(1)}_{j -
\frac{1}{2}} G(j - \frac{1}{2}) \right)} \cdot (a^{(1)})^{-2L(0)}  
\cdot v_1, (x_1, \varphi_1)) \cdots \\
Y(e^{-\sum_{j \in \Z} \left(A^{(m-1)}_j L(j) + M^{(m-1)}_{j -
\frac{1}{2}} G(j - \frac{1}{2}) \right)} \cdot (a^{(m-1)})^{-2L(0)}
\cdot v_{m-1}, \\
(x_{m-1}, \varphi_{m-1})) \cdot e^{- \sum_{j \in \Z} \left(A^{(m)}_j 
L(j) + M^{(m)}_{j - \frac{1}{2}} G(j - \frac{1}{2}) \right)} \cdot
(t^{-\frac{1}{2}}a^{(m)})^{-2L(0)} \cdot \\
Y(e^{- \sum_{j \in \Z} \left(B^{(1)}_j L(j) + N^{(1)}_{j -
\frac{1}{2}} G(j - \frac{1}{2}) \right)} \cdot (b^{(1)})^{-2L(0)}
\cdot v_m, (x,\varphi)) \cdot \\
e^{- \sum_{j \in \Z} \left(B^{(2)}_j
L(j) + N^{(2)}_{j - \frac{1}{2}} G(j - \frac{1}{2}) \right)} \cdot 
\left. (b^{(2)})^{-2L(0)} \cdot v_{m+1}\rangle \right|_{\begin{scriptsize}
\begin{array}{c}
 (x_i,\varphi_i)
= (z_i,\theta_i),\\
(x,\varphi) = (z,\theta)  \end{array} \end{scriptsize}} .
\end{multline} 
Since $V$ is a positive energy representation for $\mathfrak{ns}$ and by
the truncation condition (\ref{truncation}) for $Y(\cdot, (x,\varphi))$,
the right-hand side of (\ref{case f}) is a polynomial in $t^{1/2}$ and
thus convergent when $t^{1/2} = 1$.  If we can prove 
\begin{multline*}
e^{- \sum_{j \in \Z} \left( A_j^{(m)} L(j) + M_{j -
\frac{1}{2}}^{(m)} G(j - \frac{1}{2}) \right)} \cdot
(t^{-\frac{1}{2}} a^{(m)})^{-2L(0)} Y(u,(x,\varphi)) \cdot
(t^{-\frac{1}{2}} a^{(m)})^{2L(0)} \cdot \\
e^{\sum_{j \in \Z} \left( A_j^{(m)} L(j) + M_{j -
\frac{1}{2}}^{(m)} G(j - \frac{1}{2}) \right)} 
\end{multline*} 
\begin{equation}\label{Theta for case f}
= Y\biggl( (t^{-\frac{1}{2}} a^{(m)})^{-2L(0)} e^{- \sum_{j \in \Z}
\left( \Theta^{(1)}_j L(j) + \Theta^{(1)}_{j - \frac{1}{2}} G(j - \frac{1}{2})
\right)} \cdot e^{-2\Theta^{(1)}_0 L(0)} u, H_{t^{1/2}}^{-1}(x,\varphi)
\biggr)  
\end{equation}
where 
\begin{multline*}
H_{t^{1/2}}^{-1}(x,\varphi) \\
= (t^{-\frac{1}{2}} a^{(m)})^{2L_0(x,\varphi)} \cdot
\exp \Bigl(\sum_{j \in \Z} \Bigl(A_j^{(m)} L_j(x,\varphi) + M_{j -
\frac{1}{2}}^{(m)} G_{j - \frac{1}{2}}(x,\varphi) \Bigr) \Bigr) \cdot
(x,\varphi) 
\end{multline*} 
and the $\Theta^{(1)}_j = \Theta^{(1)}_j (t^{-1/2} \asqrt, A^{(m)},
M^{(m)}, (x,\varphi))$ are defined by (\ref{define first Theta}), then
since $e^{\Gamma(a^{(m)},A^{(m)},M^{(m)},\mathbf{0})c} = 1$, when $t^{1/2}
= 1$, the right-hand side of (\ref{case f}) is equal to  
\[\langle v', \nu_{m + 1}^Y(Q_1 \; _m\infty_0 \; Q_2) (v_1 \otimes
\cdots v_{m + 1}) \rangle .\] 
Thus in this case, the proof of the sewing axiom is reduced to the
proof of (\ref{Theta for case f}).  

Let 
\[H_{t^{1/2}}^{-1} (x,\varphi) = (\tilde{x} (t^\frac{1}{2}),
\tilde{\varphi}(t^\frac{1}{2})) .\]
Then by (\ref{conjugate by L(0)}), (\ref{exponential L(-1) and G(-1/2)
property}) and Corollary \ref{first Theta identity in End}, the
right-hand side of (\ref{Theta for case f}) is equal to
\begin{multline}
(t^{-\frac{1}{2}} a^{(m)})^{-2L(0)} Y\biggl( e^{- \sum_{j \in \Z}
\left( \Theta^{(1)}_j L(j) + \Theta^{(1)}_{j - \frac{1}{2}} G(j - \frac{1}{2})
\right)} \cdot e^{-2\Theta^{(1)}_0 L(0)} u, \biggr. \\
\biggl. ((t^{-\frac{1}{2}}
a^{(m)})^{2} \tilde{x} (t^\frac{1}{2}), (t^{-\frac{1}{2}}
a^{(m)}) \tilde{\varphi}(t^\frac{1}{2})\biggr) 
(t^{-\frac{1}{2}} a^{(m)})^{2L(0)} 
\end{multline}
\begin{multline*}
= \; (t^{-\frac{1}{2}} a^{(m)})^{-2L(0)}  Y \Bigl(
e^{((t^{-\frac{1}{2}} a^{(m)})^{2}\tilde{x}(t^\frac{1}{2}) - x) L(-1) +
((t^{-\frac{1}{2}} a^{(m)}) \tilde{\varphi} (t^\frac{1}{2}) - \varphi)
G(-\frac{1}{2})} \cdot  \\
 e^{- \sum_{j \in \Z} \left( \Theta^{(2)}_j L(j) + \Theta^{(2)}_{j -
\frac{1}{2}} G(j - \frac{1}{2}) \right)} \cdot  
e^{- 2\Theta^{(2)}_0 L(0)} u, (x,\varphi)
\Bigr) (t^{-\frac{1}{2}} a^{(m)})^{2L(0)} 
\end{multline*}
\begin{multline*}
= \; (t^{-\frac{1}{2}} a^{(m)})^{-2L(0)} Y \Biggl(\exp \Biggl(\! - \! \!
\sum_{k = -1}^{\infty} \sum_{j \in \Z}
\binom{j+1}{k+1} t^j (a^{(m)})^{-2j} x^{j - k}  \\
\biggl(\!  \Bigl( A^{(m)}_j + 2\left(\frac{j-k}{j+1} \right)
t^{-\frac{1}{2}} a^{(m)} x^{-1} \varphi M^{(m)}_{j - \frac{1}{2}}
\Bigr)  L(k) \\
+ \; x^{-1} \Bigl( \left(\frac{j-k}{j+1} \right) t^{-\frac{1}{2}}
a^{(m)} M^{(m)}_{j - \frac{1}{2}} + \varphi \frac{(j-k)}{2}
A^{(m)}_j \Bigr)  G(k + \frac{1}{2}) \biggl) \Biggr)  u, (x,\varphi)
\Biggr) 
\cdot \\ 
(t^{-\frac{1}{2}} a^{(m)})^{2L(0)} .
\end{multline*} 

But since 
\begin{multline*} 
e^{-  \sum_{j \in \Z}
\left( A_j L(j) + AM{j-\frac{1}{2}} G(j - \frac{1}{2}) \right) }
(t^{-\frac{1}{2}} a)^{-2L(0)}  \\
= (t^{-\frac{1}{2}} a)^{-2L(0)} e^{ \sum_{j \in \Z}
\left( t^{j} a^{-2j} A_j L(j) + t^{j - \frac{1}{2}}
a^{-2j + 1} M_{j-\frac{1}{2}} G(j - \frac{1}{2}) \right) }, 
\end{multline*}
proving (\ref{Theta for case f}) is equivalent to proving
\begin{multline*}
e^{- \sum_{j \in \Z} \left( t^{j}
(a^{(m)})^{-2j} A_j^{(m)} L(j) +  t^{j - \frac{1}{2}} (a^{(m)})^{-2j + 1}
M_{j - \frac{1}{2}}^{(m)} G(j - \frac{1}{2}) \right)}  Y(u,(x,\varphi))
\cdot  \\
e^{\sum_{j \in \Z} \left( t^{j}(a^{(m)})^{-2j} A_j^{(m)} L(j) + 
t^{j - \frac{1}{2}} (a^{(m)})^{-2j + 1} M_{j - \frac{1}{2}}^{(m)} G(j -
\frac{1}{2}) \right)}  
\end{multline*} 
\begin{multline*}
= \; Y \Biggl(\exp \Biggl(\! - \! \!
\sum_{k = -1}^{\infty} \sum_{j \in \Z}
\binom{j+1}{k+1} t^j (a^{(m)})^{-2j} x^{j - k}  \\
\biggl(\!  \Bigl( A^{(m)}_j + 2\left(\frac{j-k}{j+1} \right)
t^{-\frac{1}{2}} a^{(m)} x^{-1} \varphi M^{(m)}_{j - \frac{1}{2}}
\Bigr)  L(k) \\
+ \; x^{-1} \Bigl( \left(\frac{j-k}{j+1} \right) t^{-\frac{1}{2}}
a^{(m)} M^{(m)}_{j - \frac{1}{2}} + \varphi \frac{(j-k)}{2}
A^{(m)}_j \Bigr)  G(k + \frac{1}{2}) \biggl) \Biggr)  u, (x,\varphi)
\Biggr) 
\end{multline*}

By this equality, the proof of (\ref{new sewing coordinates}) is
reduced to the proof of the following bracket formula
\begin{equation*}
\Bigl[ \sum_{j \in \Z} \left( t^{j} (a^{(m)})^{-2j} A_j^{(m)} L(j) +  t^{j
- \frac{1}{2}} (a^{(m)})^{-2j + 1} M_{j - \frac{1}{2}}^{(m)} G(j -
\frac{1}{2}) \right), Y(u,(x,\varphi))
\Bigr] \hspace{.6in}
\end{equation*}
\begin{multline}\label{special case f bracket}
= Y \Biggl( \sum_{k = -1}^{\infty} \sum_{j \in \Z}
\binom{j+1}{k+1} t^j (a^{(m)})^{-2j} x^{j - k}  \\
\biggl(\!  \Bigl( A^{(m)}_j + 2\left(\frac{j-k}{j+1} \right)
t^{-\frac{1}{2}} a^{(m)} x^{-1} \varphi M^{(m)}_{j - \frac{1}{2}}
\Bigr)  L(k) \\
+ \; x^{-1} \Bigl( \left(\frac{j-k}{j+1} \right) t^{-\frac{1}{2}}
a^{(m)} M^{(m)}_{j - \frac{1}{2}} + \varphi \frac{(j-k)}{2}
A^{(m)}_j \Bigr)  G(k + \frac{1}{2})  \biggr)  u, (x,\varphi)
\Biggr) 
\end{multline} 

{}From (\ref{bracket relation for a vosa}), we have
\begin{multline*} 
\left[ Y(\tau, (x_1,\varphi_1)), Y(u,(x,\varphi)) \right] \\
= \mbox{Res}_{x_0} x^{-1} \delta \left( \frac{x_1 - x_0 -
\varphi_1 \varphi}{x} \right) Y(Y(\tau,(x_0, \varphi_1 -
\varphi))u,(x, \varphi)) .
\end{multline*} 
Let 
\begin{eqnarray*}
l_t(x_1) &=& \sum_{j \in \Z} t^{j} (a^{(m)})^{-2j} A_j^{(m)} x_1^{j + 1}
\\
g_t(x_1) &=& \sum_{j \in \Z} t^{j - \frac{1}{2}}
(a^{(m)})^{-2j + 1} M_{j - \frac{1}{2}}^{(m)} x_1^{j} .
\end{eqnarray*}
Then 
\begin{eqnarray*}
\varphi_1 \mbox{Res}_{x_1} l_t(x_1) Y(\omega, (x_1,\varphi_1)) &=&
\varphi_1 \sum_{j \in \Z} t^j (a^{(m)})^{-2j} A_j^{(m)}  L(j) \\ 
\varphi_1 \mbox{Res}_{x_1} g_t(x_1) Y(\tau, (x_1,\varphi_1)) &=&
\varphi_1 \sum_{j \in \Z} t^{j - \frac{1}{2}} (a^{(m)})^{-2j + 1} M_{j -
\frac{1}{2}}^{(m)} G(j - \frac{1}{2}) ,
\end{eqnarray*}
and by the $G(-1/2)$-derivative property
(\ref{G(-1/2)-derivative}), we have 
\begin{eqnarray*}
Y(\omega, (x_1,\varphi_1)) &=& Y( \frac{1}{2} G(-\frac{1}{2}) \tau,
(x_1, \varphi_1)) \\
&=& \frac{1}{2} \Bigl( \frac{\partial}{\partial \varphi_1} + \varphi_1
\frac{\partial}{\partial x_1} \Bigr) Y(\tau, (x_1,\varphi_1))  .
\end{eqnarray*}
Thus 
\begin{eqnarray*}
& & \hspace{-.5in} \varphi_1 \biggl[ \sum_{j \in \Z} \left( t^j
(a^{(m)})^{-2j} A_j^{(m)} L(j) + t^{j - \frac{1}{2}} (a^{(m)})^{-2j + 1}
M_{j - \frac{1}{2}}^{(m)} G(j - \frac{1}{2}) \right) , Y(u,(x,\varphi))
\biggr] \\ 
&=& \varphi_1 \mbox{Res}_{x_1} \left[ \Bigl( \frac{1}{2}l_t(x_1)
\Bigl( \frac{\partial}{\partial \varphi_1} + \varphi_1
\frac{\partial}{\partial x_1} \Bigr) + g_t(x_1) \Bigr) Y(\tau,
(x_1,\varphi_1)), Y(u,(x,\varphi)) \right] \\
&=& \varphi_1 \mbox{Res}_{x_1} \left( \Bigl( \frac{1}{2}l_t(x_1) \Bigl(
\frac{\partial}{\partial \varphi_1} + \varphi_1
\frac{\partial}{\partial x_1} \Bigr) + g_t(x_1) \Bigr) \right. \\
& & \hspace{.7in} \left. \mbox{Res}_{x_0} x^{-1} \delta \Bigl(
\frac{x_1 - x_0 - \varphi_1 \varphi}{x} \Bigr) Y(Y(\tau,(x_0,
\varphi_1 - \varphi))u,(x, \varphi)) \right)\\
&=& Y\left(\varphi_1 \mbox{Res}_{x_1}\mbox{Res}_{x_0} \Bigl(
\frac{1}{2}l_t(x_1) \Bigl( \frac{\partial}{\partial \varphi_1} +
\varphi_1 \frac{\partial}{\partial x_1} \Bigr) + g_t(x_1) \Bigr)
\right. \\
& & \left. \hspace{1in} x^{-1} \delta \Bigl( \frac{x_1 - x_0 -
\varphi_1 \varphi}{x} \Bigr) Y(\tau,(x_0, \varphi_1 - \varphi))u,(x,
\varphi) \right) .  
\end{eqnarray*}
But
\begin{multline*}
\varphi_1 \mbox{Res}_{x_1}\mbox{Res}_{x_0} \left(
\frac{1}{2}l_t(x_1) \Bigl( \frac{\partial}{\partial \varphi_1} +
\varphi_1 \frac{\partial}{\partial x_1} \Bigr) + g_t(x_1) \right) \\
x^{-1} \delta \Bigl( \frac{x_1 - x_0 - \varphi_1
\varphi}{x} \Bigr) Y(\tau,(x_0, \varphi_1 - \varphi)) 
\end{multline*}
\begin{multline*}
=  \varphi_1 \mbox{Res}_{x_1}\mbox{Res}_{x_0} \left( \frac{1}{2}
l_t(x_1) \Bigl( \frac{\partial}{\partial \varphi_1} + \varphi_1
\frac{\partial}{\partial x_1} \Bigr) + g_t(x_1) \right) \\
x_1^{-1} \delta \Bigl( \frac{x + x_0 + \varphi_1
\varphi}{x_1} \Bigr) Y(\tau,(x_0, \varphi_1 - \varphi)) 
\end{multline*}
\begin{multline*}
= \varphi_1 \mbox{Res}_{x_1}\mbox{Res}_{x_0} \left(
\frac{1}{2}l_t(x_1) \Bigl( \frac{\partial}{\partial \varphi_1} 
x_1^{-1} \delta \Bigl( \frac{x + x_0 + \varphi_1 \varphi}{x_1} \Bigr)
\Bigr) Y(\tau,(x_0,- \varphi)) \right.\\ 
 + \; \frac{1}{2}l_t(x_1) x_1^{-1} \delta \Bigl(
\frac{x + x_0}{x_1} \Bigr) \frac{\partial}{\partial \varphi_1}
Y(\tau, (x_0,\varphi_1)) \\
\left. + \; g_t(x_1) x_1^{-1} \delta \Bigl( \frac{x
+ x_0}{x_1} \Bigr) Y(\tau,(x_0,- \varphi)) \right) 
\end{multline*}
\begin{multline*}
= \varphi_1 \mbox{Res}_{x_1} \mbox{Res}_{x_0} \biggl( 
\biggl(\sum_{n \in \mathbb{Z}} n (x + x_0)^{n-1} \varphi x_1^{-n-1}
\biggr) \biggl( \frac{1}{2} \sum_{j \in \Z} t^j (a^{(m)})^{-2j} A^{(m)}_j
x_1^{j + 1}
\biggr) \biggr.\\
\biggl( \sum_{k \in \mathbb{Z}} G(k + \frac{1}{2})
x_0^{-k - 2} \biggr) 
+ \biggl(\sum_{n \in \mathbb{Z}}(x + x_0)^n x_1^{-n-1} \biggr) 
\biggl( \frac{1}{2} \sum_{j \in \Z} t^j (a^{(m)})^{-2j} A^{(m)}_j  x_1^{j
+ 1} \biggr)\\
\biggl( 2\sum_{k \in \mathbb{Z}} L(k) x_0^{-k - 2} \biggr)
+ \biggl(\sum_{n \in \mathbb{Z}}(x + x_0)^n x_1^{-n-1} \biggr)
\biggl( \sum_{j \in \Z} t^{j-\frac{1}{2}} (a^{(m)})^{-2j +1}
M^{(m)}_{j - \frac{1}{2}} x_1^{j} \biggr) \\
\biggl. \biggl( \sum_{k \in \mathbb{Z}} G(k + \frac{1}{2})
x_0^{-k - 2} - 2\varphi L(k) x_0^{-k - 2} \biggr) \biggr)
\end{multline*}
\begin{multline*}
= \varphi_1 \mbox{Res}_{x_0} \biggl( \sum_{l \in \mathbb{N}} \sum_{k \in
\mathbb{Z}} \sum_{j \in \Z} \biggl( (j + 1) \binom{j}{l} x^{j -l} x_0^l
\frac{\varphi}{2} t^j (a^{(m)})^{-2j} A^{(m)}_j G(k + \frac{1}{2})
x_0^{-k - 2} \\
+ \; \binom{j + 1}{l} x^{j-l + 1} x_0^l t^j (a^{(m)})^{-2j}
A^{(m)}_j L(k) x_0^{-k - 2} \\
+ \; \binom{j}{l} x^{j-l} x_0^l t^{j-\frac{1}{2}} (a^{(m)})^{-2j + 1}
M^{(m)}_{j - \frac{1}{2}} \Bigl(G(k + \frac{1}{2}) - 2\varphi L(k) \Bigr)
x_0^{-k - 2}
\biggr) \biggr) 
\end{multline*}
\begin{multline*}
= \varphi_1 \sum_{l \in \mathbb{N}} \sum_{j \in \Z} \left( \frac{(j +
1)}{2} \binom{j}{l} x^{j -l} \varphi t^j (a^{(m)})^{-2j} A^{(m)}_j
G(l - \frac{1}{2}) \right. \\
+ \; \binom{j + 1}{l} x^{j-l + 1} t^j (a^{(m)})^{-2j} A^{(m)}_j
L(l - 1) \\
\left. + \; \binom{j}{l} x^{j-l} t^{j-\frac{1}{2}} (a^{(m)})^{-2j +1}
M^{(m)}_{j - \frac{1}{2}} (G(l - \frac{1}{2}) - 
2\varphi L(l - 1)) \right) 
\end{multline*}
\begin{multline*}
= \varphi_1 \sum_{k = -1}^\infty \sum_{j \in \Z} \biggl( \binom{j +
1}{k+1} t^j (a^{(m)})^{-2j} x^{j-k} \\
\Bigl(A^{(m)}_j + 2 \left( \frac{j-k}{j + 1} \right)
x^{-1} t^{-\frac{1}{2}} a^{(m)} \varphi M^{(m)}_{j - \frac{1}{2}} \Bigr)
L(k) \\
+ \; \binom{j}{k + 1} t^j (a^{(m)})^{-2j} x^{j -k -1} \Bigl( \frac{(j +
1)}{2}  \varphi  A^{(m)}_j + t^{-\frac{1}{2}} a^{(m)} M^{(m)}_{j -
\frac{1}{2}} \Bigr) G(k + \frac{1}{2})  \biggr) .
\end{multline*}
This finishes the proof of (\ref{special case b bracket}) thus
concluding the proof of the sewing axiom for case (b).

(g)  Let $n \geq 2$ and fix $m \in \Z$.  Assume the sewing axiom holds
for all $Q_1 \in SK(m)$, $Q_2 \in SK(n - 1)$.  We show the sewing
axiom holds for all $Q_1 \in SK(m)$, $Q_2 \in SK(n)$.  Let 
\begin{multline*}
Q_1 = ((z_1, \theta_1),...,(z_{m-1}, \theta_{m-1});
(A^{(0)},M^{(0)}), (a^{(1)}, A^{(1)}, M^{(1)}),...,\\
(a^{(m)}, A^{(m)}, M^{(m)})) , 
\end{multline*}
\begin{multline*}
Q_2 = ((z_1', \theta_1'),...,(z_{n-1}', \theta_{n-1}');
(B^{(0)},N^{(0)}), (b^{(1)}, B^{(1)}, N^{(1)}),...,\\
(b^{(n)}, B^{(n)}, N^{(n)})) .
\end{multline*}
Then since $\nu^Y_{m + n - 1} (Q_1 \; _i\infty_0 \; Q_2)$ and each
term in the series $(\nu^Y_m (Q_1) \; _i*_0 \; \nu^Y_n (Q_2))_{t^{1/2}}$ 
are analytic in $(z_1')_B,...,(z_{n -1}')_B$, we need only prove the
sewing axiom when 
\[|(z_1')_B|,...,|(z_{n - 2}')_B| < |(z_{n - 1}')_B| . \] 
In this case $Q_2 = Q_2' \; _{n - 1}\infty_0 \; Q_2''$ where  
\begin{multline*}
Q_2' = ((z_1', \theta_1'),...,(z_{n-2}', \theta_{n-2}');
(B^{(0)},N^{(0)}), (b^{(1)}, B^{(1)}, N^{(1)}),...,\\
(b^{(n-2)}, B^{(n-2)}, N^{(n-2)}), (1, \mathbf{0}))
\in SK(n-1) 
\end{multline*} 
and 
\[Q_2'' = ((z_{n-1}', \theta_{n-1}'); \mathbf{0}, (b^{(n-1)},
B^{(n-1)}, N^{(n-1)}),(b^{(n)},B^{(n)}, N^{(n)})) \in SK(2) .\]
By our inductive assumption
\[(\nu_m^Y(Q_1) \; _i*_0 \; \nu_{n - 1}^Y(Q_2')) e^{- \Gamma(a^{(i)},
A^{(i)},M^{(i)},B^{(0)},N^{(0)})c} = \nu_{m + n - 2}^Y (Q_1 \;
_i\infty_0 \; Q_2') , \]
and by step (f)
\[\nu_{m + n -2}^Y(Q_1 \; _i\infty_0 \; Q_2') \; _{m + n - 2}*_0 \;
\nu_2^Y(Q_2'') = \nu_{m + n - 1}^Y ((Q_1 \; _i\infty_0 \; Q_2') \; _{m
+ n - 2}\infty_0 \; Q_2'') .\]
Thus using the associativity of sewing and the $t^{1/2}$-contraction,
Propositions \ref{sewing associativity} and \ref{contraction associativity},
and Proposition \ref{doubly convergent}, we have   
\begin{eqnarray*}
\lefteqn{(\nu_m^Y(Q_1) \; _i*_0 \;  \nu_{n}^Y(Q_2)) e^{-
\Gamma(a^{(i)}, A^{(i)},M^{(i)},B^{(0)},N^{(0)})c}} \\
&=& \nu_m^Y(Q_1) \; _i*_0 (\nu_{n-1}^Y(Q_2') \; _{n-1}*_0 \;
\nu_2^Y(Q_2'')) e^{- \Gamma(a^{(i)}, A^{(i)}, M^{(i)}, B^{(0)},
N^{(0)})c} \\ 
&=& ((\nu_m^Y(Q_1) \; _i*_0 \; \nu_{n - 1}^Y(Q_2')) \; _{i + n - 2}*_0
\; \nu_2^Y(Q_2'')) e^{- \Gamma(a^{(i)}, A^{(i)}, M^{(i)}, B^{(0)},
N^{(0)})c} \\  
&=& \nu_{m + n - 1}^Y( (Q_1 \; _i\infty_0 \; Q_2') \; _{i + n -
2}\infty_0 \; Q_2'') \\
&=& \nu_{m + n -1}^Y(Q_1 \; _i\infty_0 \; (Q_2' \; _{n-1}\infty_0 \;
Q_2'')) \\ 
&=& \nu_{m +  n - 1}^Y (Q_1 \; _i\infty_0 \; Q_2)  .
\end{eqnarray*} 

(h)  Finally, assuming the sewing axiom holds for $Q_1 \in SK(m-1)$ with $m
\geq 2$ and $Q_2 \in SK(n)$ for all $n$, we show that it holds for
$Q_1 \in SK(m)$ and $Q_2 \in SK(n)$.  Let  
\begin{multline*}
Q_1 = ((z_1, \theta_1),...,(z_{m-1}, \theta_{m-1});
(A^{(0)},M^{(0)}), (a^{(1)}, A^{(1)}, M^{(1)}),...,\\
(a^{(m)}, A^{(m)}, M^{(m)})) , 
\end{multline*}
\begin{multline*}
Q_2 = ((z_1', \theta_1'),...,(z_{n-1}', \theta_{n-1}');
(B^{(0)},N^{(0)}), (b^{(1)}, B^{(1)}, N^{(1)}),...,\\
(b^{(n)}, B^{(n)}, N^{(n)})) .
\end{multline*}
Then since $\nu^Y_{m + n - 1} (Q_1 \; _i\infty_0 \; Q_2)$ and every
term in the series $(\nu^Y_m (Q_1) \; _i*_0 \; \nu^Y_n
(Q_2))_{t^{1/2}}$ is analytic in $(z_1)_B,...,(z_{m - 1})_B$, we need
only prove the sewing axiom when 
\[|(z_1)_B|,...,|(z_{m - 2})_B| < |(z_{m - 1})_B| . \]  
In this case $Q_1 = Q_1' \; _{m - 1}\infty_0 \; Q_1''$ where
\begin{multline*}
Q_1' = ((z_1, \theta_1),...,(z_{m-2}, \theta_{m-2});
(A^{(0)},M^{(0)}), (a^{(1)}, A^{(1)}, M^{(1)}),...,\\
(a^{(m-2)}, A^{(m-2)}, M^{(m-2)}), (1, \mathbf{0})) \in SK(m-1) 
\end{multline*}  
and 
\[Q_1'' = ((z_{m-1}, \theta_{m-1}); \mathbf{0}, (a^{(m-1)},
A^{(m-1)}, M^{(m-1)}),(a^{(m)},A^{(m)}, M^{(m)})) \in SK(2) .\]
We have three cases:

(i) If $i = m$, then 
\[Q_1 \; _i\infty_0 \; Q_2 = (Q_1' \; _{m - 1}\infty_0 \; Q_1'') \;
_m\infty_0 \; Q_2 = Q_1' \; _{m - 1}\infty_0 \; (Q_1'' \; _2\infty_0
\; Q_2 ) . \]  
For $m > 2$, the sewing axiom follows by the inductive assumption,
associativity of sewing and the $t^{1/2}$-contraction and
Proposition \ref{doubly convergent}.  For $m = 2$ the sewing axiom
follows by special case (c), by (g) (the induction on $n$), the
associativity of sewing and the $t^{1/2}$-contraction and
Proposition \ref{doubly convergent}.     

(ii) If $i = m-1$, then 
\[Q_1 \; _i\infty_0 \; Q_2 = (Q_1' \; _{m - 1}\infty_0 \; Q_1'')  \;
_{m-1}\infty_0 \; Q_2 = Q_1' \; _{m - 1}\infty_0 \; (Q_1'' \;
_1\infty_0 \; Q_2 ) . \] 
For $m > 2$, the sewing axiom follows by the inductive assumption,
associativity of sewing and the $t^{1/2}$-contraction and
Proposition \ref{doubly convergent}.  For $m = 2$ the sewing axiom
follows by special case (e), by (g) (the induction on $n$), the
associativity of sewing and the $t^{1/2}$-contraction and
Proposition \ref{doubly convergent}.    

(iii) If $i \neq m, m - 1$, then 
\[Q_1 \; _i\infty_0 \; Q_2 =  (Q_1' \; _{m - 1}\infty_0 \; Q_1'') \;
_i\infty_0 \; Q_2 = (Q_1' \; _i\infty_0 \; Q_2)  \; _{m + n -
2}\infty_0 \; Q_1'' , \]
and the sewing axiom follows by the inductive assumption, special case
(f), the associativity of sewing and the $t^{1/2}$-contraction
and Proposition \ref{doubly convergent}. 
\end{proof} 

\section{The isomorphism between the category of $N=1$ SG-VOSAs and the
category of $N=1$ NS-VOSAs}\label{iso} 

In this section, we prove the equivalence of the notions of $N=1$ 
NS-VOSA with odd formal variables and $N=1$ SG-VOSA.  By Theorem 5.3 in
\cite{B vosas}, this then implies the equivalence of $N=1$ NS-VOSA without
odd formal variables and $N=1$ SG-VOSAs.  Again, since $N=1$ NS- (resp.,
SG-) VOSAs are isomorphic if taken over different subalgebras of the
underlying Grassmann algebra $\bigwedge_\infty$, we suppress the details
of which subalgebras are actually acting non-trivially.

Let $\mathbf{SG}(c)$ be the category of $N=1$ SG-VOSAs with central
charge $c$.  Let $\mathbf{SV}(\varphi,c)$ be the category of $N=1$
NS-VOSAs with odd formal variables and central charge $c$, and let
$\mathbf{SV}(c)$ be the category of $N=1$ NS-VOSAs without odd formal
variables and with central charge $c$.  Let $1_{SV_\varphi}$ and $1_{SG}$
be the identity functors on $\mathbf{SV}(\varphi,c)$ and $\mathbf{SG}(c)$,
respectively. The following theorem is the main result of this paper.  

\begin{thm}\label{supergeometric and superalgebraic}
Assume that Conjecture \ref{Gamma converges} is true.  
For any $c \in \mathbb{C}$, the two categories $\mathbf{SV}(\varphi,c)$ 
and $\mathbf{SG}(c)$ are isomorphic (and hence $\mathbf{SV}(c)$ 
and $\mathbf{SG}(c)$ are isomorphic).  That is there exist two 
functors $F_{SV}: \mathbf{SV}(\varphi,c) \rightarrow \mathbf{SG}(c)$ and 
$F_{SG}: \mathbf{SG}(c) \rightarrow \mathbf{SV}(\varphi,c)$ such that
$F_{SV} \circ F_{SG} = 1_{SG}$ and $F_{SG} \circ F_{SV} =
1_{SV_\varphi}$. 
\end{thm}

\begin{proof} We begin by defining $F_{SV}$ on objects by $F_{SV}(V,
Y(\cdot, (x, \varphi)), \mathbf{1}, \tau) = (V, \nu^Y)$.  Proposition
\ref{get a sgvosa} shows that $F_{SV}$ takes objects in
$\mathbf{SV}(\varphi,c)$ to objects in $\mathbf{SG}(c)$.   

Let $\gamma: (V_1, Y_1(\cdot, (x, \varphi)), \mathbf{1}_1, \tau_1)
\rightarrow (V_2, Y_2(\cdot, (x, \varphi)), \mathbf{1}_2, \tau_2)$ be an
$N=1$ NS-VOSA homomorphism. Then for $Q \in SK(n)$ given by
\[Q = ((z_1, \theta_1),...,(z_{n-1}, \theta_{n-1}); (A^{(0)},M^{(0)}),
(a^{(1)}, A^{(1)}, M^{(1)}),...,(a^{(n)}, A^{(n)}, M^{(n)})) , \] 
and for $\bar{\gamma} : \bar{V}_1 \longrightarrow \bar{V}_2$ the
unique extension of $\gamma$, we have
\[\langle v', \bar{\gamma} \circ \nu_n^{Y_1}(Q) (v_1 \otimes \cdots
\otimes v_n) \rangle = \hspace{3in} \]
\begin{multline*}
= \iota^{-1}_{1 \cdots n-1} \langle e^{- \sum_{j \in \Z} \left(
A^{(0)}_j L'(j) + M^{(0)}_{j - \frac{1}{2}} G'(j - \frac{1}{2})
\right)} v', \\
\bar{\gamma} ( Y_1(e^{- \sum_{j \in \Z} \left(A^{(1)}_j 
L(j) + M^{(1)}_{j - \frac{1}{2}} G(j - \frac{1}{2}) \right)} \cdot
(a^{(1)})^{-2L(0)} \cdot v_1, (x_1, \varphi_1)) \cdots \\
Y_1(e^{- \sum_{j \in \Z} \left(A^{(n-1)}_j L(j) + M^{(n-1)}_{j
- \frac{1}{2}} G(j - \frac{1}{2}) \right)} \cdot (a^{(n-1)})^{-2L(0)}
\cdot v_{n-1}, \\
(x_{n-1}, \varphi_{n-1})) \cdot e^{- \sum_{j \in \Z} \left(A^{(n)}_j L(j)
+ M^{(n)}_{j - \frac{1}{2}} G(j - \frac{1}{2}) \right)} \cdot \\
\left.
(a^{(n)})^{-2L(0)} \cdot v_n) \rangle \right|_{(x_i,\varphi_i) =
(z_i,\theta_i)} 
\end{multline*}
\begin{multline*}
= \iota^{-1}_{1 \cdots n-1} \langle e^{- \sum_{j \in \Z} \left(
A^{(0)}_j L'(j) + M^{(0)}_{j - \frac{1}{2}} G'(j - \frac{1}{2})
\right)} v', \\
Y_2(e^{- \sum_{j \in \Z} \left(A^{(1)}_j L(j) + M^{(1)}_{j -
\frac{1}{2}} G(j - \frac{1}{2}) \right)} \cdot (a^{(1)})^{-2L(0)}
\cdot \gamma(v_1), (x_1, \varphi_1)) \cdots \\ 
Y_2(e^{- \sum_{j \in \Z} \left(A^{(n-1)}_j L(j) + M^{(n-1)}_{j
- \frac{1}{2}} G(j - \frac{1}{2}) \right)} \cdot (a^{(n-1)})^{-2L(0)}
\cdot \gamma(v_{n-1}), \\
(x_{n-1}, \varphi_{n-1})) \cdot 
e^{- \sum_{j \in \Z} \left(A^{(n)}_j L(j)
+ M^{(n)}_{j - \frac{1}{2}} G(j - \frac{1}{2}) \right)} \cdot \\
\left.
(a^{(n)})^{-2L(0)} \cdot \gamma(v_n) \rangle \right|_{(x_i,\varphi_i)
= (z_i,\theta_i)} 
\end{multline*}
\[ = \langle v', \nu_n^{Y_2}(Q) \circ \gamma^{\otimes n}(v_1 \otimes
\cdots \otimes v_n) \rangle . \hspace{2.5in} \]
Thus defining $F_{SV}(\gamma) = \gamma$, we have that $F_{SV}(\gamma)$ 
is a homomorphism of $N=1$ SG-VOSAs, i.e., $F_{SV}$ takes morphisms in
$\mathbf{SV}(\varphi,c)$ to morphisms in $\mathbf{SG}(c)$.  By
definition, $F_{SV} (\gamma_1 \gamma_2) = \gamma_1 \gamma_2 = F_{SV}
(\gamma_1) F_{SV} (\gamma_2)$, and $F_{SV}  (id_{(V,Y(\cdot, (x,
\varphi)), \mathbf{1}, \tau)}) = id_{(V,\nu^Y)}$. Thus $F_{SV}$ is a
functor.  

We next define $F_{SG}(V,\nu) = (V, Y_{\nu}(\cdot, (x,\varphi)),
\mathbf{1}_{\nu}, \tau_{\nu})$. Proposition \ref{get a vosa} shows that
$F_{SG}$ takes objects in $\mathbf{SG}(c)$ to objects in
$\mathbf{SV}(\varphi, c)$.  Let $\gamma : (V_1, \nu) \rightarrow (V_2,
\mu)$ be an $N=1$ SG-VOSA homomorphism.  Then 
\begin{eqnarray*}
\left. \langle v', \gamma(Y_{\nu} (u, (x,\varphi))v ) \rangle
\right|_{(x,\varphi) = (z,\theta)} &=& \langle v', \bar{\gamma} \circ
\nu_2((z,\theta); \mathbf{0}, (1,\mathbf{0}), (1,\mathbf{0})) (u
\otimes v) \rangle \\
&=& \langle v', \mu_2((z,\theta); \mathbf{0}, (1,\mathbf{0}),
(1,\mathbf{0}))  \gamma^{\otimes 2} (u \otimes v) \rangle \\
&=& \left. \langle v', Y_{\mu} (\gamma(u), (x,\varphi)) \gamma(v)
\rangle \right|_{(x,\varphi) = (z,\theta)} ,
\end{eqnarray*}
\begin{eqnarray*}
\langle v', \gamma(\mathbf{1}_\nu)\rangle &=& \langle v', \bar{\gamma}
\circ \nu_0(\mathbf{0})\rangle \\
&=& \langle v', \mu_0(\mathbf{0})\rangle \\
&=& \langle v', \gamma(\mathbf{1}_\mu)\rangle ,
\end{eqnarray*}
and
\begin{eqnarray*}
\langle v', \gamma(\tau_\nu)\rangle &=& \langle v', \bar{\gamma}
\left(- \frac{\partial}{\partial \epsilon} \nu_0(\mathbf{0},
M(\epsilon, \frac{3}{2})) \right)\rangle \\
&=& - \langle v', \frac{\partial}{\partial \epsilon} \bar{\gamma} \circ
\nu_0(\mathbf{0},M(\epsilon, \frac{3}{2}))\rangle \\
&=& - \langle v', \frac{\partial}{\partial \epsilon}
\mu_0(\mathbf{0},M(\epsilon, \frac{3}{2}))\rangle \\ 
&=& \langle v', \gamma(\tau_\mu)\rangle ,
\end{eqnarray*}
for all $v' \in V'$.  Thus defining $F_{SG}(\gamma) = \gamma$, we have 
that $F_{SG}(\gamma)$ is indeed an $N=1$ NS-VOSA homomorphism,
i.e., $F_{SG}$ takes morphisms in $\mathbf{SG}(c)$ to morphisms in
$\mathbf{SV}(\varphi,c)$.  By definition, $F_{SG} (\gamma_1 \gamma_2) =
\gamma_1 \gamma_2 = F_{SG} (\gamma_1) F_{SG} (\gamma_2)$, and
$F_{SG}(id_{(V,\nu)}) = id_{(V,Y_\nu(\cdot, (x, \varphi)), \mathbf{1}_\nu,
\tau_\nu)})$. Thus $F_{SG}$ is a functor.

The fact that $F_{SV} \circ F_{SG} = 1_{SG}$ and $F_{SG} \circ
F_{SV} = 1_{SV_\varphi}$ on morphisms is trivial.  To show $F_{SV} 
\circ F_{SG} = 1_{SG}$ on objects, we need to show that given an
$N=1$ SG-VOSA, $(V,\nu)$, then $F_{SV} \circ F_{SG} (V, \nu) =
(V,\nu^{Y_\nu}) = (V, \nu)$, i.e.,
\begin{multline*}
\langle v', \nu_n((z_1, \theta_1),...,(z_{n-1}, \theta_{n-1});
(A^{(0)},M^{(0)}), (a^{(1)}, A^{(1)}, M^{(1)}),...,\\
(a^{(n)}, A^{(n)}, M^{(n)})) (v_1 \otimes \cdots
\otimes v_n) \rangle 
\end{multline*}
\begin{multline}\label{functor1}
= \iota^{-1}_{1 \cdots n-1} \langle e^{- \sum_{j \in \Z} \left(
A^{(0)}_j L'(j) + M^{(0)}_{j - \frac{1}{2}} G'(j - \frac{1}{2})
\right)} v',\\
Y_\nu(e^{- \sum_{j \in \Z} \left(A^{(1)}_j 
L(j) + M^{(1)}_{j - \frac{1}{2}} G(j - \frac{1}{2}) \right)} \cdot
(a^{(1)})^{-2L(0)} \cdot v_1, (x_1, \varphi_1)) \cdots \\
Y_\nu(e^{- \sum_{j \in \Z} \left(A^{(n-1)}_j L(j) +
M^{(n-1)}_{j - \frac{1}{2}} G(j - \frac{1}{2}) \right)} \cdot
(a^{(n-1)})^{-2L(0)} \cdot v_{n-1}, (x_{n-1}, \varphi_{n-1})) \cdot \\
e^{- \sum_{j \in \Z} \left(A^{(n)}_j L(j)
+ M^{(n)}_{j - \frac{1}{2}} G(j - \frac{1}{2}) \right)} \cdot \left.
(a^{(n)})^{-2L(0)} \cdot v_n) \rangle \right|_{(x_i,\varphi_i) =
(z_i,\theta_i)} .
\end{multline}
Since both $(V,\nu)$ and $F_{SV} \circ F_{SG} (V, \nu) =
(V,\nu^{Y_\nu})$ are $N=1$ SG-VOSAs, they both satisfy the sewing axiom. 
Thus by Proposition \ref{how to get any supersphere}, we need only prove
(\ref{functor1}) for $Q = (\mathbf{0})$, $((A^{(0)},M^{(0)}),
(a^{(1)}, A^{(1)},M^{(1)}))$, and $((z, \theta);\mathbf{0}, (1,
\mathbf{0}), (1, \mathbf{0}))$.  For $Q = (\mathbf{0})$ and $((z,
\theta);\mathbf{0}, (1, \mathbf{0}), (1, \mathbf{0}))$, equation
(\ref{functor1}) follows immediately {}from the definitions of
$\mathbf{1}_\nu$ and $Y_\nu$.  The equality
\begin{multline*}
\langle v', \nu_1((A^{(0)},M^{(0)}), (a^{(1)}, A^{(1)}, M^{(1)}))
(v) \rangle \\
= \langle e^{- \sum_{j \in \Z} \left(A^{(0)}_j L'(j) +
M^{(0)}_{j - \frac{1}{2}} G'(j - \frac{1}{2}) \right)} v', 
e^{- \sum_{j \in \Z} \left(A^{(1)}_j L(j) + M^{(1)}_{j
- \frac{1}{2}} G(j - \frac{1}{2}) \right)} \cdot \\
(a^{(1)})^{-2L(0)}
\cdot v) \rangle 
\end{multline*}
is equivalent to the following three equations
\begin{eqnarray}
\ \ \ \ \ \ \  \langle v', \nu_1((A^{(0)},M^{(0)}), (1,\mathbf{0}))(v) \rangle 
 \! \! &=& \! \! 
\langle e^{- \sum_{j \in \Z} \left(A^{(0)}_j L'(j) +
M^{(0)}_{j - \frac{1}{2}} G'(j - \frac{1}{2}) \right)} v', v) \rangle 
\label{functor2} \\
\langle v', \nu_1(\mathbf{0}, (1 , A^{(1)}, M^{(1)})) (v) \rangle \! \! &=& \! \!
\langle v', e^{- \sum_{j \in \Z} \left(A^{(1)}_j L(j) + 
M^{(1)}_{j - \frac{1}{2}} G(j - \frac{1}{2}) \right)} \cdot v) \rangle
\quad   \label{functor3} \\
\langle v', \nu_1(\mathbf{0}, (a^{(1)}, \mathbf{0})) (v)
\rangle \! \! &=& \! \! \langle v', (a^{(1)})^{-2L(0)} \cdot v) \rangle .
\label{functor4} 
\end{eqnarray}
The equality (\ref{functor4}) follows immediately {}from the grading
axiom.  By the supermeromorphicity axiom, $\langle v',
\nu_1(\mathbf{0}, (1 , A^{(1)}, M^{(1)})) (v) \rangle$ is a polynomial
in the $A^{(1)}_j$'s and $M^{(1)}_{j - 1/2}$'s for $j \in \Z$.
Thus $\langle v', \nu_1(\mathbf{0}, (1 , A^{(1)}, M^{(1)})) (v)
\rangle$ is well defined for any $(A^{(1)}, M^{(1)}) \in
\bigwedge_\infty^\infty$.  But since the series $(A^{(1)}, M^{(1)})$ is
not necessarily in $\mathcal{H}$, the formal series $\tilde{E}(A^{(1)},
M^{(1)})(x,\varphi)$ is not necessarily convergent.  In
\cite{B memoirs}, such superspheres as $(\mathbf{0}, (1 , A^{(1)},
M^{(1)}))$ were called {\it generalized superspheres with tubes}, and the
space of such formal superspheres was denoted $\tilde{SK}(1)$.  However,
since $\langle v',\nu_1(\mathbf{0}, (1 , t(A^{(1)}, M^{(1)}))) (v) \rangle$
is linear in $v$ and $v'$, there exists a linear operator $R(t)$ {}from $V$
to $\bar{V}$ such that
\begin{equation}\label{functor5}
\langle v',\nu_1(\mathbf{0}, (1 , t(A^{(1)}, M^{(1)})) (v) \rangle =
\langle v', R(t)v \rangle .
\end{equation}

We first show that, in fact, $R(t)$ is an operator {}from $V$ to itself.
Using the sewing axiom and (\ref{functor4}), we have
\begin{eqnarray}
\lefteqn{\langle v', \nu_1((\mathbf{0},(t_1^\frac{1}{2},\mathbf{0})) \;
_1\infty_0 \; (\mathbf{0}, (1, t(A^{(1)},M^{(1)}))) (v) \rangle} \nonumber \\
&=& \langle v', (\nu_1(\mathbf{0},(t_1^\frac{1}{2},\mathbf{0})) \;
_1*_0 \; \nu_1(\mathbf{0}, (1, t(A^{(1)},M^{(1)}))) (v) \rangle  \nonumber \\
&=& \sum_{k \in \frac{1}{2} \mathbb{Z}} \sum_{i^{(k)} = 1}^{\dim V_{(k)}}
\langle v', \nu_1(\mathbf{0},(t_1^\frac{1}{2},\mathbf{0}))
e_{i^{(k)}}^{(k)} \rangle \langle (e_{i^{(k)}}^{(k)})^*,
\nu_1(\mathbf{0}, (1, t(A^{(1)},M^{(1)})) (v) \rangle \nonumber \\
&=& \sum_{k \in \frac{1}{2} \mathbb{Z}} \sum_{i^{(k)} = 1}^{\dim V_{(k)}}
\langle v', t_1^{- L(0)}  e_{i^{(k)}}^{(k)}\rangle \langle
(e_{i^{(k)}}^{(k)})^*, R(t)v \rangle \nonumber \\
&=& \langle v', t_1^{- L(0)} R(t)v \rangle \nonumber\\
\qquad &=& \langle t_1^{- L'(0)} v', R(t)v \rangle . \label{R goes to V1}
\end{eqnarray}

On the other hand, {}from 
\begin{multline*}
(\mathbf{0},(t_1^\frac{1}{2},\mathbf{0})) \; _1\infty_0 \;
(\mathbf{0}, (1, t(A^{(1)},M^{(1)})) \\
= (\mathbf{0}, (1,\{t t_1^{j}A^{(1)}_j, t t_1^{j -
\frac{1}{2}} M^{(1)}_{j -\frac{1}{2}} \}_{j \in \Z} ))  \;
_1\infty_0 \; (\mathbf{0},(t_1^\frac{1}{2},\mathbf{0})) ,
\end{multline*}
we have
\begin{eqnarray}
\lefteqn{ \langle v', \nu_1((\mathbf{0},(t_1^\frac{1}{2},\mathbf{0})) \;
_1\infty_0 \; (\mathbf{0}, (1, t(A^{(1)},M^{(1)}))) (v) \rangle } \nonumber \\
&=&  \langle v', \nu_1((\mathbf{0}, (1,\{t t_1^{j}A^{(1)}_j, t
t_1^{j - \frac{1}{2}} M^{(1)}_{j -\frac{1}{2}} \}_{j \in \Z} ))  \;
_1\infty_0 \; (\mathbf{0},(t_1^\frac{1}{2},\mathbf{0}))) (v) \rangle
\nonumber \\ 
&=&  \langle v', (\nu_1(\mathbf{0}, (1,\{t t_1^{j}A^{(1)}_j, t
t_1^{j - \frac{1}{2}} M^{(1)}_{j -\frac{1}{2}} \}_{j \in \Z} ))  \;
_1*_0 \; \nu_1(\mathbf{0},(t_1^\frac{1}{2},\mathbf{0}))) (v) \rangle
\nonumber \\ 
&=&  \sum_{k \in \frac{1}{2} \mathbb{Z}} \sum_{i^{(k)} = 1}^{\dim
V_{(k)}} \langle v', \nu_1(\mathbf{0}, (1,\{t t_1^{j}A^{(1)}_j, t
t_1^{j - \frac{1}{2}} M^{(1)}_{j -\frac{1}{2}} \}_{j \in \Z} ))
e_{i^{(k)}}^{(k)} \rangle \nonumber \\
& & \hspace{2in} \langle (e_{i^{(k)}}^{(k)})^*,
\nu_1(\mathbf{0},(t_1^\frac{1}{2},\mathbf{0})) (v) \rangle  \nonumber\\
&=&  \sum_{k \in \frac{1}{2} \mathbb{Z}} \sum_{i^{(k)} = 1}^{\dim
V_{(k)}} \langle v', \nu_1(\mathbf{0}, (1,\{t t_1^{j}A^{(1)}_j, t
t_1^{j - \frac{1}{2}} M^{(1)}_{j -\frac{1}{2}} \}_{j \in \Z} ))
e_{i^{(k)}}^{(k)} \rangle \nonumber \\
& & \hspace{2.5in} \langle (e_{i^{(k)}}^{(k)})^*,t_1^{-L(0)}v
\rangle  \nonumber \\
\quad &=& \langle v', \nu_1(\mathbf{0}, (1,\{t t_1^{j}A^{(1)}_j, t t_1^{j -
\frac{1}{2}} M^{(1)}_{j -\frac{1}{2}} \}_{j \in \Z} ))t_1^{-L(0)}v
\rangle  . \label{R goes to V2}
\end{eqnarray}

When $v' \in V'_{(m)}$ and $v \in V_{(n)}$, combining (\ref{R goes to
V1}) and (\ref{R goes to V2}), we have 
\begin{equation}\label{R goes to V3}
\langle v', R(t)v \rangle t_1^{- m}  = \langle v',
\nu_1(\mathbf{0}, (1,\{t t_1^{j}A^{(1)}_j, t t_1^{j - 
\frac{1}{2}} M^{(1)}_{j -\frac{1}{2}} \}_{j \in \Z} ))v
\rangle t_1^{-n} .
\end{equation}
By the supermeromorphicity axiom, $\langle v', \nu_1(\mathbf{0}, (1, \{t
t_1^{j}A^{(1)}_j, t t_1^{j - 1/2} M^{(1)}_{j -1/2} \}_{j \in
\Z} ))v \rangle$ is a power series in $t_1^{1/2}$.  Then {}from (\ref{R
goes to V3}), we conclude that $\langle v', R(t)v \rangle = 0$ when $m > n$.
{}From the positive energy axiom, it follows that $R(t)$ is an operator 
{}from $V$ to itself. 

{}From Proposition \ref{t and s composition}, the sewing axiom and
(\ref{functor5}), we have
\begin{eqnarray*}
\lefteqn{\langle v', R(t_1 + t_2)v \rangle} \\ 
&=& \langle v',\nu_1(\mathbf{0}, (1, (t_1 + t_2) (A^{(1)}, M^{(1)}))) 
(v) \rangle \\
&=& \langle v',\nu_1 \left( (\mathbf{0}, (1, t_1(A^{(1)}, M^{(1)})))
\; _1\infty_0 \;  (\mathbf{0}, (1, t_2(A^{(1)}, M^{(1)}))) \right) (v)
\rangle \\ 
&=& \langle v',\left( \nu_1(\mathbf{0}, (1, t_1(A^{(1)}, M^{(1)})))
\; _1*_0 \; \nu_1 (\mathbf{0}, (1, t_2(A^{(1)}, M^{(1)}))) \right) (v)
\rangle \\ 
&=& \langle v', R(t_1) R(t_2)v \rangle
\end{eqnarray*}
for all $v \in V$ and $v' \in V'$.  Thus
\begin{equation}\label{functor6}
R(t_1 + t_2) = R(t_1) R(t_2) .
\end{equation}
Since $\langle v',\nu_1(\mathbf{0}, (1, t(A^{(1)}, M^{(1)}))) (v)
\rangle$ is a polynomial in $t$, so is $\langle v', R(t)v \rangle$.
Hence $R(t) \in (\mbox{End} \; V)[[t]]$.  This fact and equation 
(\ref{functor6}) imply that there exists $X \in \mbox{End} \; V$
such that $R(t) = e^{tX}$.  Substituting this into equation
(\ref{functor5}), taking the derivative with respect to $t$ and then
letting $t = 0$, we obtain
\[ \frac{\partial}{\partial t} \langle v',\left.\nu_1(\mathbf{0}, (1, 
t(A^{(1)}, M^{(1)}))) (v) \rangle \right|_{t = 0} = \langle v', Xv
\rangle . \]
But {}from the definition of $G(j - 1/2)$ and $L(j)$ for $j \in
\Z$,
\begin{multline*}
\frac{\partial}{\partial t} \langle v', \left. \nu_1(\mathbf{0}, (1, 
t(A^{(1)}, M^{(1)}))) (v) \rangle \right|_{t = 0} \\
= \Bigl. \sum_{j \in \Z} \biggl( - A^{(1)}_j \biggl( -
\frac{\partial}{\partial(tA^{(1)}_j)} \langle v',\nu_1(\mathbf{0}, (1,
t(A^{(1)}, M^{(1)}))) (v) \rangle \Bigr|_{t(A^{(1)}, M^{(1)}) =
\mathbf{0}} \biggr) \biggr.  \\ 
\biggl. \Bigl. - M^{(1)}_{j - \frac{1}{2}} \biggl( -
\frac{\partial}{\partial(tM^{(1)}_{j - \frac{1}{2}})} \langle
v',\nu_1(\mathbf{0}, (1, t(A^{(1)}, M^{(1)}))) (v) \rangle
\Bigr|_{t(A^{(1)}, M^{(1)}) = \mathbf{0}} \biggr) \biggr) .
\end{multline*} 
Thus $X = - \sum_{j \in \Z} \left( - A^{(1)}_j L(j) + M^{(1)}_{j -
1/2} G(j - 1/2) \right)$ and hence
\begin{eqnarray*}
\langle v', \nu_1(\mathbf{0}, (1 , A^{(1)}, M^{(1)})) (v) \rangle &=&
\langle v', e^X v \rangle \\
&=& \langle v', e^{- \sum_{j \in \Z} \left(A^{(1)}_j L(j) + 
M^{(1)}_{j - \frac{1}{2}} G(j - \frac{1}{2}) \right)} \cdot v) \rangle
\end{eqnarray*}
giving (\ref{functor3}).

Similarly, by the supermeromorphicity axiom, $\langle v',
\nu_1((A^{(0)}, M^{(0)}), (1,\mathbf{0})) (v) \rangle$ is a polynomial 
in the $A^{(0)}_j$'s and $M^{(0)}_{j - 1/2}$'s for $j \in \Z$.
Thus $\langle v',\nu_1((A^{(0)}, M^{(0)}), (1,\mathbf{0})) (v)
\rangle$ is well defined for any $(A^{(0)}, M^{(0)}) \in
\bigwedge_\infty^\infty$.  In general $((A^{(0)}, M^{(0)}),
(1,\mathbf{0})) \in \tilde{SK}(1)$, i.e., it is a generalized supersphere
with tubes.  Since $\langle v', \nu_1((A^{(0)}, M^{(0)}), (1,\mathbf{0})) (v)
\rangle$ is linear in $v$ and $v'$, there exists a linear operator $R'(t)$
{}from $V'$ to $\bar{V}^*$ such that 
\begin{equation}\label{functor7}
\langle v',\nu_1(t(A^{(0)}, M^{(0)}), (1,\mathbf{0}))  (v) \rangle =
\langle R'(t)v', v \rangle .
\end{equation}

We first show that, in fact, $R'(t)$ is an operator {}from $V'$ to itself.
Using the sewing axiom and (\ref{functor4}), we have
\begin{eqnarray}
\lefteqn{\langle v', \nu_1((t(A^{(0)}, M^{(0)}), (1,\mathbf{0})) \;
_1\infty_0 \; (\mathbf{0},(t_1^\frac{1}{2},\mathbf{0}))) (v) \rangle} \nonumber\\ 
&=& \langle v', (\nu_1(t(A^{(0)}, M^{(0)}), (1,\mathbf{0})) \; _1*_0 \;
\nu_1(\mathbf{0},(t_1^\frac{1}{2},\mathbf{0}))) (v) \rangle  \nonumber\\
&=& \sum_{k \in \frac{1}{2} \mathbb{Z}} \sum_{i^{(k)} = 1}^{\dim V_{(k)}}
\langle v', \nu_1(t(A^{(0)}, M^{(0)}), (1,\mathbf{0}))  e_{i^{(k)}}^{(k)}
\rangle \langle (e_{i^{(k)}}^{(k)})^*,
\nu_1(\mathbf{0},(t_1^\frac{1}{2},\mathbf{0})) (v) \rangle \nonumber\\  
&=& \sum_{k \in \frac{1}{2} \mathbb{Z}} \sum_{i^{(k)} = 1}^{\dim V_{(k)}}
\langle R'(t) v', e_{i^{(k)}}^{(k)}\rangle \langle
(e_{i^{(k)}}^{(k)})^*, t_1^{- L(0)}  v \rangle \nonumber \\ 
&=& \langle R'(t) v', t_1^{- L(0)} v \rangle \nonumber\\
\qquad &=& \langle R'(t) v', t_1^{- L(0)} v \rangle . \label{R' goes to V1}
\end{eqnarray}

On the other hand, {}from 
\begin{multline*}
(t(A^{(0)}, M^{(0)}),(1,\mathbf{0})) \; _1\infty_0 \;
(\mathbf{0},(t_1^\frac{1}{2},\mathbf{0})) \\
= (\mathbf{0},(t_1^\frac{1}{2},\mathbf{0}))  \;
_1\infty_0 \; ((\{t t_1^j A^{(0)}_j, t t_1^{j - \frac{1}{2}}
M^{(0)}_{j -\frac{1}{2}} \}_{j \in \Z} ),(1,\mathbf{0}))  ,
\end{multline*}
we have
\begin{eqnarray}
\lefteqn{\langle v', \nu_1((t(A^{(0)}, M^{(0)}),(1,\mathbf{0})) \; _1\infty_0
\; (\mathbf{0},(t_1^\frac{1}{2},\mathbf{0}))) (v) \rangle } \nonumber\\
&=& \langle v', \nu_1((\mathbf{0},(t_1^\frac{1}{2},\mathbf{0})) \; 
_1\infty_0 \; ((\{t t_1^j A^{(0)}_j, t t_1^{j - \frac{1}{2}}
M^{(0)}_{j -\frac{1}{2}} \}_{j \in \Z} ),(1,\mathbf{0}))) (v) \rangle
\nonumber\\ 
&=& \langle v', (\nu_1(\mathbf{0},(t_1^\frac{1}{2},\mathbf{0})) \;
_1*_0 \; \nu_1((\{t t_1^j A^{(0)}_j, t t_1^{j - \frac{1}{2}}
M^{(0)}_{j -\frac{1}{2}} \}_{j \in \Z} ),(1,\mathbf{0}))) (v) \rangle
\nonumber\\ 
&=& \sum_{k \in \frac{1}{2} \mathbb{Z}} \sum_{i^{(k)} = 1}^{\dim V_{(k)}}
\langle v', \nu_1(\mathbf{0},(t_1^\frac{1}{2},\mathbf{0})) 
e_{i^{(k)}}^{(k)} \rangle \nonumber\\ 
& & \hspace{1in} \langle (e_{i^{(k)}}^{(k)})^*, 
\nu_1((\{t t_1^j A^{(0)}_j, t t_1^{j - \frac{1}{2}}
M^{(0)}_{j -\frac{1}{2}} \}_{j \in \Z} ),(1,\mathbf{0})) (v) \rangle
\nonumber\\  
&=& \sum_{k \in \frac{1}{2} \mathbb{Z}} \sum_{i^{(k)} = 1}^{\dim V_{(k)}}
\langle v', t_1^{-L(0)} e_{i^{(k)}}^{(k)} \rangle \nonumber  \\
& & \hspace{1in} \langle
(e_{i^{(k)}}^{(k)})^*, \nu_1((\{t t_1^j A^{(0)}_j, t t_1^{j - \frac{1}{2}}
M^{(0)}_{j -\frac{1}{2}} \}_{j \in \Z} ),(1,\mathbf{0})) (v) 
\rangle  \nonumber \\
\quad &=& \langle t_1^{-L'(0)} v', \nu_1((\{t t_1^j A^{(0)}_j, t t_1^{j -
\frac{1}{2}} M^{(0)}_{j -\frac{1}{2}} \}_{j \in \Z} ),(1,\mathbf{0}))
(v) \rangle  . \label{R' goes to V2}
\end{eqnarray}

When $v' \in V'_{(m)}$ and $v \in V_{(n)}$, combining (\ref{R' goes to
V1}) and (\ref{R' goes to V2}), we have 
\begin{equation}\label{R' goes to V3}
\langle R'(t) v', v \rangle t_1^{- n}  = \langle v',
\nu_1((\{t t_1^j A^{(0)}_j, t t_1^{j - \frac{1}{2}} M^{(0)}_{j
-\frac{1}{2}} \}_{j \in \Z} ),(1,\mathbf{0})) (v) \rangle t_1^{-m} . 
\end{equation}
By the supermeromorphicity axiom, $\langle v', \nu_1((\{t t_1^j
A^{(0)}_j, t t_1^{j - 1/2} M^{(0)}_{j - 1/2} \}_{j \in
\Z} ),(1,\mathbf{0})) (v) \rangle$ is a power series in
$t_1^{1/2}$.  Then {}from  (\ref{R' goes to V3}), we conclude
that $\langle R'(t)v', v \rangle = 0$ when $m < n$. {}From the positive
energy axiom, it follows that $R'(t)$ is an operator {}from $V'$ to
itself.  Thus the adjoint $R$ {}from $V$ to $V$ is well defined. {}From
Proposition \ref{t and s composition}, the sewing axiom and
(\ref{functor7}), we have
\begin{eqnarray*}
\lefteqn{ \langle R'(t_1 + t_2) v', v \rangle} \\
&=& \langle v',\nu_1 ((t_1 + t_2) (A^{(0)}, M^{(0)}), (1,\mathbf{0}))
(v) \rangle \\ 
&=& \langle v',\nu_1 \left( ((t_1(A^{(0)}, M^{(0)}), (1,\mathbf{0}))  
\; _1\infty_0 \; ((t_2(A^{(0)}, M^{(0)}), (1,\mathbf{0})) \right) (v)
\rangle \\ 
&=& \langle v',\left( \nu_1 ((t_1(A^{(0)}, M^{(0)}), (1,\mathbf{0})) 
\; _1*_0 \; \nu_1 ((t_2(A^{(0)}, M^{(0)}), (1,\mathbf{0})) \right) (v) 
\rangle \\ 
&=& \sum_{k \in \frac{1}{2} \mathbb{Z}} \sum_{i^{(k)} = 1}^{\dim V_{(k)}}
\langle R'(t_1) v', e_{i^{(k)}}^{(k)} \rangle \langle R'(t_2) 
(e_{i^{(k)}}^{(k)})^*, v \rangle \\
&=&  \sum_{k \in \frac{1}{2} \mathbb{Z}} \sum_{i^{(k)} = 1}^{\dim V_{(k)}}
\langle R'(t_1) v', e_{i^{(k)}}^{(k)} \rangle \langle
(e_{i^{(k)}}^{(k)})^*, R(t_2) v \rangle \\ 
&=&  \langle R'(t_1) v', R(t_2) v \rangle \\ 
&=& \langle R'(t_2) R'(t_1) v', v \rangle
\end{eqnarray*}
for all $v \in V$ and $v' \in V'$.  Thus
\begin{equation}\label{functor8}
R'(t_1 + t_2) = R'(t_1) R'(t_2) .
\end{equation}
Since $\langle v', \nu_1(t(A^{(0)}, M^{(0)}), (1,\mathbf{0})) (v)
\rangle$ is a polynomial in $t$, so is $\langle R'(t) v', v \rangle$. 
Hence $R'(t) \in (\mbox{End} \; V')[[t]]$.  This fact and equation
(\ref{functor8}) imply that there exists $X \in \mbox{End} \; V'$
such that $R'(t) = e^{tX}$.  Substituting this into equation
(\ref{functor7}), taking the derivative with respect to $t$ and then
letting $t = 0$, we obtain
\[\frac{\partial}{\partial t} \langle v',\left. \nu_1 (t(A^{(0)},
M^{(0)}), (1,\mathbf{0}))  (v) \rangle \right|_{t = 0} = \langle Xv',
v \rangle . \] 
But {}from the definition of $G(- j + 1/2)$ and $L(-j)$ for $j
\in \Z$,
\begin{eqnarray*}
\lefteqn{\frac{\partial}{\partial t} \langle v',\left. \nu_1 (t(A^{(0)},
M^{(0)}), (1,\mathbf{0})) (v) \rangle \right|_{t = 0}} \\
&=& \Bigl. \sum_{j \in \Z} \biggl( - A^{(0)}_j \biggl( -
\frac{\partial}{\partial(tA^{(0)}_j)} \langle v',\nu_1 (t(A^{(0)},
M^{(0)}), (1,\mathbf{0})) (v) \rangle \Bigr|_{t(A^{(0)}, M^{(0)}) =
\mathbf{0}} \biggr) \biggr.  \\
& & \biggl. \Bigl. - M^{(0)}_{j - \frac{1}{2}} \biggl( -
\frac{\partial}{\partial(tM^{(0)}_{j - \frac{1}{2}})} \langle
v',\nu_1 (t(A^{(0)}, M^{(0)}), (1,\mathbf{0})) (v) \rangle
\Bigr|_{t(A^{(0)}, M^{(0)}) = \mathbf{0}} \biggr) \biggr) \\
&=& - \sum_{j \in \Z} \langle v', A^{(0)}_j L(-j) v \rangle - \sum_{j
\in \Z} \langle v', M^{(0)}_{j - \frac{1}{2}} G(- j + \frac{1}{2}) v
\rangle . 
\end{eqnarray*}
Thus $X = - \sum_{j \in \Z} \left(A^{(0)}_j L'(j) + M^{(0)}_{j -
1/2} G'(j - 1/2) \right)$, and therefore
\begin{eqnarray*}
\langle v', \nu_1 ((A^{(0)}, M^{(0)}), (1,\mathbf{0})) (v) \rangle &=& 
\langle e^X v', v \rangle \\
&=& \langle e^{- \sum_{j \in \Z} \left(A^{(0)}_j L'(j) + 
M^{(0)}_{j - \frac{1}{2}} G'(j - \frac{1}{2}) \right)} \cdot v', v
\rangle 
\end{eqnarray*}
giving (\ref{functor2}).  Hence $F_{SV} \circ F_{SG} = 1_{SG}$.

Finally, that $F_{SG} \circ F_{SV}  = 1_{SV_\varphi}$ on objects
is equivalent to $(V,Y(\cdot,(x,\varphi)), \mathbf{1}, \tau) =
(V,Y_{\nu^Y}(\cdot,(x,\varphi)), \mathbf{1}_{\nu^Y}, \tau_{\nu^Y})$.
But this is obvious since by definition
\begin{eqnarray*}
\left. \langle v', Y_{\nu^Y}(v_1,(x,\varphi))v_2 \rangle
\right|_{(x,\varphi) = (z,\theta)} &=& \langle v', \nu_2^Y((z,\theta);
\mathbf{0}, (1, \mathbf{0}), (1,\mathbf{0})) (v_1 \otimes v_2) \rangle
\\
&=& \left. \langle v', Y(v_1,(x,\varphi))v_2 \rangle
\right|_{(x,\varphi)  = (z,\theta)} ,
\end{eqnarray*}
\begin{eqnarray*}
\langle v', \mathbf{1}_{\nu^Y}\rangle &=& \langle v',
\nu_0^Y(\mathbf{0})\rangle \\
&=& \langle v', \mathbf{1} \rangle ,
\end{eqnarray*} 
and 
\begin{eqnarray*}
\langle v', \tau_{\nu^Y}\rangle &=& \langle v',
- \frac{\partial}{\partial \epsilon} \nu_0^Y(\mathbf{0}, M(\epsilon,
\frac{3}{2}))\rangle \\
&=& \langle e^{\epsilon G'(\frac{3}{2})} v', \mathbf{1} \rangle \\
&=& \langle G'(\frac{3}{2}) v', \mathbf{1} \rangle \\
&=& \langle v',  G(-\frac{3}{2}) \mathbf{1} \rangle \\
&=& \langle v', \tau \rangle .
\end{eqnarray*} 
Therefore $F_{SG} \circ F_{SV}  = 1_{SV_\varphi}$. 
\end{proof}

\begin{rema}\label{nonsuper for chapter 5} {\em The results of
this paper subsume the correspondence between geometric vertex
operator algebras and vertex operator algebras developed in \cite{H
book}, assuming the convergence of the central series $\Gamma$, i.e.,
assuming Conjecture \ref{Gamma converges}.}
\end{rema}

\end{document}